\documentclass{amsproc}
\usepackage{breqn,float,amssymb,pdflscape,rotating,setspace}
\makeatletter
\@namedef{subjclassname@2020}{%
  \textup{2020} Mathematics Subject Classification}
\makeatother
\newtheorem{theorem}{Theorem}[section]
\newtheorem{lemma}[theorem]{Lemma}
\newtheorem{proposition}[theorem]{Proposition}

\theoremstyle{definition}

\newtheorem{example}[theorem]{Example}

\theoremstyle{remark}

\numberwithin{equation}{section}
\allowdisplaybreaks
\begin{document}

\title{A Table of Generating Functions}

%    Remove any unused author tags.

%    author one information
\author{Robert Reynolds}
\address[Robert Reynolds]{Department of Mathematics and Statistics, York University, Toronto, ON, Canada, M3J1P3}
\email[Corresponding author]{milver73@gmail.com}
\thanks{}

%    author two information
%\author{ Allan Stauffer}
%\address[Allan Stauffer]{Department of Mathematics and Statistics, York University, Toronto, ON, Canada, M3J1P3}
%\email{stauffer@yorku.ca}
%\thanks{This research is supported by NSERC Canada under Grant 504070}



\date{}

\dedicatory{}

\begin{abstract}
This is a compendium of generating functions involving single, double sums and definite integrals. These generating functions also involve special functions in both the summand function and closed form solution.
\end{abstract}
\maketitle
\section{Preliminaries}
\subsection{Lagrange Series}
\begin{theorem}
The Lagrange series has two basic forms;
\begin{equation}\label{eq:lag1}
f(x)=f(0)+\sum_{n=1}^{\infty}\frac{y^n}{n!}\left[D^{n-1}\left(f'(x) \right)\phi^{n}(x) \right]_{x=0}
\end{equation}
and 
\begin{equation}\label{eq:lag2}
\frac{f(x)}{1-y \phi'(x)}=\sum_{n=0}^{\infty}\frac{y^n}{n!}\left[D^{n}\left(f(x)\right) \phi^{n}(x)\right]_{x=0},
\end{equation}
where $y=x/\phi'(x), D=d/dx$, and the primes denote derivatives.
\end{theorem}
In this compendium we apply Theorem (1.1) given by equation (18) in \cite{riordan} to derive the propositions that follow. We apply the method in section (7) in \cite{reyn_plos} to the propositions to derive the corresponding examples. The functions used in this work are; the hypergeometric function $\, _2F_1(a,b ;c;z)$ listed in section (15.2) in \cite{dlmf}, the Euler gamma function $\Gamma(z)$ listed in section (5) in \cite{dlmf}, the Kummer confluent hypergeometric function $\, _1F_1(a;b;z)$ listed in section (13) in \cite{dlmf}, the incomplete gamma function $\Gamma(a,z)$ listed in section (8) in \cite{dlmf}, the Appell hypergeometric function of two variables $F_1(a;b_{1},b_{2} ;c;x,y)$ listed in section (16) in \cite{dlmf}, the generalized hypergeometric function $\,_3F_2\left(\{a,b,c\};\{a_{1},b_{1}\};z\right)$ listed in section (16.2) in \cite{dlmf}, the Gegenbauer polynomial $C_n^{m}(x)$ listed in section (18) in \cite{dlmf}, the generalized hypergeometric function $\, _2F_2\left(\{a,b\};\{c,d\};e\right)$ listed in section (17) in \cite{dlmf}, the generalized hypergeometric function $\, _1F_2\left(a;\{b,c\};z\right)$ listed in section (16) in \cite{dlmf}, the Hermite polynomial $H_k(x)$ listed in section (18) in \cite{dlmf}, the error function $\text{erf}(x)$ listed in section (7.2(i)) in \cite{dlmf}, the imaginary error function $\text{erfi}(x)$ listed in section (7.2) in \cite{dlmf}. All variables are valid over the complex plane unless specified otherwise.
\section{Introduction}
Generating functions are useful in various branches of mathematics, including combinatorics, number theory, and analysis. They offer a versatile and effective way to encode and manipulate sequences, leading to valuable insights and solutions to complex problems. One of the primary applications of generating functions is in counting and combinatorics. By representing combinatorial structures as generating functions, we can determine the number of objects or arrangements with specific properties by examining the coefficients of the generating function. Generating functions also come in handy when dealing with recurrence relations. By transforming a recurrence relation into an equation involving generating functions, we can manipulate the generating function algebraically to find a closed-form expression for the sequence. This approach simplifies the solution of recursive problems. In the area of partition theory, generating functions have proven their worth. They provide efficient computations of partition functions and offer insights into the properties of integer partitions. By analyzing the generating functions, we can gain a deeper understanding of the behaviour of partitions.\\\\
Generating functions are also valuable for asymptotic analysis. By studying the singularities of generating functions in the complex plane, we can determine growth rates, extract coefficients, and derive asymptotic estimates. This analysis helps us understand the long-term behaviour of sequences. Probability theory benefits from generating functions as well. Moment generating functions and probability generating functions encode information about the distribution of random variables. These functions offer insights into the statistical properties and moments of the variables.  Readers interested in literature on generating functions can peruse the following works \cite{buhring,kalla,hansen,lint,riordan,dlmf,mourad,warnke,gessel}.\\\\  
The development of standard tables for binomial coefficient summations and the use of the Lagrange inversion theorem and contour integration contribute to the expansion of generating function techniques and their application in solving functional equations involving special functions.
\section{Generating functions involving the Hypergeometric functions}
\subsection{Entry 1}
\begin{proposition}
%For all $|Re(\alpha)|<1, |Re(z)|<1,\beta\in\mathbb{C}$ then,
\begin{multline}\label{eq:gf1}
\sum_{n=0}^{\infty}\frac{\Gamma (n \beta +1) \left(x (x+1)^{-\beta }\right)^n \,
   _2F_1(-n,-\alpha ;n (\beta -1)+1;z)}{n! \Gamma (\beta  n-n+1)}\\
   =\frac{(x+1) (x
   z+1)^{\alpha }}{-\beta  x+x+1}
\end{multline}
\end{proposition}
where $f(x)=(x+1) (x z+1)^{\alpha }$ and $\phi(x)=(x+1)^{\beta }$.
\begin{example}
\begin{multline}\label{eq:gf2}
\sum_{n=0}^{\infty}\frac{\Gamma (n \beta +1) \left(x (x+1)^{-\beta }\right)^n \,
   _2F_1(-n,-\alpha ;n (\beta -1)+1;z)}{\Gamma (n+2) \Gamma (\beta 
   n-n+1)}\\
   =\frac{(x+1) (x z+1)^{\alpha +1} \, _2F_1\left(1,\alpha -\beta
   +2;\alpha +2;\frac{x z+1}{1-z}\right)}{(\alpha +1) x (z-1)}\\
   -\frac{(x+1)^{\beta
   } \, _2F_1\left(1,\alpha -\beta +2;\alpha +2;\frac{1}{1-z}\right)}{(\alpha +1)
   x (z-1)}
\end{multline}
\end{example}
\begin{example}
\begin{multline}\label{eq:gf3}
\sum_{n=0}^{\infty}\frac{\Gamma (n \beta +1) \left(x (x+1)^{-\beta }\right)^n \, _2F_1(-n,-\alpha ;n (\beta -1)+1;z)}{(n+2)! \Gamma (\beta  n-n+1)}\\
=\int_{0}^{x}\frac{\Gamma (\alpha +1) (s+1)^{-2 \beta } ((\beta
   -1) s-1) (x+1)^{2 \beta } }{x^2 (z-1)}\\
   \left(\frac{(s+1)^{\beta -1} \, _2F_1\left(1,\alpha -\beta +2;\alpha +2;\frac{1}{1-z}\right)}{\Gamma (\alpha +2)}\right. \\ \left.
   -\frac{(s z+1)^{\alpha +1} \,
   _2F_1\left(1,\alpha -\beta +2;\alpha +2;\frac{s z+1}{1-z}\right)}{\Gamma (\alpha +2)}\right)ds
\end{multline}
\end{example}
\subsection{Entry 2}
\begin{proposition}
\begin{multline}\label{eq:gf4}
\sum_{n=0}^{\infty}\frac{\Gamma (a+1) \left(x (c x+1)^{-\beta }\right)^n \, _2F_1(-n,-n \beta
   ;a-n+1;c)}{n! \Gamma (a-n+1)}=-\frac{(x+1)^a (c x+1)}{\beta  c x-c
   x-1}
\end{multline}
\end{proposition}
where $f(x)=(x+1)^a (c x+1)$ and $\phi(x)=(c x+1)^{\beta }$.
\begin{example}
\begin{multline}\label{eq:gf5}
\sum_{n=0}^{\infty}\frac{\Gamma (a+1) \left(x (x z+1)^{-\beta }\right)^n \, _2F_1(-n,-n \beta
   ;a-n+1;z)}{\Gamma (n+2) \Gamma (a-n+1)}\\
   =\frac{(x z+1)^{\beta } \,
   _2F_1\left(1,a-\beta +2;a+2;\frac{z}{z-1}\right)}{(a+1) x
   (z-1)}\\
   -\frac{(x+1)^{a+1} (x z+1) \, _2F_1\left(1,a-\beta +2;a+2;\frac{(x+1)
   z}{z-1}\right)}{(a+1) x (z-1)}
\end{multline}
\end{example}
\subsection{Entry 3}
\begin{proposition}
\begin{equation}\label{eq:gf6}
\sum_{n=0}^{\infty}\frac{\Gamma (n \beta +1) \left(x (x+1)^{-\beta }\right)^n \, _1F_1(-n;\beta  n-n+1;-z)}{n! \Gamma (\beta 
   n-n+1)}=-\frac{(x+1) e^{x z}}{\beta  x-x-1}
\end{equation}
\end{proposition}
where $f(x)=e^{x z}(x+1)$ and $\phi(x)=(x+1)^{\beta }$.
\begin{example}
\begin{multline}\label{eq:gf7}
\sum_{n=0}^{\infty}\frac{\Gamma (n \beta +1) \left(x (x+1)^{-\beta }\right)^n \, _1F_1(-n;n (\beta -1)+1;-z)}{\Gamma (n+2) \Gamma
   (\beta  n-n+1)}\\
   =\frac{e^{-z} (x+1)^{\beta } (-z)^{\beta -1} \Gamma (1-\beta ,-z)}{x}\\
   -e^{-z} (-((x+1) z))^{\beta -1}
   \Gamma (1-\beta ,-((x+1) z))\\
   -\frac{e^{-z} (-((x+1) z))^{\beta -1} \Gamma (1-\beta ,-((x+1) z))}{x}
\end{multline}
\end{example}
\subsection{Entry 4}
\begin{proposition}
\begin{multline}\label{eq:gf8}
\sum_{n=0}^{\infty}\frac{\Gamma (m+n+1) \left(x (x+1)^{-\beta }\right)^n \Gamma ((m+n) \beta +1) }{\Gamma (n+1) (m+n)! \Gamma (-n+(m+n) \beta +1)}\\
\, _2F_1(-n,-\alpha ;m+(m+n)
   (\beta -1)+1;z)\\
   =\frac{(x+1)^{\beta  m+1} (x z+1)^{\alpha
   }}{(1-\beta ) x+1}
\end{multline}
\end{proposition}
where $f(x)=(x+1) x^m (x z+1)^{\alpha }$ and $\phi(x)=(x+1)^{\beta }$.
\begin{example}
\begin{multline}\label{eq:gf9}
\sum_{n=0}^{\infty}\frac{\Gamma (m+n+1) \left(x (x+1)^{-\beta }\right)^n \Gamma ((m+n) \beta +1) }{\Gamma (n+1) (m+n+1)! \Gamma (-n+(m+n) \beta +1)}\\
\, _2F_1(-n,-\alpha ;m+(m+n)
   (\beta -1)+1;z)\\
   =\frac{(x+1)^{\beta +\beta  m} F_1(m+1;\beta
   ,-\alpha ;m+2;-x,-x z)}{m+1}
\end{multline}
\end{example}
\subsection{Entry 5.}
\begin{proposition}
\begin{multline}\label{eq:gf10}
\sum_{n=0}^{\infty}\frac{(-1)^n \left(x (1-x)^{-\beta }\right)^n \Gamma ((m+n) \beta +1) \, _1F_1(-n;m+(m+n) (\beta
   -1)+1;z)}{\Gamma (n+1) \Gamma (-n+(m+n) \beta +1)}\\
   =\frac{(-1)^{-m} (x-1) e^{x z} (1-x)^{\beta  m}}{(\beta -1)
   x+1}
\end{multline}
\end{proposition}
where $f(x)=(x-1) x^m e^{x z}$ and $\phi(x)=(1-x)^{\beta }$.
\begin{example}
\begin{multline}\label{eq:gf11}
\sum_{n=0}^{\infty}\frac{(-1)^n \left(x (1-x)^{-\beta }\right)^n \Gamma ((m+n) \beta +1) \, _1F_1(-n;m+(m+n) (\beta
   -1)+1;z)}{(m+n+1) \Gamma (n+1) \Gamma (-n+(m+n) \beta +1)}\\
   =-\int_{0}^{x}\frac{(-1)^{-m} s^m (1-s)^{-\beta } e^{s z}
   (1-x)^{\beta } \left(x (1-x)^{-\beta }\right)^{-m}}{x}ds
\end{multline}
\end{example}
\subsection{Entry 6.}
\begin{proposition}
\begin{equation}\label{eq:gf12}
\sum_{n=0}^{\infty}\frac{\Gamma (m+1) \left(x e^{-b x}\right)^n \, _1F_1(-n;m-n+1;-a-b n)}{n! \Gamma (m-n+1)}=\frac{e^{a x}
   (x+1)^m}{1-b x}
\end{equation}
\end{proposition}
where $f(x)=(x+1)^m e^{x (a+b)}$ and $\phi(x)=1$.
\begin{example}
\begin{multline}\label{eq:gf13}
\sum_{n=0}^{\infty}\frac{\Gamma (m+1) \left(x e^{-b x}\right)^n \, _1F_1(-n;m-n+1;-a-b n)}{\Gamma (n+2) \Gamma
   (m-n+1)}\\
   =\frac{(b-a)^{-m-1} e^{-a+b x+b} \Gamma (m+1,b-a)}{x}\\
   -\frac{(x+1)^{m+1} e^{-a+b x+b} (-((x+1) (a-b)))^{-m-1}
   \Gamma (m+1,-((a-b) (x+1)))}{x}
\end{multline}
\end{example}
\subsection{Entry 7.}
\begin{proposition}
%For all $\alpha,\beta,x,z\in\mathbb{C},m\in\mathbb{Z^{+}}$ then,
\begin{multline}\label{eq:gf14}
\sum_{n=0}^{\infty}\frac{(-1)^n \Gamma (\alpha +1) \left(x (x z+1)^{-\beta }\right)^n \, _2F_1(-n,-((m+n) \beta );-n+\alpha
   +1;-z)}{n! \Gamma (-n+\alpha +1)}\\
   =\frac{(-1)^{-m} (1-x)^{\alpha } (x z+1)^{\beta  m+1}}{\beta  x z-x z-1}
\end{multline}
\end{proposition}
where $f(x)=x^m (1-x)^{\alpha }$ and $\phi(x)=(x z+1)^{\beta }$.
\begin{example}
\begin{multline}\label{eq:gf15}
\sum_{n=0}^{\infty}\frac{(-1)^{n+1} \Gamma (\alpha +1) \left(x (x z+1)^{-\beta }\right)^n \, _2F_1(-n,-((m+n) \beta );-n+\alpha
   +1;-z)}{(n+1)! \Gamma (-n+\alpha +1)}\\
   =(-1)^{-m} \left(x (x z+1)^{-\beta }\right)^{-m} \left(\frac{x^{m-1} (x
   z+1)^{\beta -\beta  m} \, _2F_1\left(1,\alpha +(m-1) \beta +2;\alpha +2;\frac{z}{z+1}\right)}{(\alpha +1)
   (z+1)}\right. \\ \left.
   +\frac{(x-1) x^{m-1} (1-x)^{\alpha } (x z+1)^{\beta +\beta  (m-1)-\beta  m+1} \, _2F_1\left(1,\alpha +(m-1)
   \beta +2;\alpha +2;\frac{z-x z}{z+1}\right)}{(\alpha +1) (z+1)}\right)
\end{multline}
\end{example}
\subsection{Entry 9.}
\begin{proposition}
\begin{multline}\label{eq:gf16}
\sum_{n=0}^{\infty}\frac{\Gamma (m \beta +n \beta +1) \left(x (x+1)^{-\beta } e^{-x z}\right)^n \, _1F_1(-n;n (\beta -1)+m \beta
   +1;(-m-n) z)}{\Gamma (n+1) \Gamma (\beta  n-n+m \beta +1)}\\
   =-\frac{(x+1) \left((x+1)^{\beta } e^{x z}\right)^m}{x
   (\beta +x z+z-1)-1}
\end{multline}
\end{proposition}
where $f(x)=x^m $ and $\phi(x)=(x+1)^{\beta } e^{x z}$.
\begin{example}
\begin{multline}\label{eq:gf17}
\sum_{n=0}^{\infty}\frac{\Gamma (m \beta +n \beta +1) \, _1F_1(-n;n (\beta -1)+m
   \beta +1;-((m+n) z)) \left(x (x+1)^{-\beta } e^{-x
   z}\right)^{m+n}}{\Gamma (n+2) \Gamma (\beta  n-n+m \beta
   +1)}\\
=\left(\frac{(x+1)^{\beta }}{x}\right)^{1-m} e^{-((m-1) (x+1) z)}
   (z-m z)^{\beta +\beta  (-m)-1} \Gamma (m \beta -\beta +1,z-m
   z)\\
-(x+1)^{\beta  (m-1)+1} \left(\frac{(x+1)^{\beta }}{x}\right)^{1-m}
   e^{-((m-1) (x+1) z)} (-((m-1) (x+1) z))^{\beta +\beta  (-m)-1}\\
 \Gamma
   (m \beta -\beta +1,-((m-1) (x+1) z))
\end{multline}
\end{example}

%
%\begin{example}
%\begin{multline}\label{eq:gf18}
%\sum_{n=0}^{\infty}\frac{\Gamma (m \beta +n \beta +1) \left(x (x+1)^{-\beta } e^{-x z}\right)^n \, _1F_1(-n;n (\beta -1)+m \beta
%   +1;(-m-n) z)}{\Gamma (n+1) \Gamma (\beta  n-n+m \beta +1)}\\
%   =-\frac{(x+1) \left((x+1)^{\beta } e^{x z}\right)^m}{x
%   (\beta +x z+z-1)-1}
%\end{multline}
%\end{example}
%%
\subsection{Entry 10.}
\begin{proposition}
%For all $\alpha,\beta,x,z\in\mathbb{C},m\in\mathbb{Z^{+}}$ then,
\begin{multline}\label{eq:gf18}
\sum_{n=0}^{\infty}\frac{(-1)^{n+1} \left(x (1-x)^{\beta }\right)^n \Gamma (1-(m+n) \beta ) \, _2F_1(-n,-\alpha ;m-(m+n) (\beta
   +1)+1;z)}{\Gamma (n+1) \Gamma (-n-(m+n) \beta +1)}\\
   =\frac{(-1)^{-m} (x-1) x^m \left(x (1-x)^{\beta }\right)^{-m}
   (1-x z)^{\alpha }}{\beta  x+x-1}
\end{multline}
\end{proposition}
where $f(x)=x^m (1-x z)^{\alpha }$ and $\phi(x)=(1-x)^{-\beta }$.
\begin{example}
\begin{multline}\label{eq:gf19}
\sum_{n=0}^{\infty}\frac{(-1)^{n+1} \left(x (1-x)^{\beta }\right)^n \Gamma (1-(m+n) \beta ) \, _2F_1(-n,-\alpha ;m-(m+n) (\beta
   +1)+1;z)}{\Gamma (n+2) \Gamma (-n-(m+n) \beta +1)}\\
   =\frac{(-1)^{-m} (1-x)^{-\beta } }{(\alpha +1) x (z-1)}\left(\, _2F_1\left(1,\alpha -m
   \beta +\beta +2;\alpha +2;\frac{1}{1-z}\right)\right. \\ \left.
   -(1-x)^{\beta +\beta  (-m)+1} (1-x z)^{\alpha +1} \,
   _2F_1\left(1,\alpha -m \beta +\beta +2;\alpha +2;\frac{x z-1}{z-1}\right)\right)
\end{multline}
\end{example}
\subsection{Entry 11.}
\begin{proposition}
\begin{multline}\label{eq:gf20}
\sum_{n=0}^{\infty}\frac{\Gamma (-m-n+\alpha +1) \left(x (x+1) (x z+1)^{-\beta }\right)^n }{\Gamma (n+1) \Gamma (m-2 (m+n)+\alpha +1)}\, _2F_1(-n,1-(m+n) \beta ;m-2
   (m+n)+\alpha +1;z)\\
   =-\frac{x^m (x+1)^{\alpha +1} \left(x (x+1) (x
   z+1)^{-\beta }\right)^{-m}}{\beta  x (x+1) z-(2 x+1) (x z+1)}
\end{multline}
\end{proposition}
where $f(x)=\frac{x^m (x+1)^{\alpha }}{x z+1}$ and $\phi(x)=\frac{(x z+1)^{\beta }}{x+1}$.
\subsection{Entry 12.}
\begin{proposition}
\begin{multline}\label{eq:gf21}
\sum_{n=0}^{\infty}\frac{\Gamma (m+n+\alpha +1) \left(\frac{x (x z+1)^{-\beta }}{x+1}\right)^n \, _2F_1(-n,-((m+n) \beta
   );m+\alpha +1;z)}{\Gamma (n+1) \Gamma (m+\alpha +1)}\\
   =-\frac{x^m (x+1)^{\alpha +1} (x z+1) \left(\frac{x (x
   z+1)^{-\beta }}{x+1}\right)^{-m}}{\beta  x^2 z+\beta  x z-x z-1}
\end{multline}
\end{proposition}
where $f(x)=x^m (x+1)^{\alpha }$ and $\phi(x)=(x+1) (x z+1)^{\beta }$.
\subsection{Entry 13.}
\begin{proposition}
\begin{multline}\label{eq:gf22}
\sum_{n=0}^{\infty}\frac{(m+n) (m+n+1)_m \Gamma ((2 m+n) \alpha +1) \left(x (x+1)^{-\alpha }\right)^{2 m+n} }{(2 m+n)! \Gamma (-m-n+(2 m+n) \alpha +2)}\\
\,
   _3F_2(1,1,-m-n+1;2,2 \alpha  m-m-n+n \alpha +2;z)\\
   =\prod_{p=1}^{\infty}\frac{2 (x+1)
   x^{m+1}}{((1-\alpha ) x+1) \left((x z+1)^{2^{-p}}+1\right)}\\
   =\frac{(x+1) x^m \log (x z+1)}{z ((1-\alpha )
   x+1)}
\end{multline}
\end{proposition}
where $f(x)=x^m \log (x z+1)$ and $\phi(x)=(x+1)^{\alpha }$ and the product is recorded in equation (2.6) in \cite{ruffa}.
\begin{example}
\begin{multline}\label{eq:gf23}
\sum_{n=0}^{\infty}\frac{\left(x (x+1)^{-\alpha }\right)^n \Gamma ((m+n) \alpha +1) \, _3F_2(1,1,1-n;2,\alpha  n-n+m \alpha
   +2;z)}{(m+n+1) \Gamma (n) \Gamma (-n+(m+n) \alpha +2)}\\
   =\int_{0}^{x}\frac{s^m (s+1)^{-\alpha } (x+1)^{\alpha } \left(x
   (x+1)^{-\alpha }\right)^{-m} \log (s z+1)}{x z}ds
\end{multline}
\end{example}
\subsection{Entry 14.}
\begin{proposition}
\begin{multline}\label{eq:gf24}
\sum_{n=0}^{\infty}\left(\frac{x}{x+1}\right)^n \binom{m+n+\alpha }{n-1} \, _3F_2(1,1,1-n;2,m+\alpha +2;z)\\
=\frac{1}{z}x^m
   \left(\frac{x}{x+1}\right)^{-m} (x+1)^{\alpha +1} \log (x z+1)
\end{multline}
\end{proposition}
where $f(x)=x^m (x+1)^{\alpha } \log (x z+1)$ and $\phi(x)=x+1$.
\begin{example}
\begin{multline}\label{eq:gf25}
\sum_{n=0}^{\infty}\frac{\left(\frac{x}{x+1}\right)^n \binom{m+n+\alpha }{n-1} \, _3F_2(1,1,1-n;2,m+\alpha +2;z)}{m+n+1}\\
=\frac{s^m
   \left(\frac{x}{x+1}\right)^{-m-1} (s+1)^{\alpha -1} \log (s z+1)}{z}
\end{multline}
\end{example}
\subsection{Entry 15.}
\begin{proposition}
%For all $\alpha ,\beta ,\gamma ,\delta ,x,z\in\mathbb{C}, m\in\mathbb{Z^{+}}$ then, 
\begin{multline}\label{eq:gf26}
\sum_{n=0}^{\infty}\frac{(-1)^{n+1} \Gamma ((m+n) \gamma +\delta +1) \left(x (1-x)^{-\gamma } (x z+1)^{-\alpha }\right)^n}{\Gamma (n+1) \Gamma (-n+(m+n) \gamma +\delta
   +1)}\\
    \,
   _2F_1(-n,-((m+n) \alpha )-\beta ;m+(m+n) (\gamma -1)+\delta +1;-z)\\
   =\frac{(-1)^{-m} x^m (1-x)^{\delta +1} (x z+1)^{\beta +1} \left(x (1-x)^{-\gamma } (x z+1)^{-\alpha
   }\right)^{-m}}{x (\gamma +(\alpha -1) (x-1) z+\gamma  x z-1)+1}
\end{multline}
\end{proposition}
where $f(x)=x^m (1-x)^{\delta } (x z+1)^{\beta }$ and $\phi(x)=(1-x)^{\gamma } (x z+1)^{\alpha }$.
\begin{example}
\begin{multline}\label{eq:gf27}
\sum_{n=0}^{\infty}\frac{(-1)^{n+1} \Gamma ((m+n) \gamma +\delta +1) \left(x (1-x)^{-\gamma } (x z+1)^{-\alpha }\right)^n }{(m+n+1) \Gamma (n+1) \Gamma (-n+(m+n) \gamma
   +\delta +1)}\\
   \,
   _2F_1(-n,-((m+n) \alpha )-\beta ;m+(m+n) (\gamma -1)+\delta +1;-z)\\
   =\frac{(-1)^{-m} x^m (1-x)^{\gamma } (x z+1)^{\alpha } \left(x (1-x)^{-\gamma } (x z+1)^{-\alpha
   }\right)^{-m} }{m+1}\\
   F_1(m+1;\gamma -\delta ,\alpha -\beta ;m+2;x,-x z)
\end{multline}
\end{example}
\subsection{Entry 16. The Gegenbauer polynomial}
\begin{proposition}
\begin{multline}\label{eq:gf28}
\sum_{n=0}^{\infty}\left(x \left(x^2-2 \beta  x+1\right)^{-\alpha }\right)^n C_n^{(-((m+n) \alpha ))}(\beta )\\
=-\frac{x^m
   \left(x^2-2 \beta  x+1\right) \left(x \left(x^2-2 \beta  x+1\right)^{-\alpha }\right)^{-m}}{2 \alpha  x^2-x^2-2
   \alpha  \beta  x+2 \beta  x-1}
\end{multline}
\end{proposition}
where $f(x)=x^m $ and $\phi(x)=\left(x^2-2 \beta  x+1\right)^{\alpha }$.
\begin{example}
\begin{multline}\label{eq:gf29}
\sum_{n=0}^{\infty}\frac{\left(x \left(x^2-2 \beta  x+1\right)^{-\alpha }\right)^n C_n^{(-m-n) \alpha }(\beta
   )}{m+n+1}\\
   =\frac{\left(\left(\left(\sqrt{\beta ^2-1}-\beta \right) x+1\right) \left(1-\left(\sqrt{\beta ^2-1}+\beta
   \right) x\right)\right)^{\alpha } \left(x^2-2 \beta  x+1\right)^{\alpha  m} }{m+1}\\
   F_1\left(m+1;\alpha ,\alpha ;m+2;x
   \left(\beta +\sqrt{\beta ^2-1}\right),\frac{x}{\beta +\sqrt{\beta ^2-1}}\right)
\end{multline}
\end{example}
\subsection{Entry 17.}
\begin{proposition}
\begin{multline}\label{eq:gf30}
\sum_{n=0}^{\infty}\left(x \left(x^2-2 x z+1\right)^{-\alpha }\right)^n C_n^{(-((m+n) \alpha )-\beta )}(z)\\
=-\frac{x^m \left(x^2-2
   x z+1\right)^{\beta +1} \left(x \left(x^2-2 x z+1\right)^{-\alpha }\right)^{-m}}{2 \alpha  x^2-x^2-2 \alpha  x z+2
   x z-1}
\end{multline}
\end{proposition}
where $f(x)=x^m \left(x^2-2 x z+1\right)^{\beta }$ and $\phi(x)=\left(x^2-2 x z+1\right)^{\alpha }$.
\begin{example}
\begin{multline}\label{eq:gf31}
\sum_{n=0}^{\infty}\frac{\left(x \left(x^2-2 x z+1\right)^{-\alpha }\right)^n C_n^{-(m+n) \alpha -\beta }(z)}{m+n+1}\\
=\frac{x^m
   \left(x^2-2 x z+1\right)^{\beta } \left(x \left(x^2-2 x z+1\right)^{-\alpha }\right)^{-m}
   \left(\frac{x+\sqrt{z^2-1}-z}{\sqrt{z^2-1}-z}\right)^{\alpha -\beta }
   \left(\frac{-x+\sqrt{z^2-1}+z}{\sqrt{z^2-1}+z}\right)^{\alpha -\beta } }{m+1}\\
   F_1\left(m+1;\alpha -\beta ,\alpha -\beta
   ;m+2;\frac{x}{z-\sqrt{z^2-1}},\frac{x}{z+\sqrt{z^2-1}}\right)
\end{multline}
\end{example}
\subsection{Entry 18.}
\begin{proposition}
\begin{equation}\label{eq:gf32}
\sum_{n=0}^{\infty}\frac{\left((-1)^n+1\right) \left(x e^{-b x^2}\right)^n (a+b (m+n))^{n/2}}{2 \Gamma
   \left(\frac{n+2}{2}\right)}=\frac{e^{a x^2} x^m \left(x e^{-b x^2}\right)^{-m}}{1-2 b x^2}
\end{equation}
\end{proposition}
where $f(x)=e^{a x^2} x^m$ and $\phi(x)=e^{b x^2}$.
%
%\begin{proposition}
%\begin{multline}\label{eq:gf33}
%\sum_{n=0}^{\infty}
%\end{multline}
%\end{proposition}
%
\subsection{Entry 19.}
\begin{proposition}
\begin{equation}\label{eq:gf34}
\sum_{n=0}^{\infty}\left(\frac{x e^{-b x^2}}{x+1}\right)^n \, _2F_2\left(\frac{1}{2}-\frac{n}{2},-\frac{n}{2};\frac{1}{2},1;a+b
   n\right)=-\frac{(x+1) e^{a x^2}}{2 b x^3+2 b x^2-1}
\end{equation}
\end{proposition}
where $f(x)=e^{a x^2} x^m$ and $\phi(x)=(x+1) e^{b x^2}$.
\begin{example}
\begin{multline}\label{eq:gf35}
\sum_{n=0}^{\infty}\frac{\left(\frac{x e^{-b x^2}}{x+1}\right)^n \,
   _2F_2\left(\frac{1}{2}-\frac{n}{2},-\frac{n}{2};\frac{1}{2},1;a+b n\right)}{n+1}=\int_{0}^{x}\frac{(x+1) e^{s^2 (a-b)+b
   x^2}}{(s+1) x}ds
\end{multline}
\end{example}
\subsection{Entry 20.}
\begin{proposition}
\begin{multline}\label{eq:gf36}
\sum_{n=0}^{\infty}\frac{(-1)^{n+1} x^n \, _3F_2(1,1,1-n;2,-n+\alpha +2;-z)}{\Gamma (n) \Gamma (-n+\alpha +2)}=\frac{(-1)^{-m}
   (1-x)^{\alpha } \log (x z+1)}{z \Gamma (\alpha +1)}
\end{multline}
\end{proposition}
where $f(x)=x^m (1-x)^{\alpha } \log (x z+1)$ and $\phi(x)=1$.
\subsection{Entry 21. The hypergeometric function in terms of infinite double sums}
\begin{proposition}
\begin{equation}\label{eq:gf37}
\, _2F_1(a,b;c;x+z)=\sum_{n=0}^{\infty}\sum_{j=0}^{\infty}\frac{\Gamma (c) x^n \Gamma (a+j) \Gamma (b+j)
   z^{j-n}}{n! \Gamma (a) \Gamma (b) \Gamma (c+j) \Gamma (j-n+1)}
\end{equation}
\end{proposition}
where $f(x)=\, _2F_1(a,b;c;x+z)$ and $\phi(x)=1$.
\subsection{Entry 22.}
\begin{proposition}\label{eq:gf38}
%\begin{multline}
%\, _2F_1(a,b;c;x+z)=\sum_{n=0}^{\infty}\sum_{j=0}^{\infty}\frac{\Gamma (c) x^n \Gamma (a+j) \Gamma (b+j) z^{j-n} \Gamma (m+n+1)}{\Gamma (a) \Gamma
%   (b) \Gamma (n+1) (m+n)! \Gamma (c+j) \Gamma (j-n+1)}
%\end{multline}
%%\end{proposition}
%
%\subsection{Entry 22.}
%\begin{proposition}
\begin{equation}\label{eq:gf39}
\, _2F_1\left(a,b;c;(x+z)^{\alpha }\right)=\sum_{n=0}^{\infty}\sum_{j=0}^{\infty}\frac{\Gamma (c) x^n \Gamma (a+j) \Gamma (b+j) \Gamma (j \alpha +1)
   z^{\alpha  j-n}}{n! \Gamma (a) \Gamma (b) \Gamma (j+1) \Gamma (c+j) \Gamma (-n+j \alpha +1)}
\end{equation}
\end{proposition}
where $f(x)=\, _2F_1\left(a,b;c;(x+z)^{\alpha }\right)$ and $\phi(x)=1$.
\subsection{Entry 23.}
\begin{proposition}
\begin{multline}\label{eq:gf40}
\frac{(x+z) \, _2F_1(a,b;c;x+z)}{z}=\sum_{n=0}^{\infty}\sum_{j=0}^{\infty}\frac{\Gamma (c) z^j \Gamma (a+j)
   \Gamma (b+j) \Gamma (j+n+1) \left(\frac{x}{x+z}\right)^n}{n! \Gamma (a) \Gamma
   (b) \Gamma (j+1)^2 \Gamma (c+j)}
\end{multline}
\end{proposition}
where $f(x)=\, _2F_1(a,b;c;x+z)$ and $\phi(x)=x+z$.
\subsection{Entry 24.}
\begin{proposition}
\begin{multline}\label{eq:gf41}
\frac{x^m z^{-m} (x+z) \left(x (x+z)^{-\alpha }\right)^{-m} \, _2F_1\left(a,b;c;(x+z)^{\beta }\right)}{-\alpha 
   x+x+z}\\
   =\sum_{n=0}^{\infty}\sum_{j=0}^{\infty}\frac{\Gamma (c) \Gamma (a+j) \Gamma (b+j) \left(x (x+z)^{-\alpha }\right)^n \Gamma ((m+n) \alpha +j \beta
   +1) z^{\beta  j+(\alpha -1) (m+n)}}{\Gamma (a) \Gamma (b) \Gamma (j+1) \Gamma (n+1) \Gamma (c+j) \Gamma (m+(m+n)
   (\alpha -1)+j \beta +1)}
\end{multline}
\end{proposition}
where $f(x)=x^m \, _2F_1\left(a,b;c;(x+z)^{\beta }\right)$ and $\phi(x)=(x+z)^{\alpha }$.
\subsection{Entry 25.}
\begin{proposition}
\begin{multline}\label{eq:gf42}
\frac{x^m \left(x e^{\alpha  (-x)}\right)^{-m} \, _2F_1\left(a,b;c;e^{x \beta }\right)}{1-\alpha 
   x}\\
   =\sum_{n=0}^{\infty}\sum_{j=0}^{\infty}\frac{\Gamma (c) \Gamma (a+j) \Gamma (b+j) \Gamma (m+n+1) \left(x e^{\alpha  (-x)}\right)^n (\beta  j+\alpha 
   (m+n))^n}{\Gamma (a) \Gamma (b) \Gamma (j+1) \Gamma (n+1) (m+n)! \Gamma (c+j)}
\end{multline}
\end{proposition}
where $f(x)=x^m \, _2F_1\left(a,b;c;e^{x \beta }\right)$ and $\phi(x)=e^{\alpha  x}$.
\subsection{Entry 26.}
\begin{proposition}
\begin{multline}\label{eq:gf43}
(1-x)^{\alpha } e^{\frac{1}{x z+1}}=\sum_{n=0}^{\infty}\sum_{j=0}^{\infty}\frac{(-1)^n \Gamma (\alpha +1) x^n \, _2F_1(j,-n;-n+\alpha +1;-z)}{\Gamma
   (j+1) \Gamma (n+1) \Gamma (-n+\alpha +1)}
\end{multline}
\end{proposition}
where $f(x)=(1-x)^{\alpha } e^{\frac{1}{x z+1}}$ and $\phi(x)=1$.
\subsection{Entry 27.}
\begin{proposition}
For all $a,b,c,x,z,\alpha \in\mathbb{C},m\in\mathbb{Z^{+}}$ then,
\begin{multline}\label{eq:gf44}
\left(\frac{1}{1-x}\right)^{-m} (1-x)^{\alpha +1} \, _1F_2\left(a;b,c;\frac{1}{x z+1}\right)\\
=\sum_{n=0}^{\infty}\sum_{j=0}^{\infty}\frac{\Gamma (b)
   \Gamma (c) (-1)^{2 m+n} \left(\frac{1}{1-x}\right)^n x^n \Gamma (a+j) \Gamma (m+n+\alpha +1) \, _2F_1(j,-n;m+\alpha
   +1;-z)}{\Gamma (a) \Gamma (j+1) \Gamma (n+1) \Gamma (b+j) \Gamma (c+j) \Gamma (m+\alpha +1)}
\end{multline}
\end{proposition}
where $f(x)=x^m (1-x)^{\alpha } \, _1F_2\left(a;b,c;\frac{1}{x z+1}\right)$ and $\phi(x)=1-x$.
\subsection{Generalized Abel Theorem}
In this section we look at a generalized form of Abel's theorem where the original form is recorded in the works of Gessel et al. \cite{gessel} and Abel \cite{abel}.
\begin{proposition}
\begin{equation}\label{eq:gf45}
\sum_{n=0}^{\infty}\frac{\left(x e^{-b x}\right)^n (a+b n)^{n-m}}{(n-m)!}=\frac{e^{a x} x^m}{1-b x}
\end{equation}
\end{proposition}
where $f(x)=e^{a x} x^m$ and $\phi(x)=e^{b x}$.
\begin{example}
\begin{multline}\label{eq:gf46}
\sum_{n=0}^{\infty}\frac{\left(x e^{-b x}\right)^n (a+b (m+n))^n}{(n+1)!}=\frac{x^{m-1} e^{-b (m-1) x} \left(x e^{-b
   x}\right)^{-m} \left(e^{x (a+b (m-1))}-1\right)}{a+b (m-1)}
\end{multline}
\end{example}
\begin{example}
Generalized Gould form.
\begin{multline}\label{eq:gf47}
\sum_{n=0}^{\infty}\frac{(\alpha +\beta  (m+n))^n (a-\beta  (m+n))^{k-m-n}}{n!
   (k-m-n)!}=\frac{e^{\frac{a+\alpha }{\beta }} \left(\frac{1}{\beta
   }\right)^{m-k} \Gamma \left(k-m+1,\frac{a+\alpha }{\beta
   }\right)}{(k-m)!}
\end{multline}
\end{example}
\begin{proof}
Use equation (\ref{eq:gf45}) and apply the contour integral method in \cite{reyn4}. Similar form is recorded in equation (6.10.9) in \cite{hansen}.
\end{proof}
\subsection{Entry 29.}
\begin{proposition}
\begin{multline}\label{eq:gf47}
\sum_{n=0}^{\infty}\frac{x^n (x+1)^{-n} \Gamma (b+m+n+1) \, _2F_1(-a,-n;b+m+1;u)}{\Gamma (n+1) \Gamma (b+m+1)}=(u x+1)^a
   (x+1)^{b+m+1}
\end{multline}
\end{proposition}
where $f(x)=(x+1)^b x^m (u x+1)^a$ and $\phi(x)=x+1$.
\begin{example}
\begin{multline}\label{eq:gf48}
\sum_{n=0}^{\infty}\frac{\left(\frac{x}{x+1}\right)^n \Gamma (b+m+n+1) \, _2F_1(-a,-n;b+m+1;z)}{(m+n+1) \Gamma (n+1) \Gamma
   (b+m+1)}\\
   =\frac{(x+1)^{m+1} F_1(m+1;1-b,-a;m+2;-x,-x z)}{m+1}
\end{multline}
\end{example}
\subsection{Entry 30.}
\begin{proposition}
\begin{multline}
\sum _{n=1}^{\infty }
   \frac{(-1)^n n \left((1-x)^{-a} x\right)^n \Gamma (a n)}{n! \Gamma (2+(-1+a)
   n)} \left(-i m n \left(\, _3F_2\left(1,1,1-n;2,2-n+a n;-\frac{i
   m}{b}\right)\right.\right. \\ \left.\left.
-\, _3F_2\left(1,1,1-n;2,2-n+a n;\frac{i m}{b}\right)\right)+2 b (1+(-1+a) n) \log (b)\right)\\
=\frac{b }{a}\left(-2 \log (b)-\frac{(-1+x) \log \left(b^2+m^2 x^2\right)}{1+(-1+a) x}\right)\\
\end{multline}
\end{proposition}
where $f(x)=(1-x)^{a n} \log \left(b^2+m^2 x^2\right)$ and $\phi(x)=(1-x)^a$.
\subsection{Entry 31.}
\begin{proposition}
\begin{multline}
\frac{b (1-x) x^m \log \left(b^2+a^2 x^2\right)}{1-x+c x}
=\sum
   _{p=0}^{\infty } \frac{(-1)^{m+p} \left((1-x)^{-c} x\right)^{m+p} }{\Gamma (1+p) \Gamma (2-p+c (m+p))}\Gamma (1+c
   (m+p)) \\
\left(-i a p \, _3F_2\left(1,1,1-p;2,2+c m-p+c p;-\frac{i a}{b}\right)\right. \\ \left.
+i
   a p \, _3F_2\left(1,1,1-p;2,2+c m-p+c p;\frac{i a}{b}\right)+b (1-p+c (m+p))
   \log \left(b^2\right)\right)
\end{multline}
\end{proposition}
where $f(x)=x^m \log \left(a^2 x^2+b^2\right)$ and $\phi(x)=(1-x)^c$.
\subsection{Entry 32.}
\begin{proposition}
\begin{multline}
\sum _{n=0}^{\infty } \frac{e^{-b+i n \pi } \left((1-x)^{-a}
   x\right)^n \Gamma (1+a n)}{2 \Gamma (1+n) \Gamma (1-n+a n)}\\ \left(-\, _2F_2\left(\frac{1}{2}-\frac{n}{2},-\frac{n}{2};\frac{1}{2}-\frac{n}{2}+\frac{a
   n}{2},1-\frac{n}{2}+\frac{a n}{2};-z\right)\right. \\ \left.
+e^{2 b} \,
   _2F_2\left(\frac{1}{2}-\frac{n}{2},-\frac{n}{2};\frac{1}{2}-\frac{n}{2}+\frac{a n}{2},1-\frac{n}{2}+\frac{a
   n}{2};z\right)\right)\\
=\frac{(1-x) \sinh \left(b+x^2 z\right)}{1-x+a x}
\end{multline}
\end{proposition}
where $f(x)=(1-x)^{a n} \sinh \left(b+x^2 z\right)$ and $\phi(x)=(1-x)^{-a}$.
%
%\subsection{Entry 33.}
%\begin{proposition}
%
%\end{proposition}
%%
%where $f(x)=$ and $\phi(x)=$.
%%
%\subsection{Entry 34.}
%\begin{proposition}
%
%\end{proposition}
%%
%where $f(x)=$ and $\phi(x)=$.
%
%
\section{Extended generating functions of Prudnikov}
In this section we apply the method of section (7) in \cite{reyn_plos} to equations in \cite{prud2}.
\subsection{Prudnikov equation (5.12.1)}
\begin{theorem}
\begin{equation}\label{eq:gf49}
\sum_{k=0}^{\infty}\frac{\left(z e^{-2 a z}\right)^k H_k(a k+x)}{k!}=\frac{e^{2 x z-z^2}}{1-2
   a z}
\end{equation}
\end{theorem}
\begin{example}
\begin{multline}\label{eq:gf50}
\sum_{k=0}^{\infty}\frac{\left(z e^{-2 a z}\right)^k H_k(a k+x)}{(k+1)!}=\frac{\sqrt{\pi }
   e^{(a-x)^2+2 a z} (\text{erf}(a-x+z)-\text{erf}(a-x))}{2 z}
\end{multline}
\end{example}
\subsection{Prudnikov equation (5.12.26)}
\begin{theorem}
\begin{equation}\label{eq:gf51}
\sum_{n=0}^{\infty}\frac{(a k+1)^{k/2} \left(\sqrt{z} e^{a z}\right)^k
   H_k\left(\frac{x}{\sqrt{a k+1}}\right)}{k!}=\frac{e^{2 x \sqrt{z}-z}}{2 a
   z+1}
\end{equation}
\end{theorem}
\begin{example}
\begin{multline}\label{eq:gf52}
\sum_{n=0}^{\infty}\frac{(a k+1)^{k/2} \left(\sqrt{z} e^{a z}\right)^k
   H_k\left(\frac{x}{\sqrt{a k+1}}\right)}{(k+1)!}\\
   =\frac{\sqrt{\pi }
   e^{-\frac{x^2}{a-1}-a z} \left(\text{erfi}\left(\frac{(a-1)
   \sqrt{z}+x}{\sqrt{a-1}}\right)-\text{erfi}\left(\frac{x}{\sqrt{a-1}}\right)\right)}{2 \sqrt{a-1} \sqrt{z}}
\end{multline}
\end{example}
\subsection{Prudnikov equation (5.12.28)}
\begin{theorem}
\begin{equation}\label{eq:gf53}
\sum_{k=0}^{\infty}\frac{(b k+1)^{k/2} \left(z \left(-e^{2 a z+b
   z^2}\right)\right)^k H_k\left(\frac{a k+x}{\sqrt{b
   k+1}}\right)}{k!}=\frac{e^{-2 x z-z^2}}{2 a z+2 b z^2+1}
\end{equation}
\end{theorem}
\begin{example}
\begin{multline}\label{eq:gf54}
\sum_{k=0}^{\infty}\frac{(b k+1)^{k/2} \left(z \left(-e^{z (2 a+b z)}\right)\right)^k
   H_k\left(\frac{a k+x}{\sqrt{b k+1}}\right)}{(k+1)!}\\
   =\frac{\sqrt{\pi }
   e^{-\frac{(a-x)^2}{b-1}-2 a z-b z^2} \left(\text{erfi}\left(\frac{a+(b-1)
   z-x}{\sqrt{b-1}}\right)-\text{erfi}\left(\frac{a-x}{\sqrt{b-1}}\right)\right)}
   {2 \sqrt{b-1} z}
\end{multline}
\end{example}
\begin{example}
From Eq. (5.11.4.16) in \cite{prud2}
\begin{multline}
\sum _{k=0}^{\infty } \frac{\left(\frac{u}{\sqrt{1+u}}\right)^k L_k^{a-\frac{k}{2}}(x)}{k+n}=u^{-n} (1+u)^{n/2}
   \int_0^u e^{-x z} z^{-1+n} (1+z)^{a-\frac{n}{2}} \, dz
\end{multline}
\end{example}
\begin{example}
From Eq. (5.11.4.16) in \cite{prud2}
\begin{multline}
\sum _{k=0}^{\infty } \left(\frac{u}{\sqrt{1+u}}\right)^k \left(-H_{\frac{1}{2} (-1+k)}+H_{\frac{k}{2}}\right)
   L_k^{a-\frac{k}{2}}(x)=\int_0^u \frac{2 e^{-x z} \sqrt{1+u} (1+z)^a}{\sqrt{1+u} z+u \sqrt{1+z}} \, dz
\end{multline}
\end{example}
\begin{example}
From Eq. (5.11.4.17) in \cite{prud2}
\begin{multline}
\sum _{k=0}^{\infty } \frac{\left(\frac{u}{\sqrt{1+u}}\right)^k L_k^{a-\frac{k}{2}}(x)}{(2 a+k)
   (k+n)}\\
=\int_0^u \frac{1}{2} u^{-n} (1+u)^{n/2} z^{-1+n} (x z)^{-2 a} (1+z)^{-1+a-\frac{n}{2}} (2+z) \Gamma (2
   a,0,x z) \, dz
\end{multline}
\end{example}
\begin{example}
From Eq. (5.11.4.16) in \cite{prud2}
\begin{multline}
\sum _{k=0}^{\infty } \frac{\left(\frac{u}{\sqrt{1+u}}\right)^k \left(-\psi
   ^{(0)}\left(\frac{1+k}{2}\right)+\psi ^{(0)}\left(\frac{2+k}{2}\right)\right) L_k^{a-\frac{k}{2}}(x)}{2
   a+k}\\
=\int_0^u \frac{\sqrt{1+u} (x z)^{-2 a} (1+z)^{-1+a} (2+z) \Gamma (2 a,0,x z)}{\sqrt{1+u} z+u \sqrt{1+z}} \,
   dz
\end{multline}
\end{example}
\begin{example}
From Eq. (5.14.3.3) in \cite{prud2}
\begin{multline}
\sum _{k=0}^{\infty } \frac{\left(\frac{p}{\sqrt{1+p}}\right)^k
   P_k^{\left(a-\frac{k}{2},a-\frac{k}{2}\right)}(x)}{(2 a+k) (k+n)}\\
=\frac{(1+p)^{n/2}}{4 a n (1+n)} \left(2 (1+n)
   F_1\left(n;\frac{1}{2} (-2 a+n),2 a;1+n;-p,\frac{1}{2} p (-1+x)\right)\right. \\ \left.
-n p F_1\left(1+n;\frac{1}{2} (2-2 a+n),2
   a;2+n;-p,\frac{1}{2} p (-1+x)\right)\right)
\end{multline}
\end{example}
\begin{example}
From Eq. (58.81) in \cite{hansen}
\begin{multline}
\sum _{k=0}^{\infty } \frac{2^{-k} \left(\left(-1+w^2\right) z\right)^k I_{a+k}(z)}{(k+1)!}\\
=-\frac{2^{1-a} a
   w^{-a} \left(-2 w^a z^a+2 (w z)^a+2^a z \left(w^a I_{-1+a}(z)-w I_{-1+a}(w z)\right) \Gamma
   (a)\right)}{\left(-1+w^2\right) z^2 \Gamma (1+a)}
\end{multline}
\end{example}
\begin{example}
From Eq. (58.8.7) in \cite{hansen}
\begin{equation}
\sum _{k=0}^{\infty } \frac{2^{-k} \left(\left(-1+w^2\right) z\right)^k I_{a-k}(z)}{(k+1)!}=\frac{2
   \left(-I_{1+a}(z)+w^{1+a} I_{1+a}(w z)\right)}{(-1+w) (1+w) z}
\end{equation}
\end{example}
\begin{example}
From Eq. (58.8.7) in \cite{hansen}
\begin{multline}
\sum _{k=0}^{\infty } \frac{2^{-k} \left(\left(-1+w^2\right) z\right)^k I_{a-k}(z)}{(k+2)!}=\frac{-2
   \left(-1+w^2\right) z I_{1+a}(z)-4 I_{2+a}(z)+4 w^{2+a} I_{2+a}(w z)}{\left(-1+w^2\right)^2 z^2}
\end{multline}
\end{example}
\begin{example}
From Eq. (58.8.7) in \cite{hansen}
\begin{multline}
\sum _{k=0}^{\infty } \frac{2^{-k} \left(\left(-1+w^2\right) z\right)^k
   I_{a-k}(z)}{(k+3)!}\\
=-\frac{\left(-1+w^2\right)^2 z^2 I_{1+a}(z)+4 \left(-1+w^2\right) z I_{2+a}(z)+8 I_{3+a}(z)-8
   w^{3+a} I_{3+a}(w z)}{\left(-1+w^2\right)^3 z^3}
\end{multline}
\end{example}
\begin{example}
From Eq. (58.8.7) in \cite{hansen}
\begin{multline}
\sum _{k=0}^{\infty } \frac{2^{-k} \left(\left(-1+w^2\right) z\right)^k
   I_{a-k}(z)}{(k+4)!}\\
=-\frac{1}{3
   \left(-1+w^2\right)^4 z^4}\left(\left(\left(-1+w^2\right) z \left(24+\left(-1+w^2\right)^2 z^2\right) I_{1+a}(z)\right)\right. \\ \left.
+6
   \left(\left(-1+w^2\right) \left(16+8 a+z^2-w^2 z^2\right) I_{2+a}(z)-8 I_{4+a}(z)+8 w^{4+a} I_{4+a}(w z)\right)\right)
\end{multline}
\end{example}
\begin{example}
From Eq. (58.8.7) in \cite{hansen}
\begin{multline}
\sum _{k=0}^{\infty } \frac{2^{-k} \left(\left(-1+w^2\right) z\right)^k
   I_{a-k}(z)}{(k+5)!}\\
=\frac{1}{12
   \left(-1+w^2\right)^5 z^6}\left(\left(-1+w^2\right) z \left(384 (3+a)-48 \left(-1+w^2\right) z^2-\left(-1+w^2\right)^3
   z^4\right) I_{1+a}(z)\right. \\ \left.
-8 \left(-1+w^2\right) \left(96 (2+a) (3+a)-12 \left(-4-a+(2+a) w^2\right)
   z^2\right.\right. \\ \left.\left.
+\left(-1+w^2\right)^2 z^4\right) I_{2+a}(z)-384 z \left(I_{5+a}(z)-w^{5+a} I_{5+a}(w z)\right)\right)
\end{multline}
\end{example}
\begin{example}
From Eq. (48.17.1) errata in \cite{hansen}
\begin{equation}
\sum _{k=0}^{\infty } \left(y (1+y)^{-1-a}\right)^k L_k^{c+a k}(x)=\frac{e^{-x y} (1+y)^{1+c}}{1-a y}
\end{equation}
\end{example}
\begin{example}
From Eq. (48.17.1) errata in \cite{hansen}
\begin{multline}
\sum _{k=0}^{\infty } \frac{\left(y (1+y)^{-1-a}\right)^k L_k^{c+a k}(x)}{k+b}=y^{-b} (1+y)^{(1+a) b} \int_0^y
   e^{-t x} t^{-1+b} (1+t)^{-((1+a) b)+c} \, dt
\end{multline}
\end{example}
\begin{example}
From Eq. (48.17.2) in \cite{hansen}
\begin{multline}
\sum _{k=0}^{\infty } \frac{4^{-k} \left(1-v^2\right)^k L_k^{c+k}(x)}{k+b}=\left(-1+v^2\right)^{-b} \int_1^v
   2^{1+c} e^{\frac{(-1+t) x}{1+t}} (1+t)^{-c} \left(-1+t^2\right)^{-1+b} \, dt
\end{multline}
\end{example}
\begin{example}
From Eq. (48.17.4) in \cite{hansen}
\begin{multline}
\sum _{k=0}^{\infty } \frac{4^{-k} \left(-1+u^2\right)^k L_k^{c-2 k}(x)}{k+b}=\left(-1+u^2\right)^{-b} \int_1^u
   \frac{2^{-c} e^{-\frac{1}{2} (-1+t) x} (1+t)^c \left(-1+t^2\right)^b}{-1+t} \, dt
\end{multline}
\end{example}
\begin{example}
From Eq. (17.10.1) in \cite{hansen}. Stirling number of the second kind.
\begin{multline}
\sum _{k=1}^{\infty } k^n \cos (k x) (\csc (t+x) \sin (t))^k\\
=\sum _{k=0}^n \cos (t+k t+k x) k! (\csc (t+x) \sin
   (t))^k \left(\csc ^2(t+x) \sin ^2(x)\right)^{\frac{1}{2} (-1-k)} \mathcal{S}_n^{(k)}
\end{multline}
\end{example}
\begin{example}
From Eq. (17.10.1) in \cite{hansen}
\begin{multline}
\sum _{k=1}^{\infty } \frac{k^n \cos (k x) (\csc (t+x) \sin (t))^k}{k+b}\\
=\sin ^{-b}(t) \sin ^b(t+x) \int_0^t
   \left(\sum _{k=0}^n \cos (z+k (x+z)) \csc (x) k! \left(\csc ^2(x+z) \sin ^2(x)\right)^{\frac{1-k}{2}} \right. \\ \left.
\sin ^b(z)
   (\csc (x+z) \sin (z))^{-1+k} \sin ^{-b}(x+z) \mathcal{S}_n^{(k)}\right) \, dz
\end{multline}
\end{example}
\begin{example}
From Eq. (48.19.13) in
\begin{multline}
\sum _{k=0}^{\infty } \frac{4^{-k} \left(-1+v^2\right)^k L_{k+m}^{-1-2 k-2 m}(x) (1+m)_k}{(k+1)!}=\int_1^v
   \frac{2 e^{\frac{1}{2} (1-t) x} t^{-2 m} L_m^{-1-2 m}(t x)}{-1+v^2} \, dt
\end{multline}
\end{example}
\begin{example}
From Eq. (48.19.13) in \cite{hansen}
\begin{multline}
\sum _{k=0}^{\infty } \frac{4^{-k} \left(-1+v^2\right)^k L_{k+m}^{-1-2 k-2 m}(x) (1+m)_k}{k!
   (k+b)}\\
=\left(-1+v^2\right)^{-b} \int_1^v 2 e^{-\frac{1}{2} x (-1+z)} z^{-2 m} \left(-1+z^2\right)^{-1+b} L_m^{-1-2
   m}(x z) \, dz
\end{multline}
\end{example}
\section{Extended Brychkov series.}
\begin{example}
From Eq. (6.2.1.84) in \cite{brychkov}
\begin{multline}
\sum _{k=0}^{\infty } \frac{(-1)^k (a)_{2 k} \psi ^{(0)}(a+2 k) \tan ^{2 k}(u)}{(2 k+1)!}\\
=\cot (u)
   \int_0^u \frac{\left(2 \cos (a s) \psi ^{(0)}(a)-\cos (a s) \log \left(1+\tan ^2(s)\right)-2 s \sin (a
   s)\right) \tan (s)}{\left(1+\tan ^2(s)\right)^{a/2} \sin (2 s)} \, ds
\end{multline}
\end{example}
\begin{example}
From Eq. (6.2.1.84) in \cite{brychkov}
\begin{multline}
\sum _{k=0}^{\infty } \frac{(-1)^k (a)_{2 k} \psi ^{(0)}(a+2 k) \tan ^{2 k}(u)}{(2 k)! (k+b)}\\
=-\frac{2
  }{\tan ^{2 b}(u)} \int_0^u \frac{\left(\cos (a s) \log \left(1+\tan ^2(s)\right)-2 \cos (a s) \psi ^{(0)}(a)+2 s \sin (a s)\right)
   \tan ^{2 b}(s)}{\sin (2 s) \left(1+\tan ^2(s)\right)^{a/2}} \, ds
\end{multline}
\end{example}
\begin{example}
From Eq. (6.14.5.1) in  \cite{brychkov}
\begin{multline}
\sum _{k=0}^{\infty } \frac{\left(\frac{e^w w}{z}\right)^k L_k^{-k+l}((1+k) z)}{(1+k) (k+n)}\\
=\frac{e^{-n w} w^{-n}}{1+l} \int_0^w e^{(-1+n) t} t^{-2+n} (1+t) \left(-z+(t+z) \left(\frac{t+z}{z}\right)^l\right) \,
   dt
\end{multline}
\end{example}
\begin{example}
From Eq. (6.14.5.6(i)) in  \cite{brychkov}
\begin{multline}
\sum _{k=0}^{\infty } \frac{\left(w-\sqrt{-1+v+w^2}\right)^k L_k(z)
   P_k(w)}{k+a}\\
= \int_1^v \frac{e^{\frac{\left(-1+s+w^2-w
   \sqrt{-1+s+w^2}\right) z}{s}} \left(-w+\sqrt{-1+s+w^2}\right)^a J_0\left(\frac{\sqrt{1-w^2}
   \left(w-\sqrt{-1+s+w^2}\right) z}{s}\right)}{2 \sqrt{s} \left(-1+s+w^2-w \sqrt{-1+s+w^2}\right)\left(-w+\sqrt{-1+v+w^2}\right)^{a}} \, ds
\end{multline}
\end{example}
\begin{example}
From Eq. (6.14.5.6(ii)) in  \cite{brychkov}
\begin{multline}
\sum _{k=0}^{\infty } \frac{\left(w+\sqrt{-1+v+w^2}\right)^{-a} \left(-2^{a+k}
   w^{a+k}+\left(w+\sqrt{-1+v+w^2}\right)^{a+k}\right) L_k(z) P_k(w)}{a+k}\\
=
   \int_1^v \frac{e^{\frac{\sqrt{-1+s+w^2} \left(w+\sqrt{-1+s+w^2}\right) z}{s}} \left(w+\sqrt{-1+s+w^2}\right)^a
   J_0\left(\frac{\sqrt{1-w^2} \left(w+\sqrt{-1+s+w^2}\right) z}{s}\right)}{2 \sqrt{s} \left(-1+s+w^2+w
   \sqrt{-1+s+w^2}\right)\left(w+\sqrt{-1+v+w^2}\right)^{a}} \, ds
\end{multline}
\end{example}
\begin{example}
From Eq. (6.14.5.2) in  \cite{brychkov}
\begin{multline}
\sum _{k=0}^{\infty } \frac{(1+k)^{-1+k} \left(-e^w w\right)^k L_k\left(\frac{z}{1+k}\right)}{k!
   (k+a)}=\frac{(w \exp (w))^{-a} }{\sqrt{z}}\int_0^w e^{(-1+a) t} t^{-\frac{3}{2}+a} (1+t) I_1\left(2 \sqrt{z t}\right) \,
   dt
\end{multline}
\end{example}
\begin{example}
From Eq. (6.15.4.1) in  \cite{brychkov}
\begin{multline}
\sum _{k=0}^{\infty } \frac{2^{-k} \left(\frac{e^w w}{z}\right)^k C_k^{(-k+\lambda )}(1+(1+k) z)}{(1+k)
   (1-\lambda )_k (k+a)}\\
=\frac{2 z \lambda  e^{-a w} w^{-a} }{1+2 \lambda }\int_0^w e^{(-1+a) t} t^{a-2} (1+t)
   \left(1-\, _1F_1\left(-\frac{1}{2}-\lambda ;-2 \lambda ;-\frac{2 t}{z}\right)\right) \, dt
\end{multline}
\end{example}
\begin{example}
From Eq. (6.13.3.1) in  \cite{brychkov}
\begin{multline}
\sum _{k=0}^{\infty } \frac{(1+k)^{-1+k} \left(e^w w\right)^k H_{2 k}\left(\frac{z}{\sqrt{1+k}}\right)}{(2
   k)! (k+a)}\\
=e^{-a w} w^{-a} \int_0^w \frac{e^{(-1+a) t} t^{-2+a} (1+t) \left(1-\cosh \left(2 \sqrt{t} z\right)+2
   \sqrt{t} z \sinh \left(2 \sqrt{t} z\right)\right)}{2 z^2} \, dt
\end{multline}
\end{example}
\begin{example}
From Eq. (6.2.1.86) in  \cite{brychkov}
\begin{multline}
\sum _{k=0}^{\infty } \frac{(-1)^k (a)_{1+2 k} \psi ^{(0)}(1+a+2 k) \tan ^{1+2 k}(u)}{(1+2 k)! (k+b)}\\
=\cos ^{-1+2 b}(u) \sin
   ^{1-2 b}(u) \int_0^u \frac{\sin ^{-2+2 b}(t) \left(2 t \cos (a t)-\left(\log \left(\sec ^2(t)\right)-2 \psi ^{(0)}(a)\right) \sin (a
   t)\right)}{\cos ^{2 b}(t) \sec ^2(t)^{a/2}} \, dt
\end{multline}
\end{example}
\begin{example}
From Eq. (6.11.4.2) in  \cite{brychkov}
\begin{multline}
\sum _{k=0}^{\infty } \frac{(1+k)^{-1+\frac{k}{2}} \left(-e^w w\right)^k (1+k+z)^{k/2} P_k\left(\frac{2+2 k+z}{2 \sqrt{(1+k)
   (1+k+z)}}\right)}{k! (k+b)}\\
=e^{-b w} w^{-b} \int_0^w e^{-\frac{1}{2} t (2-2 b+z)} t^{-1+b} (1+t) \left(I_0\left(\frac{t
   z}{2}\right)+I_1\left(\frac{t z}{2}\right)\right) \, dt
\end{multline}
\end{example}
\begin{example}
From Eq. (6.13.3.2) in  \cite{brychkov}
\begin{multline}
\sum _{k=0}^{\infty } \frac{(1+k)^{-\frac{1}{2}+k} \left(e^w w\right)^k }{(1+2 k)!
   (k+b)}H_{1+2 k}\left(\frac{z}{\sqrt{1+k}}\right)\\
=\frac{2 }{z (w \exp (w))^b}\int_0^w e^{(-1+b) t} t^{-2+b} (1+t) \sinh ^2\left(\sqrt{t} z\right) \, dt
\end{multline}
\end{example}
\begin{example}
From Eq. (6.14.2.6) in  \cite{brychkov}. An extension of Laplace 's first integral form see pp. 171 in \cite{rainville}.
\begin{multline}
\sum _{k=0}^{\infty } \frac{\left(w-\sqrt{-1+v+w^2}\right)^k L_k(z) P_k(w)}{k+b}\\
=\int_1^v \frac{e^{-\frac{\sqrt{-1+t+w^2}
   \left(w-\sqrt{-1+t+w^2}\right) z}{t}} \left(w-\sqrt{-1+t+w^2}\right)^b J_0\left(\frac{\sqrt{1-w^2} \left(w-\sqrt{-1+t+w^2}\right)
   z}{t}\right)}{2 \sqrt{t} \left(-1+t+w^2-w \sqrt{-1+t+w^2}\right) \left(w-\sqrt{-1+v+w^2}\right)^b} \, dt
\end{multline}
where $Re(v)>1$
\end{example}
\begin{example}
From Eq. (6.14.2.6) in  \cite{brychkov}. An extension of Laplace 's first integral form.
\begin{multline}
\sum _{k=0}^{\infty } \frac{\left(\left(w+\sqrt{-1+v+w^2}\right)^{b+k}-(2 w)^{b+k}\right) L_k(z) P_k(w)}{b+k}\\
=\int_1^v
   \frac{e^{\frac{\sqrt{-1+t+w^2} \left(w+\sqrt{-1+t+w^2}\right) z}{t}} \left(w+\sqrt{-1+t+w^2}\right)^b
   J_0\left(\frac{\sqrt{1-w^2} \left(w+\sqrt{-1+t+w^2}\right) z}{t}\right)}{2 \sqrt{t} \left(-1+t+w^2+w \sqrt{-1+t+w^2}\right)} \,
   dt
\end{multline}
\end{example}
\begin{example}
From Eq. (6.15.4.2) in  \cite{brychkov}
\begin{multline}
\sum _{k=0}^{\infty } \frac{(1+k)^{-1+k} \left(-e^w w\right)^k C_{2 k}^{(-k+\lambda
   )}\left(\frac{z}{\sqrt{1+k}}\right)}{(1-\lambda )_k (k+b)}\\
=\frac{e^{-b w} w^{-b}}{2 z^2 (1-\lambda )} \int_0^w e^{(-1+b) t} t^{-2+b}
   (1+t) \left(-1+\, _1F_1\left(-1+\lambda ;-\frac{1}{2};t z^2\right)\right) \, dt
\end{multline}
\end{example}
\begin{example}
From Eq. (6.15.4.3) errata in  \cite{brychkov}
\begin{equation}
\sum _{k=0}^{\infty } \frac{(1+k)^{-\frac{1}{2}+k} \left(-e^w w\right)^k C_{1+2 k}^{(-k+\lambda
   )}\left(\frac{z}{\sqrt{1+k}}\right)}{(1-\lambda )_k}=\frac{e^{-w} \left(-1+\, _1F_1\left(\lambda ;\frac{1}{2};w
   z^2\right)\right)}{w z}
\end{equation}
\end{example}
\begin{example}
From Eq. (6.15.4.3) in  \cite{brychkov}
\begin{multline}
\sum _{k=0}^{\infty } \frac{(1+k)^{-\frac{1}{2}+k} \left(-e^w w\right)^k C_{1+2 k}^{(-k+\lambda
   )}\left(\frac{z}{\sqrt{1+k}}\right)}{(1-\lambda )_k (k+b)}\\
=\frac{e^{-b w} w^{-b} }{z}\int_0^w e^{(-1+b) t} t^{-2+b}
   (1+t) \left(-1+\, _1F_1\left(\lambda ;\frac{1}{2};t z^2\right)\right) \, dt
\end{multline}
\end{example}
\begin{example}
From Eq. (6.15.4.4) in  \cite{brychkov}
\begin{multline}
\sum _{k=0}^{\infty } \frac{\left(-e^w w\right)^k (k+z)^k C_{2 k}^{(-k+\lambda
   )}\left(\sqrt{\frac{-1+z}{k+z}}\right)}{(1+k) (1-\lambda )_k (k+b)}\\
=\frac{\left(e^{-b w} w^{-b}\right) }{(-1+z) (1+2 \lambda
   )}\int_0^w e^{(-1+b) t}
   t^{-2+b} (1+t) \left(1-\, _1F_1\left(-\frac{1}{2}-\lambda ;-\frac{1}{2};t-t z\right)\right) \, dt
\end{multline}
\end{example}
\begin{example}
From Eq. (6.15.4.5) in  \cite{brychkov}
\begin{multline}
\sum _{k=0}^{\infty } \frac{\left(-e^w w\right)^k (k+z)^{\frac{1}{2}+k} C_{1+2 k}^{(-k+\lambda
   )}\left(\sqrt{\frac{-1+z}{k+z}}\right)}{(1+k) (1-\lambda )_k (k+b)}\\
=\frac{2 e^{-b w} w^{-b} \lambda  }{\sqrt{-1+z}
   (1+2 \lambda )}\int_0^w
   e^{(-1+b) t} t^{-2+b} (1+t) \left(-1+\, _1F_1\left(-\frac{1}{2}-\lambda ;\frac{1}{2};t-t z\right)\right) \, dt
\end{multline}
\end{example}
\begin{example}
From Eq. (5.2.13.29) in \cite{prud1} and \cite{reyn_plos}.
\begin{equation}
\sum _{k=0}^{\infty } \frac{\left((-1+y) y^{-v}\right)^k (u+k v-1)!}{k!
   (u+k v-k)!}=\frac{y^u}{u}
\end{equation}
\end{example}
\begin{example}
From Eq. (5.2.13.29) in \cite{prud1} and \cite{reyn_plos}.
\begin{equation}
\sum _{k=0}^{\infty } \frac{\left((-1+y) y^{-v}\right)^k (u+k v-1)!}{k! (u+1+k v-k)!}=\frac{y^u (u y+v (1+u-u
   y))}{u (1+u) (u+v)}
\end{equation}
\end{example}
\begin{example}
From Eq. (5.2.13.29) in \cite{prud1} and \cite{reyn_plos}.
\begin{multline}
\sum _{k=0}^{\infty } \frac{\left((-1+y) y^{-v}\right)^k (u+k v-1)!}{k! (u+2+k v-k)!}\\
=\frac{\exp
   \left(\frac{(2+u) (-\log (1-y)+v \log (y))}{-1+v}\right) }{u (u+v)(u+2 v)}\left(-v^2 \, _2F_1\left(\frac{2+u}{1-v},\frac{u+2
   v}{1-v};\frac{1+u+v}{1-v};1\right)\right. \\ \left.
+y^{\frac{u+2 v}{1-v}} \left(\frac{(u+2 v) (1-y)^{\frac{1+u}{-1+v}} y}{(1+u) (2+u)
   (1+u+v)} \left(-u
   (2+u) (-1+v)^2 (1-y)^{\frac{v}{-1+v}}\right.\right.\right. \\ \left.\left.\left.
+(u+v) (1-y)^{\frac{1}{-1+v}} (2+u+(-1+v) y)\right)\right.\right. \\ \left.\left.
+v^2 \, _2F_1\left(\frac{2+u}{1-v},\frac{u+2 v}{1-v};\frac{1+u+v}{1-v};y\right)\right)\right)
\end{multline}
\end{example}
\begin{example}
Bell number.
\begin{equation}
\sum _{n=0}^{\infty } \frac{B_n x^n}{(n+a) n!}=x^{-a} \int_0^x e^{-1+e^t} t^{-1+a} \, dt
\end{equation}
\end{example}
\begin{example}
From pp. 116 in \cite{rainville}
\begin{multline}
\sum _{n=-\infty }^{\infty } \frac{\left(\left(-t+\sqrt{-1+t^2}\right)^n-(-1)^n e^{i \pi  \alpha } \left(1-2
   t \left(t+\sqrt{-1+t^2}\right)\right)^{\alpha }\right) J_n(z)}{n+2 \alpha }\\
=-\int_1^t \frac{e^{-\sqrt{-1+r^2} z}
   \left(1-2 r \left(r+\sqrt{-1+r^2}\right)\right)^{\alpha } \left(1+2 t
   \left(-t+\sqrt{-1+t^2}\right)\right)^{\alpha }}{\sqrt{-1+r^2}} \, dr
\end{multline}
\end{example}
\begin{example}
From pp. 184 in \cite{rainville}
\begin{multline}
\sum _{n=0}^{\infty } \frac{\left(x-\sqrt{-1+x^2+\rho ^2}\right)^n \, _1F_2(-n;1,1;y)
   P_n(x)}{n+b}\\
=\int_1^{\rho } \frac{\left(-x+\sqrt{-1+t^2+x^2}\right)^b }{\left(-x+\sqrt{-1+x^2+\rho ^2}\right)^b
   \left(1-t^2-x^2+x \sqrt{-1+t^2+x^2}\right)}\\
J_0\left(\frac{\sqrt{2}
   \sqrt{x-\sqrt{-1+t^2+x^2}} \sqrt{-t+\sqrt{-1+t^2+x^2}} \sqrt{y}}{t}\right)\\
 J_0\left(\frac{\sqrt{2}
   \sqrt{x-\sqrt{-1+t^2+x^2}} \sqrt{t+\sqrt{-1+t^2+x^2}} \sqrt{y}}{t}\right) \, dt
\end{multline}
where $Re(x)>1$
\end{example}
\begin{example}
From pp. 184 in \cite{rainville}
\begin{multline}
\sum _{n=0}^{\infty } \frac{\left(-2^n x^n \left(x+\sqrt{x^2}\right)^b+\left(x+\sqrt{-1+x^2+\rho
   ^2}\right)^{b+n}\right) \, _1F_2(-n;1,1;y) P_n(x)}{(b+n) \left(x+\sqrt{-1+x^2+\rho ^2}\right)^b}\\
=\int_1^{\rho } \frac{\left(x+\sqrt{-1+t^2+x^2}\right)^b }{\left(x+\sqrt{-1+x^2+\rho ^2}\right)^b
   \left(-1+t^2+x \left(x+\sqrt{-1+t^2+x^2}\right)\right)}\\
I_0\left(\sqrt{2} \sqrt{\frac{1}{t^2}}
   \sqrt{-t+\sqrt{-1+t^2+x^2}} \sqrt{x+\sqrt{-1+t^2+x^2}} \sqrt{y}\right)\\
 I_0\left(\sqrt{2} \sqrt{\frac{1}{t^2}}
   \sqrt{t+\sqrt{-1+t^2+x^2}} \sqrt{x+\sqrt{-1+t^2+x^2}} \sqrt{y}\right) \, dt
\end{multline}
where $Re(x)<-1$
\end{example}
\begin{example}
\begin{multline}
\sum _{k=0}^{\infty } \frac{\left(e^{-u y} u\right)^k B_k(x+k y)}{k! (k+a)}=e^{a u y} u^{-a}
   \int_0^{\frac{u}{2}} 2^a e^{2 t (x-a y)} t^a (-1+\coth (t)) \, dt
\end{multline}
\end{example}
\begin{example}
From PP. 282 question (31) in \cite{rainville}.
\begin{multline}
\sum _{n=0}^{\infty } \frac{\left(x-\sqrt{-1+x^2+\rho ^2}\right)^n C_n^{(v)}(x)
   \left(\frac{1}{2}+v\right)_n}{\left((a+n)^2-b^2\right) (2 v)_n}\\
=-\int_1^{\rho } \frac{\left(-\left(-x+\sqrt{-1+t^2+x^2}\right)^{-2
   b}+\left(-x+\sqrt{-1+x^2+\rho ^2}\right)^{-2 b}\right)}{b \left(-1+t^2+x^2-x \sqrt{-1+t^2+x^2}\right)} \,
   \\
2^{-\frac{3}{2}+v}
   \left(-x+\sqrt{-1+t^2+x^2}\right)^{a+b} \left(1+t+x \left(-x+\sqrt{-1+t^2+x^2}\right)\right)^{\frac{1}{2}-v}\\
  \left(-x+\sqrt{-1+x^2+\rho ^2}\right)^{-a+b} dt
\end{multline}
\end{example}
\begin{example}
From Eq. (5.2.9.4(i)) in \cite{prud1}
\begin{multline}
\sum _{k=1}^{\infty } \frac{\left(a e^{-a}\right)^k k^{-1+k}}{(k+1)!
   (k+b)}=\frac{1}{a}\left(\frac{1+a+\left(1+a+a^2\right)
   b}{b^2}-\frac{e^a}{-1+b}-\frac{e^{a b} E_{2-b}(a b)}{b^2}\right)
\end{multline}
\end{example}
\begin{example}
\begin{equation}
\sum _{k=1}^{\infty } \frac{\left(a e^{-a}\right)^k k^{-1+k}}{(k+2)!}=\frac{-1+6 a^2+4 a^3+e^{2 a}+a \left(6-8
   e^a\right)}{8 a^2}
\end{equation}
\end{example}
\begin{example}
\begin{multline}
\sum _{k=1}^{\infty } \frac{\left(a e^{-a}\right)^k k^{-1+k}}{(k+3)!}=\frac{8+3 a (-19+6 a (11+a (11+6 a)))-324
   a^2 e^a+81 a e^{2 a}-8 e^{3 a}}{648 a^3}
\end{multline}\end{example}
\begin{example}
\begin{multline}
\sum _{k=1}^{\infty } \frac{\left(a e^{-a}\right)^k k^{-1+k}}{(k+4)!}\\
=\frac{1}{82944 a^4}\left(-81+4 a (175+6 a (-115+12 a (25+a
   (25+12 a))))\right. \\ \left.
   -13824 a^3 e^a+5184 a^2 e^{2 a}-1024 a e^{3 a}+81 e^{4 a}\right)
\end{multline}\end{example}
\begin{example}
\begin{multline}
\sum _{k=1}^{\infty } \frac{\left(a e^{-a}\right)^k k^{-1+k}}{(k+5)!}\\
=\frac{1}{1296000000 a^5}\left(82944+5 a (-170181+20 a (39743+30 a
   (-3799+60 a (137+a (137+60 a)))))\right. \\ \left.
-54000000 a^4 e^a+27000000 a^3 e^{2 a}-8000000 a^2 e^{3 a}+1265625 a e^{4 a}-82944\right)
   e^{5 a}
\end{multline}\end{example}
\begin{example}
\begin{multline}
\sum _{n=0}^{\infty } \frac{\left(1+(-1)^n\right) (a+b (m+n))^{n/2} \left(e^{-b x^2} x\right)^n}{2 \Gamma
   \left(\frac{4+n}{2}\right)}=\frac{-e^{2 b x^2}+e^{a x^2} x^m \left(e^{-b x^2} x\right)^{-m}}{(a+b (-2+m))
   x^2}
\end{multline}
\end{example}
\begin{example}
From Eq. (5.2.18.4) in \cite{prud1}
\begin{equation}
\sum _{k=1}^{\infty } \frac{k^{-k} x^k}{(a+k)^2-b^2}=\int _0^1\int _0^x\frac{t^{-t z} x^{-a+b} z^{a+b}
   \left(z^{-2 b}-x^{-2 b}\right)}{2 b}dzdt
\end{equation}
\end{example}
\begin{example}
\begin{multline}
\sum _{k=1}^{\infty } k^{-k} x^k (-\psi ^{(0)}(-b+k)+\psi ^{(0)}(b+k))=\int _0^1\int _0^x\frac{t^{-t z}
   x^{1-b} z^{-b} \left(-x^{2 b}+z^{2 b}\right)}{-x+z}dzdt
\end{multline}
\end{example}
\begin{example}
\begin{equation}
\sum _{n=1}^{\infty } \text{csch}(a n \pi )=\frac{i \pi -\psi _{e^{a
   \pi }}^{(0)}(1)+\psi _{e^{a \pi }}^{(0)}\left(1-\frac{i}{a}\right)}{a
   \pi }
\end{equation}
\end{example}
\begin{example}
\begin{equation}
\sum _{n=1}^{\infty } \text{csch}(n z)=\frac{i \pi -\psi _{e^z}^{(0)}(1)+\psi _{e^z}^{(0)}\left(1-\frac{i \pi
   }{z}\right)}{z}
\end{equation}
\end{example}
\begin{example}
From PP. 168 question (1) in \cite{rainville}.
\begin{multline}
\sum _{n=0}^{\infty } \frac{\left(x-\sqrt{-1+x^2+\rho ^2}\right)^n P_n(x) (1-c)_n (c)_n}{(n!)^2
   (n+a)}\\
=-\int_1^{\rho } \frac{s \left(-x+\sqrt{-1+s^2+x^2}\right)^a }{\left(1-s^2-x^2+x \sqrt{-1+s^2+x^2}\right) \left(-x+\sqrt{-1+x^2+\rho
   ^2}\right)^a}\\
\, _2F_1\left(1-c,c;1;\frac{1}{2}
   \left(1-s+x-\sqrt{-1+s^2+x^2}\right)\right)\\
 \, _2F_1\left(1-c,c;1;\frac{1}{2}
   \left(1-s-x+\sqrt{-1+s^2+x^2}\right)\right) \, ds
\end{multline}
\end{example}
\begin{example}
\begin{multline}
\sum _{n=0}^{\infty } \frac{\left(x+\sqrt{-1+x^2+\rho ^2}\right)^n P_n(x) (1-c)_n (c)_n}{(n!)^2
   (n+a)}\\
=\int_1^{\rho }\, _2F_1\left(1-c,c;1;\frac{1}{2} \left(1-s-x-\sqrt{-1+s^2+x^2}\right)\right)\\
 \,
   _2F_1\left(1-c,c;1;\frac{1}{2} \left(1-s+x+\sqrt{-1+s^2+x^2}\right)\right)\\
 \frac{s\, ds }{\left(-1+s^2+x^2+x
   \sqrt{-1+s^2+x^2}\right) \left(x+\sqrt{-1+s^2+x^2}\right)^{-a} \left(x+\sqrt{-1+x^2+\rho ^2}\right)^a} 
\end{multline}
\end{example}
\begin{example}
From PP. 271 equation (6) in \cite{rainville}
\begin{multline}
\sum _{n=0}^{\infty } \frac{\left(x-\sqrt{-1+x^2+\rho ^2}\right)^n P_n^{(\alpha ,\beta )}(x)}{n+b}\\
=-\int_1^{\rho
   } \frac{2^{\alpha +\beta } \left(\frac{1}{1+t+x-\sqrt{-1+t^2+x^2}}\right)^{\beta }
   \left(-x+\sqrt{-1+t^2+x^2}\right)^b \left(\frac{1}{1+t-x+\sqrt{-1+t^2+x^2}}\right)^{\alpha }}{\left(1-t^2-x^2+x
   \sqrt{-1+t^2+x^2}\right) \left(-x+\sqrt{-1+x^2+\rho ^2}\right)^b} \, dt
\end{multline}
\end{example}
\begin{example}
\begin{multline}
\sum _{n=0}^{\infty } \frac{\left(x+\sqrt{-1+x^2+\rho ^2}\right)^n P_n^{(\alpha ,\beta )}(x)}{n+b}\\
=\int_1^{\rho
   } \frac{2^{\alpha +\beta } \left(\frac{1}{1+s-x-\sqrt{-1+s^2+x^2}}\right)^{\alpha }
   \left(x+\sqrt{-1+s^2+x^2}\right)^b \left(\frac{1}{1+s+x+\sqrt{-1+s^2+x^2}}\right)^{\beta }}{\left(-1+s^2+x^2+x
   \sqrt{-1+s^2+x^2}\right) \left(x+\sqrt{-1+x^2+\rho ^2}\right)^b} \, ds
\end{multline}
where $Re(x)<1$.
\end{example}
\begin{example}
PP. 272 question (2) in \cite{rainville}. Extended Brafman form.
\begin{multline}
\sum _{n=0}^{\infty } \frac{\left(x-\sqrt{-1+x^2+\rho ^2}\right)^n P_n^{(\alpha ,\beta )}(x) (1+\alpha +\beta
   -\gamma )_n (\gamma )_n}{(1+\alpha )_n (1+\beta )_n (n+b)}\\
=-\int_1^{\rho } \frac{s
   \left(-x+\sqrt{-1+s^2+x^2}\right)^b }{\left(1-s^2-x^2+x \sqrt{-1+s^2+x^2}\right) \left(-x+\sqrt{-1+x^2+\rho
   ^2}\right)^b}\\
\, _2F_1\left(1+\alpha +\beta -\gamma ,\gamma ;1+\alpha ;\frac{1}{2}
   \left(1-s-x+\sqrt{-1+s^2+x^2}\right)\right)\\
 \, _2F_1\left(1+\alpha +\beta -\gamma ,\gamma ;1+\beta ;\frac{1}{2}
   \left(1-s+x-\sqrt{-1+s^2+x^2}\right)\right) \, ds
\end{multline}
\end{example}
\begin{example}
\begin{multline}
\sum _{n=0}^{\infty } \frac{\left(x+\sqrt{-1+x^2+\rho ^2}\right)^n P_n^{(\alpha ,\beta )}(x) (1+\alpha +\beta
   -\gamma )_n (\gamma )_n}{(1+\alpha )_n (1+\beta )_n (n+b)}\\
=\int_1^{\rho } \frac{s \left(x+\sqrt{-1+s^2+x^2}\right)^b
   }{\left(x+\sqrt{-1+x^2+\rho ^2}\right)^b \left(-1+s^2+x^2+x
   \sqrt{-1+s^2+x^2}\right)}\\
\, _2F_1\left(1+\alpha +\beta -\gamma ,\gamma ;1+\alpha ;\frac{1}{2} \left(1-s-x-\sqrt{-1+s^2+x^2}\right)\right)\\
 \,
   _2F_1\left(1+\alpha +\beta -\gamma ,\gamma ;1+\beta ;\frac{1}{2}
   \left(1-s+x+\sqrt{-1+s^2+x^2}\right)\right) \, ds
\end{multline}
\end{example}
\begin{example}
From PP. 280 equation (21) in \cite{rainville}.
\begin{multline}
\sum _{n=0}^{\infty } \frac{\left(x-\sqrt{-1+x^2+\rho ^2}\right)^n C_n^{(v)}(x) (2 v-\gamma )_n (\gamma )_n}{(2
   v)_n \left(\frac{1}{2}+v\right)_n (n+b)}\\
=-\int_1^{\rho } \frac{s \left(-x+\sqrt{-1+s^2+x^2}\right)^b }{\left(1-s^2-x^2+x
   \sqrt{-1+s^2+x^2}\right) \left(-x+\sqrt{-1+x^2+\rho ^2}\right)^b}\\
\, _2F_1\left(2
   v-\gamma ,\gamma ;\frac{1}{2}+v;\frac{1}{2} \left(1-s+x-\sqrt{-1+s^2+x^2}\right)\right)\\
 \, _2F_1\left(2 v-\gamma
   ,\gamma ;\frac{1}{2}+v;\frac{1}{2} \left(1-s-x+\sqrt{-1+s^2+x^2}\right)\right) \, ds
\end{multline}
\end{example}
\begin{example}
\begin{multline}
\sum _{n=0}^{\infty } \frac{\left(x+\sqrt{-1+x^2+\rho ^2}\right)^n C_n^{(v)}(x) (2 v-\gamma )_n (\gamma )_n}{(2
   v)_n \left(\frac{1}{2}+v\right)_n (n+b)}\\
=\left(x+\sqrt{-1+x^2+\rho ^2}\right)^b \int_1^{\rho } \frac{s
   \left(x+\sqrt{-1+s^2+x^2}\right)^b }{-1+s^2+x^2+x \sqrt{-1+s^2+x^2}}\\
\, _2F_1\left(2 v-\gamma ,\gamma ;\frac{1}{2}+v;\frac{1}{2}
   \left(1-s-x-\sqrt{-1+s^2+x^2}\right)\right)\\
 \, _2F_1\left(2 v-\gamma ,\gamma ;\frac{1}{2}+v;\frac{1}{2}
   \left(1-s+x+\sqrt{-1+s^2+x^2}\right)\right) \, ds
\end{multline}
\end{example}
\begin{example}
From PP. 469 Eq. (7) in \cite{magnus} and \cite{reyn4}.
\begin{multline}
\sum _{n=-\infty }^{\infty } e^{d n} \Phi \left(-e^{2 (m+d n)},-k,\frac{1}{2} (1+\log (a))\right)\\
=\sum
   _{n=0}^{\infty } \frac{2^{-2-k} e^{-\frac{m (d+2 i n \pi )}{d}} \pi  \left(\left(-\frac{2 i n \pi }{d}+\log
   (a)\right)^k+e^{\frac{4 i m n \pi }{d}} \left(\frac{2 i n \pi }{d}+\log (a)\right)^k\right)
   \text{sech}\left(\frac{n \pi ^2}{d}\right)}{d}
\end{multline}
\end{example}
\begin{example}
\begin{multline}
\sum _{n=-\infty }^{\infty } e^{-d n} \log \left(\frac{1}{2} (1-\tanh (m+d n))\right)\\
=-\sum _{n=0}^{\infty }
   \frac{e^m \pi  \text{sech}\left(\frac{n \pi ^2}{d}\right) \left(d \cos \left(\frac{2 m n \pi }{d}\right)+2 n \pi 
   \sin \left(\frac{2 m n \pi }{d}\right)\right)}{d^2+4 n^2 \pi ^2}
\end{multline}
\end{example}
\begin{example}
\begin{multline}
\sum _{n=-\infty }^{\infty } q^{-n} \log \left(1+e^{2 m} q^{2 n}\right)\\
=-\sum _{n=0}^{\infty } \frac{e^m \pi 
   \text{sech}\left(\frac{n \pi ^2}{\log (q)}\right) \left(\cos \left(\frac{2 m n \pi }{\log (q)}\right) \log (q)+2 n
   \pi  \sin \left(\frac{2 m n \pi }{\log (q)}\right)\right)}{4 n^2 \pi ^2+\log ^2(q)}
\end{multline}
\end{example}
\begin{example}
\begin{multline}
\sum _{n=-\infty }^{\infty } e^{n s} \left(i \tanh ^{-1}\left(\sqrt[4]{-1} e^{\frac{m}{2}+n s}\right)+\tanh
   ^{-1}\left((-1)^{3/4} e^{\frac{m}{2}+n s}\right)\right)\\
=\sum _{n=0}^{\infty } \frac{(-1)^{3/4} e^{-\frac{m}{2}} \pi
    \text{sech}\left(\frac{n \pi ^2}{2 s}\right) \left(-s \cos \left(\frac{m n \pi }{s}\right)+2 n \pi  \sin
   \left(\frac{m n \pi }{s}\right)\right)}{2 \left(4 n^2 \pi ^2+s^2\right)}
\end{multline}
\end{example}
\begin{example}
From Eq. (5.3.1.11) in \cite{prud2}.
\begin{equation}
\sum _{k=0}^{\infty } \frac{(a)_k \zeta (k+a) t^k}{(k+1)!}=\frac{\zeta (-1+a,1-t)-\zeta (-1+a)}{(-1+a) t}
\end{equation}
\end{example}
\begin{example}
\begin{multline}
\sum _{k=0}^{\infty } \frac{t^k (a)_k \zeta (a+k)}{(2+k)!}=\frac{\zeta (-2+a,1-t)-\zeta (-2+a)-(-2+a) t \zeta
   (-1+a)}{(-2+a) (-1+a) t^2}
\end{multline}
\end{example}
\begin{example}
\begin{multline}
\sum _{k=0}^{\infty } \frac{t^k (a)_k \zeta (a+k)}{(3+k)!}\\
=-\frac{-2 \zeta (-3+a,1-t)+2 \zeta (-3+a)+(-3+a) t
   (2 \zeta (-2+a)+(-2+a) t \zeta (-1+a))}{2 (-3+a) \left(2-3 a+a^2\right) t^3}
\end{multline}
\end{example}
\begin{example}
\begin{multline}
\sum _{k=0}^{\infty } \frac{t^k (a)_k \zeta (a+k)}{(4+k)!}\\
=-\frac{1}{6 (-4+a) \left(-6+11 a-6 a^2+a^3\right)
   t^4}\left(-6 \zeta (-4+a,1-t)+6 \zeta (-4+a)\right. \\ \left.
+(-4+a) t
   (6 \zeta (-3+a)+(-3+a) t (3 \zeta (-2+a)+(-2+a) t \zeta (-1+a)))\right)
\end{multline}
\end{example}
\begin{example}
\begin{multline}
\sum _{k=0}^{\infty } \frac{t^k (a)_k \zeta (a+k)}{(5+k)!}\\
=-\frac{1}{24 (-5+a)
   \left(24-50 a+35 a^2-10 a^3+a^4\right) t^5}\left(-24 \zeta (-5+a,1-t)+24 \zeta (-5+a)\right. \\ \left.
+(-5+a)
   t (24 \zeta (-4+a)+(-4+a) t (12 \zeta (-3+a)\right. \\ \left.
+(-3+a) t (4 \zeta (-2+a)+(-2+a) t \zeta (-1+a))))\right)
\end{multline}
\end{example}
\begin{example}
From Eq. (21) in \cite{brafman}
\begin{multline}
\sum _{n=0}^{\infty } \frac{(1+a)_n \, _2F_1\left(-n,1+a+b+n;1+a;\frac{1-x}{2}\right) z^n}{(n+1)!}\\
=\int_0^z \frac{2^{a+b} \left(1-s+\sqrt{1+s^2-2 s x}\right)^{-a}
   \left(1+s+\sqrt{1+s^2-2 s x}\right)^{-b}}{\sqrt{1+s^2-2 s x} z} \, ds
\end{multline}
\end{example}
\begin{example}
\begin{equation}
\sum _{n=0}^{\infty } \frac{(b)_n P_n(x) t^n}{(n+1)!}=\int_0^t \frac{\, _2F_1\left(\frac{b}{2},\frac{1+b}{2};1;\frac{\left(-1+x^2\right) s^2}{(-1+x s)^2}\right) (1-x s)^{-b}}{t}
   \, ds
\end{equation}
\end{example}
\begin{example}
\begin{multline}
\sum _{n=0}^{\infty } \frac{(a)_n H_n(x) t^n}{(n+1)!}=\frac{1}{t}\int_0^t \, _2F_0\left(\frac{1}{2}+\frac{a}{2},\frac{a}{2};;-\frac{4 s^2}{(1-2 x s)^2}\right) (1-2 x s)^{-a} \,
   ds
\end{multline}
\end{example}
\begin{example}
From Eq. (10.65.1) in \cite{dlmf}. Spherical Bessel function of the first kind.
\begin{multline}
\sum _{n=1}^{\infty } \frac{j_{n-1}(z) t^n}{(n+1)!}
=\frac{1}{t z^4}z^2 (\cos (z)+z \sin (z))+\cosh ^{-1}(\cos (z))^2\\
 \left(t z \cos (z)+\cos \left(\sqrt{z (-2 t+z)}\right)+\sqrt{z (-2
   t+z)} \sin \left(\sqrt{z (-2 t+z)}\right)\right)
\end{multline}
\end{example}
\begin{example}
From \cite{frontczak}.
\begin{multline}
\sum _{n=0}^{\infty } \frac{F_n(x) z^n}{(n+1)!}=-\frac{1}{z}+\frac{e^{\frac{x z}{2}} \left(\sqrt{4+x^2} \cosh \left(\frac{1}{2} \sqrt{4+x^2} z\right)-x \sinh \left(\frac{1}{2}
   \sqrt{4+x^2} z\right)\right)}{\sqrt{4+x^2} z}\
\end{multline}
\end{example}
\begin{example}
From \cite{medina}
\begin{multline}
\sum _{n=1}^{\infty } \frac{z^n H_{j n}}{(n+1)!}=\int_0^z \left(\frac{e^s (\gamma +\Gamma (0,s)+\log (s))}{j z}+\frac{\left(e^s s\right) \sum _{k=1}^{j-1} \frac{\,
   _2F_2\left(1,1;2,1+\frac{k}{j};-s\right)}{k}}{z}\right) \, ds
\end{multline}
\end{example}
\begin{example}
From Eq. (5.6.1.1) in \cite{prud2}
\begin{multline}
\sum _{n=0}^{\infty } \frac{H_n(x) t^n}{(n+b) \left\lfloor
   \frac{n}{2}\right\rfloor !}=-t^{-b} \int_0^t \frac{s^{-1+b}
   \left(-e^{\frac{4 s^2 x^2}{1+4 s^2}}-4 e^{\frac{4 s^2 x^2}{1+4 s^2}} s^2-2
   e^{\frac{4 s^2 x^2}{1+4 s^2}} s x\right)}{\left(1+4 s^2\right)^{3/2}} \,
   ds
\end{multline}
\end{example}
\section{Definite integrals involving logarithm, quotient Pochammer and q-Pochammer functions}
\begin{example}
\begin{multline}
\int_0^{\infty } \frac{x^m \log ^k(a x)}{(1+x)_{1+n}} \, dx=\sum
   _{p=0}^n \frac{i (-1)^p i^k e^{i m \pi } (1+p)^m (2 \pi )^{1+k} \binom{n}{p}
   \Phi \left(e^{2 i m \pi },-k,\frac{\pi -i \log (a (1+p))}{2 \pi
   }\right)}{n!}
\end{multline}
\end{example}
\begin{example}
\begin{multline}
\int_0^{\infty } \frac{\sqrt{x} \log (\log (a x))}{(1+x)_{1+n}} \, dx=\frac{\pi  }{n!}\sum _{p=0}^n (-1)^p \sqrt{1+p}
   \binom{n}{p} \log \left(\frac{\Gamma \left(\frac{\pi -i \log (a (1+p))}{4 \pi }\right)^2}{4 \pi  i \Gamma
   \left(\frac{3 \pi -i \log (a (1+p))}{4 \pi }\right)^2}\right)
\end{multline}
\end{example}
\begin{example}
\begin{multline}
\int_0^{\infty } \frac{x^{-m} \left(1-x^{2 m}\right)}{\log (x) (1+x)_{1+n}} \, dx\\
=\frac{e^{-i m \pi }}{n!}\sum _{p=0}^n (-1)^p
    (1+p)^{-m} \binom{n}{p} \left(\Phi \left(e^{-2 i m \pi },1,\frac{\pi -i \log (1+p)}{2 \pi
   }\right)\right. \\ \left.
-e^{2 i m \pi } (1+p)^{2 m} \Phi \left(e^{2 i m \pi },1,\frac{\pi -i \log (1+p)}{2 \pi
   }\right)\right)
\end{multline}
\end{example}
\begin{example}
\begin{multline}
\int_0^{\infty } \frac{\log (\log (a x))}{(1+x)_{1+n}} \, dx=\frac{2 \pi  }{i n!}\sum _{p=0}^n (-1)^p \binom{n}{p}
   \log \left(\frac{(a (1+p))^{\frac{\pi -2 i \log (2 \pi )}{2 (2 \pi )}} \Gamma \left(\frac{\pi -i \log (a (1+p))}{2
   \pi }\right)}{\sqrt{2 \pi }}\right)
\end{multline}
\end{example}
\begin{example}
\begin{multline}
\int_0^{\infty } \frac{\sqrt[4]{x} \log (\log (a x))}{(1+x)_{1+n}} \, dx
=-\frac{\pi }{n!}\left(\frac{1}{2}+\frac{i}{2}\right)\sum _{p=0}^n(-1)^{\frac{1}{4}+p} \sqrt[4]{1+p}  \binom{n}{p}\\
 \left(\pi -2 i \log
   (2 \pi )+(2+2 i) \log \left(\frac{\Gamma \left(\frac{\pi -i \log (a (1+p))}{8 \pi }\right)}{2 \Gamma
   \left(\frac{5}{8}-\frac{i \log (a (1+p))}{8 \pi }\right)}\right)\right. \\ \left.
-(2-2 i) \log \left(\frac{\Gamma
   \left(\frac{3}{8}-\frac{i \log (a (1+p))}{8 \pi }\right)}{2 \Gamma \left(\frac{7}{8}-\frac{i \log (a (1+p))}{8 \pi
   }\right)}\right)\right)
\end{multline}
\end{example}
\begin{example}
\begin{multline}
\int_0^{\infty } \frac{x^m}{\left(a^2+b^2-2 i a \log (x)-\log ^2(x)\right) (1+x)_{1+n}} \, dx\\
=\sum _{p=0}^n
   \frac{(-1)^p e^{i m \pi } (1+p)^m \binom{n}{p} \left(\Phi \left(e^{2 i m \pi },1,\frac{a-i b+\pi -i \log (1+p)}{2
   \pi }\right)-\Phi \left(e^{2 i m \pi },1,\frac{a+i b+\pi -i \log (1+p)}{2 \pi }\right)\right)}{2 b n!}
\end{multline}
\end{example}
\begin{example}
\begin{multline}
\int_0^{\infty } \frac{x^m \log ^k(a x)}{(1+b x)_{1+n}} \, dx\\
=-\sum _{p=0}^n \frac{(-1)^{1+p} b^{-1-m} e^{i m
   \pi } (1+p)^m (2 i \pi )^{1+k} \binom{n}{p} \Phi \left(e^{2 i m \pi },-k,\frac{\pi -i \log \left(\frac{a
   (1+p)}{b}\right)}{2 \pi }\right)}{n!}
\end{multline}
\end{example}
\begin{example}
\begin{multline}
\int_0^{\infty } \frac{\log (\log (a x))}{(1+b x)_{1+n}} \, dx\\
=\log \left(\prod _{p=0}^n (1+p)^{-\frac{(-1)^p
   \binom{n}{p} \log (2 i \pi )}{b n!}} \left(2 \pi  \Gamma \left(\frac{\pi -i \log \left(\frac{a (1+p)}{b}\right)}{2
   \pi }\right)\right)^{-\frac{2 i (-1)^p \pi  \binom{n}{p}}{b n!}}\right)
\end{multline}
\end{example}
\begin{example}
\begin{multline}
\int_0^{\infty } \frac{x^m (-x)^{(q)} \log ^k(a x)}{(1+x)_{1+n}} \, dx\\
=-\sum _{p=0}^n \frac{i (-1)^{1+p} i^k
   e^{i m \pi } (1+p)^m (2 \pi )^{1+k} \binom{n}{p} (1+p)^{(q)} \Phi \left(e^{2 i m \pi },-k,\frac{\pi -i \log (a
   (1+p))}{2 \pi }\right)}{n!}
\end{multline}
\end{example}
\begin{example}
\begin{multline}
\int_0^{\infty } x^m\frac{ (x)_q }{(1+x)_{1+n}} \log ^k(a x)\, dx\\
=\sum _{p=0}^n \frac{i (-1)^{1+p} i^k e^{i m
   \pi } (1+p)^m (2 \pi )^{1+k} \binom{n}{p} (1+p)^{(q)} \Phi \left(e^{2 i m \pi },-k,\frac{\pi -i \log (a (1+p))}{2
   \pi }\right)}{n!}
\end{multline}
\end{example}
\begin{example}
\begin{multline}
\int_0^{\infty } \frac{x^m (-x)^{(q)} \log ^k(a x)}{(x+\beta )_n} \, dx\\
=-\sum _{p=0}^{n-1} \frac{i (-1)^{1+p}
   i^k e^{i m \pi } (2 \pi )^{1+k} (p+\beta )^m \binom{-1+n}{p} (p+\beta )^{(q)} \Phi \left(e^{2 i m \pi },-k,\frac{\pi
   -i \log (a (p+\beta ))}{2 \pi }\right)}{(-1+n)!}
\end{multline}
where $n>q$
\end{example}
\begin{example}
\begin{multline}
\int_0^{\infty } \frac{\left(x^m-x^{-m}\right) (x)_q}{\log (x) (x+\beta )_n} \, dx\\
=\sum _{p=0}^n \frac{(-1)^p
   e^{-i m \pi } (p+\beta )^{-m} \binom{-1+n}{p} (p+\beta )^{(q)} }{(-1+n)!}\\
\left(\Phi \left(e^{-2 i m \pi },1,\frac{\pi -i \log
   (p+\beta )}{2 \pi }\right)-e^{2 i m \pi } (p+\beta )^{2 m} \Phi \left(e^{2 i m \pi },1,\frac{\pi -i \log (p+\beta
   )}{2 \pi }\right)\right)
\end{multline}
\end{example}
\begin{example}
\begin{multline}
\int_0^{\infty } \frac{(1-x) (x)_q}{\sqrt{x} \log (a x) (x+\beta )_n} \, dx\\
=\sum _{p=0}^n \frac{i (-1)^p
   (1+p+\beta ) \binom{-1+n}{p} (p+\beta )^{(q)} \left(\psi ^{(0)}\left(\frac{\pi -i \log (a (p+\beta ))}{4 \pi
   }\right)-\psi ^{(0)}\left(\frac{3}{4}-\frac{i \log (a (p+\beta ))}{4 \pi }\right)\right)}{2 \sqrt{p+\beta }
   (-1+n)!}
\end{multline}
\end{example}
\begin{example}
\begin{multline}
\int_0^{\infty } \frac{x^s (x;q)_p}{(1+b+x)_{1+n}} \, dx=\sum _{j=0}^p \sum _{k=0}^n \frac{(-1)^{1+j+k}
   (1+b+k)^{j+s} \pi  q^{\frac{1}{2} (-1+j) j} \binom{n}{k} \csc (\pi  (j+s)) \binom{p}{j}_q}{n!}
\end{multline}
\end{example}
\begin{example}
\begin{multline}
\int_0^{\infty } \frac{x^m \log ^k(a x) (x;q)_p}{(1+b+x)_{1+n}} \, dx\\
=\sum
   _{j=0}^p \sum _{l=0}^n \frac{(-1)^{j+l} e^{i (j+m) \pi } (1+b+l)^{j+m} (2 i
   \pi )^{1+k} q^{\frac{1}{2} (-1+j) j} \binom{n}{l} }{n!}\\
\Phi \left(e^{2 i (j+m) \pi
   },-k,\frac{\pi -i \log (a (1+b+l))}{2 \pi }\right) \binom{p}{j}_q
\end{multline}
where $n>p$
\end{example}
\begin{example}
\begin{equation}
\int_0^{\infty } \frac{x^t \left(x^s;q\right){}_n}{\left(x^s;q\right){}_{n+1}} \, dx=\frac{\pi 
   \left(-q^n\right)^{-\frac{1+t}{s}} \csc \left(\frac{\pi  (1+t)}{s}\right)}{s}
\end{equation}
\end{example}
\begin{example}
\begin{multline}
\int_0^{\infty } \frac{x^m \log ^k(a x) \left(x^{\alpha };q\right){}_n}{\left(x^{\alpha };q\right){}_{1+n}} \,
   dx\\
=\frac{e^{\frac{i (1+m) \pi }{\alpha }} (2 \pi )^{1+k} \left(-q^n\right)^{-\frac{1+m}{\alpha }}
   \left(\frac{i}{\alpha }\right)^{-1+k} \Phi \left(e^{\frac{2 i (1+m) \pi }{\alpha }},-k,\frac{\pi -i \alpha  \log
   \left(a \left(-q^n\right)^{-1/\alpha }\right)}{2 \pi }\right)}{\alpha ^2}
\end{multline}
\end{example}
\begin{example}
\begin{multline}
\int_0^{\infty } \frac{x^{-1-m+r} \log ^k\left(\frac{a}{x}\right)}{\left(1+b x^r\right)^n} \, dx\\
=\sum _{j=0}^{n-1} \sum
   _{l=0}^j \frac{(-1)^{-j+n} b^{-1+\frac{m}{r}} e^{\frac{i m \pi }{r}} m^l (2 i \pi )^{1-j+k+l} r^{-1-k-l} \binom{j}{l} \Gamma
   (1+k)}{\Gamma (1-j+k+l) \Gamma (n)}\\
 \Phi \left(e^{\frac{2 i m \pi }{r}},j-k-l,\frac{\pi -i r \log \left(a b^{1/r}\right)}{2 \pi }\right)
   S_{-1+n}^{(j)}
\end{multline}
\end{example}
\begin{example}
\begin{multline}
\int_0^1 (1-x)^m \left(\frac{x}{1-x}\right)^p \, dx\\
=\sum _{j=0}^m \sum _{q=0}^j \frac{(-1)^m (-m+p) \pi  \csc (p
   \pi ) (-1)^{m-j} S_m^{(j)} \binom{j}{q} (1-m)^{j-q} p^q}{(1+m)!}
\end{multline}
\end{example}
\begin{example}
\begin{multline}
\int_0^1 (1-x)^r \left(\frac{x}{1-x}\right)^m \log ^k\left(\frac{a x}{1-x}\right) \, dx\\
=\sum _{j=0}^r \sum
   _{q=0}^j \frac{k! (-1)^{-j+2 r} e^{i m \pi } (2 i \pi )^{k-q} (1+m-r)^{j-q} \binom{j}{q}}{\Gamma (1+k-q) \Gamma (2+r)}\\
 \left(-2 i \pi  (m-r) \Phi
   \left(e^{2 i m \pi },-k+q,\frac{\pi -i \log (a)}{2 \pi }\right)\right. \\ \left.
+(-k+q) \Phi \left(e^{2 i m \pi },1-k+q,\frac{\pi -i
   \log (a)}{2 \pi }\right)\right) S_r^{(j)}
\end{multline}
\end{example}
\begin{example}
\begin{multline}
\int_0^1 \left(\sum _{k=1}^p q^{\frac{1}{2} (-1+k) k} x^n \left(-1+x^{-1-n}\right)^{k w} \binom{p}{k}_q\right)
   \, dx=\sum _{k=1}^p \frac{k \pi  q^{\frac{1}{2} (-1+k) k} w \csc (k \pi  w) \binom{p}{k}_q}{1+n}
\end{multline}
\end{example}
\begin{example}
\begin{multline}
\int_0^1 \left(\sum _{j=1}^p q^{\frac{1}{2} (-1+j) j} x^n \left(-1+x^{-1-n}\right)^{j m} \log ^k\left(a
   \left(-1+x^{-1-n}\right)^j\right) \binom{p}{j}_q\right) \, dx\\
=-\sum _{j=1}^p \frac{e^{i j m \pi } (i j)^k (2 \pi )^k
   q^{\frac{1}{2} (-1+j) j} }{1+n}\\
\left(k \Phi \left(e^{2 i j m \pi },1-k,\frac{1}{2}-\frac{i \log (a)}{2 j \pi }\right)+2 i
   j m \pi  \Phi \left(e^{2 i j m \pi },-k,\frac{1}{2}-\frac{i \log (a)}{2 j \pi }\right)\right)
   \binom{p}{j}_q
\end{multline}
\end{example}
\begin{example}
\begin{multline}
\int_0^1 x^n \left(-1+x^{-1-n}\right)^w \left(-\left(-1+x^{-1-n}\right)^w;q\right){}_p \, dx\\
=\sum _{k=0}^p
   \frac{(1+k) \pi  q^{\frac{1}{2} (-1+k) k} w \csc ((1+k) \pi  w) \binom{p}{k}_q}{1+n}
\end{multline}
\end{example}
\begin{example}
\begin{multline}
\int_0^{\infty } x^{-1+m} \log ^{-1+k}(a x) (k+m \log (a x)) \log \left(x^s
   \left(1+\frac{1}{x}\right)_s\right) \, dx\\
=\sum _{t=1}^s -e^{i m \pi } (2 i \pi )^{1+k} t^{-m} \Phi \left(e^{2
   i m \pi },-k,\frac{\pi -i \log \left(\frac{a}{t}\right)}{2 \pi }\right)
\end{multline}where $Re(a)<0,Re(m)<0$
\end{example}
\begin{example}
Logarithm of the $q$-Pochhammer symbol.
\begin{multline}
\int_0^{\infty } x^{-1+m} \log ^{-1+k}(a x) (k+m \log (a x)) \log \left(\frac{(-b;x)_{1+n}}{1+b}\right)
   \, dx\\
=-\sum _{j=1}^n i b^{-\frac{m}{j}} e^{\frac{i m \pi }{j}} \left(\frac{i}{j}\right)^k (2 \pi )^{1+k} \Phi
   \left(e^{\frac{2 i m \pi }{j}},-k,\frac{\pi -i j \log \left(a b^{-1/j}\right)}{2 \pi }\right)
\end{multline}
where $Re(a)<0,Re(m)<0$
\end{example}
\begin{example}
\begin{multline}
\int_0^{\infty } x^{-1+m} \log ^{-1+k}(a x) (k+m \log (a x)) \log \left(\frac{(1+c) (-b;x)_{1+n}}{(1+b)
   (-c;x)_{1+n}}\right) \, dx\\
=\sum _{j=1}^n i b^{-\frac{m}{j}} c^{-\frac{m}{j}} e^{\frac{i m \pi }{j}}
   \left(\frac{i}{j}\right)^k (2 \pi )^{1+k} \left(-c^{m/j} \Phi \left(e^{\frac{2 i m \pi }{j}},-k,\frac{\pi -i j \log
   \left(a b^{-1/j}\right)}{2 \pi }\right)\right. \\ \left.
+b^{m/j} \Phi \left(e^{\frac{2 i m \pi }{j}},-k,\frac{\pi -i j \log \left(a
   c^{-1/j}\right)}{2 \pi }\right)\right)
\end{multline}
\end{example}
\begin{example}
\begin{multline}
\int_0^{\infty } \frac{\log ^k\left(\frac{a}{x}\right) (x)_n}{1+b x^r} \, dx\\
=-\sum _{m=0}^n(-1)^{-m+n} i^{1+k}
   b^{-\frac{1+m}{r}} e^{\frac{i \pi  (-1-m+r)}{r}} (2 \pi )^{1+k} r^{-1-k}\\
 \Phi \left(e^{\frac{2 i \pi 
   (-1-m+r)}{r}},-k,\frac{\pi -i r \log \left(a b^{1/r}\right)}{2 \pi }\right) S_n^{(m)}
\end{multline}
\end{example}
\begin{example}
\begin{multline}
\int_0^{\infty } \frac{\log ^k\left(\frac{a}{x}\right) (x;q)_n}{1+b x^r} \, dx\\
=-\sum _{m=0}^n (-1)^m i^{1+k}
   b^{-\frac{1+m}{r}} e^{\frac{i \pi  (-1-m+r)}{r}} (2 \pi )^{1+k} q^{\frac{1}{2} (-1+m) m} r^{-1-k}\\
 \Phi
   \left(e^{\frac{2 i \pi  (-1-m+r)}{r}},-k,\frac{\pi -i r \log \left(a b^{1/r}\right)}{2 \pi }\right)
   \binom{n}{m}_q
\end{multline}
\end{example}
\begin{example}
From Eq. (2.2.4.8) in \cite{prud1}.
\begin{multline}
\int_0^z \frac{\log ^k\left(a x \left(-x^u+z^u\right)^{-1/u}\right) \left(x
   \left(-x^u+z^u\right)^{-1/u}\right){}_n}{x} \, dx\\
=-\sum _{m=0}^n (-1)^{-m+n} e^{\frac{i m \pi }{u}} (2 \pi )^{1+k}
   \left(\frac{i}{u}\right)^{1+k} \Phi \left(e^{\frac{2 i m \pi }{u}},-k,\frac{\pi -i u \log (a)}{2 \pi }\right)
   S_n^{(m)}
\end{multline}
where $Re(k)<0,Re(a)<0$
\end{example}
\section{The Polylogarithm function}
In this section we develop a table of generating functions involving Hurwitz-Lerch zeta function in terms of definite integrals of the polylogarithm function at low-order. In this section we used the binomial coefficient $\binom{n}{j} $, Pochhammer symbol $ (n)_j$. We used Eq. (2.5.2.1) in \cite{prud3} and the method in \cite{reyn4}.
\begin{example}
\begin{multline}
\int_0^{\infty } x^{-1+m} \log ^k(a x) \text{Li}_n(-c x) \, dx\\
=\sum _{j=0}^{\infty } \frac{(-1)^{1-n} i^{1+j}
   c^{-m} e^{\frac{1}{2} i (k+2 m) \pi } m^{-j-n} (2 \pi )^{1-j+k} \Phi \left(e^{2 i m \pi },j-k,\frac{\pi -i \log
   \left(\frac{a}{c}\right)}{2 \pi }\right) (1-j+k)_j (n)_j}{\Gamma (1+j)}
\end{multline}
\end{example}
\begin{example}
\begin{multline}
\int_0^{\infty } x^{-1+m} \log ^k(-c x) \text{Li}_n(-c x) \, dx\\
=-\sum _{j=0}^{\infty } \frac{c^{-m} e^{-i m \pi
   } (-m)^{-j-n} (2 i \pi )^{1-j+k} \binom{-1+j+n}{j} k! \text{Li}_{j-k}\left(e^{2 i m \pi }\right)}{(-j+k)!}
\end{multline}
\end{example}
\begin{example}
\begin{multline}
\int_0^{\infty } x^{-1+m} \log ^k(-c x) \text{Li}_n(-c x) \, dx\\
=-\sum _{j=0}^{\infty } c^{-m} e^{-i (j+m) \pi }
   (-m)^{-j-n} \binom{-1+j+n}{j} \Gamma (1+k)\\
 \left(i^{2 j} \zeta (1-j+k,-m)-i^{2 k} \zeta (1-j+k,1+m)\right)
\end{multline}
\end{example}
\begin{example}
\begin{multline}
\int_0^{\infty } \frac{\log ^k(a x) \text{Li}_n(-c x)}{\sqrt{x}} \, dx\\
=\sum _{j=0}^{\infty } \frac{(-1)^{-j-n}
   i^{-j+k} 2^{1-j+2 k+n} \pi ^{1-j+k} \binom{-1+j+n}{j} k! }{\sqrt{c} (-j+k)!}\\
\left(\zeta \left(j-k,\frac{\pi -i \log
   \left(\frac{a}{c}\right)}{4 \pi }\right)-\zeta \left(j-k,\frac{3}{4}-\frac{i \log \left(\frac{a}{c}\right)}{4 \pi
   }\right)\right)
\end{multline}
\end{example}
\begin{example}
\begin{multline}
\int_0^{\infty } \frac{\log ^k(a x) \text{Li}_n(-c x)}{x^{3/2}} \, dx\\
=\sum _{j=0}^{\infty } \frac{i^{-j+k}
   2^{1-j+2 k+n} \sqrt{c} \pi ^{1-j+k} \binom{-1+j+n}{j} k! \left(-\zeta \left(j-k,\frac{\pi -i \log
   \left(\frac{a}{c}\right)}{4 \pi }\right)+\zeta \left(j-k,\frac{3}{4}-\frac{i \log \left(\frac{a}{c}\right)}{4 \pi
   }\right)\right)}{(-j+k)!}
\end{multline}
\end{example}
\begin{example}
\begin{multline}
\int_0^{\infty } \frac{\log ^k(-c x) \text{Li}_n(-c x)}{\sqrt{x}} \, dx\\
=-\sum _{j=0}^{\infty } \frac{i
   (-1)^{-j-n} 2^{1+k+n} \left(1-2^{1-j+k}\right) (i \pi )^{1-j+k} \binom{-1+j+n}{j} k! \zeta (j-k)}{\sqrt{c}
   (-j+k)!}
\end{multline}
\end{example}
\begin{example}
\begin{multline}
\int_0^{\infty } \frac{\log ^k(-c x) \text{Li}_n(-c x)}{x^{3/2}} \, dx\\
=\sum _{j=0}^{\infty } \frac{i^{2-j+k}
   2^{1+k+n} \left(1-2^{1-j+k}\right) \sqrt{c} \pi ^{1-j+k} \binom{-1+j+n}{j} k! \zeta (j-k)}{(-j+k)!}
\end{multline}
\end{example}
\section{Logarithmic function generating functions}
\begin{example}
\begin{multline}
\int_0^{\infty } \frac{x^m \sqrt{\log (x)} \log (\log (x))}{1+c x} \, dx\\
=-\sum _{j=0}^{\infty } \frac{i^j
   2^{\frac{1}{2}-j} c^{-1-m} e^{\frac{1}{4} i (1+4 m) \pi } m^{-j} \pi ^{\frac{3}{2}-j} (0)_j
   \left(\frac{3}{2}-j\right)_j }{\Gamma (1+j)}\\
\left(\Phi \left(e^{2 i m \pi },-\frac{1}{2}+j,\frac{\pi -i \log
   \left(\frac{1}{c}\right)}{2 \pi }\right) \left(\pi +2 i H_{\frac{1}{2}-j}+i (-4+\log (4)-2 \log (\pi ))\right)\right. \\ \left.
+2 i
   \Phi'\left(e^{2 i m \pi },-\frac{1}{2}+j,\frac{\pi -i \log \left(\frac{1}{c}\right)}{2
   \pi }\right)\right)
\end{multline}
where $Re(m)<0$
\end{example}
\begin{example}
\begin{multline}
\int_0^{\infty } \frac{x^m \log (\log (x))}{(1+c x) \sqrt{\log (x)}} \, dx\\
=-\sum _{j=0}^{\infty } \frac{i^j
   2^{-\frac{1}{2}-j} c^{-1-m} e^{\frac{1}{4} i (-1+4 m) \pi } m^{-j} \pi ^{\frac{1}{2}-j} (0)_j
   \left(\frac{1}{2}-j\right)_j }{\Gamma (1+j)}\\
\left(\Phi \left(e^{2 i m \pi },\frac{1}{2}+j,\frac{\pi -i \log
   \left(\frac{1}{c}\right)}{2 \pi }\right) \left(\pi +2 i H_{-\frac{1}{2}-j}+2 i \log \left(\frac{2}{\pi
   }\right)\right)\right. \\ \left.
+2 i \Phi'\left(e^{2 i m \pi },\frac{1}{2}+j,\frac{\pi -i \log
   \left(\frac{1}{c}\right)}{2 \pi }\right)\right)
\end{multline}
\end{example}
\begin{example}
\begin{multline}
\int_0^{\infty } \frac{x^m \log ^k(x)}{1-c^2 x^2} \, dx\\
=\sum _{j=0}^{\infty } \frac{i^{1+j} 2^{-j+k} (-c)^{-m}
   c^{-1-m} e^{\frac{1}{2} i (k+2 m) \pi } m^{-j} \pi ^{1-j+k}}{\Gamma (1+j)}\\
 \left(-c^m \Phi \left(e^{2 i m \pi },j-k,\frac{\pi -i
   \log \left(-\frac{1}{c}\right)}{2 \pi }\right)+(-c)^m \Phi \left(e^{2 i m \pi },j-k,\frac{\pi -i \log
   \left(\frac{1}{c}\right)}{2 \pi }\right)\right)\\
 (0)_j (1-j+k)_j
\end{multline}
where $Im(c)\neq 0$
\end{example}
\begin{example}
\begin{multline}
\sum _{j=0}^{\infty } \sum _{p=0}^{\infty } \sum _{q=0}^{2 p+1} \frac{(-1)^j e^{m+2 j m} m^{1+2 p-q}
   \binom{1+2 p}{q} \log ^{k-q}\left(a e^{1+2 j}\right)}{(1+2 p)! (k-q)!}\\
=\frac{2^{1+k} e^{2 m} \Phi
   \left(-e^{2 m},-k,1+\frac{\log (a)}{2}\right)-\log ^k(a)}{2 k!}
\end{multline}
\end{example}
\begin{example}
From Eq. (6.3.4.2(i)) in \cite{prud3} and \cite{reyn4}.
\begin{multline}
\sum _{j=0}^{\infty } \sum _{p=0}^{2 j} \frac{4^j m^{2 j-p} \binom{2 j}{p} E_{2 j}(x) \log ^{k-p}(a)}{j!
   (k-p)!}\\
=\frac{2^k e^{m-m (-1+2 x)} \left(\Phi \left(-e^{2 m},-k,\frac{1}{2} (2-2 x+\log (a))\right)+e^{2 m (-1+2 x)}
   \Phi \left(-e^{2 m},-k,\frac{1}{2} (2 x+\log (a))\right)\right)}{k!}
\end{multline}
\end{example}
\begin{example}
From Eq. (5.3.1.5(ii)) in \cite{prud2} and \cite{reyn4}.
\begin{multline}
\sum _{p=0}^{\infty } \sum _{q=0}^{2 p} \frac{\left(-m^2\right)^p (1+k-q)_q (1+2 p-q)_q \zeta (2 p)}{(2 \pi  a
   m)^q \Gamma (1+q)}\\
=\frac{1}{4} a^{-1-k} \left(a^k (k+2 a m \pi )+2 a e^{2 m \pi } k \Phi \left(e^{2 m \pi
   },1-k,1+a\right)\right. \\ \left.
+4 a e^{2 m \pi } m \pi  \Phi \left(e^{2 m \pi },-k,1+a\right)\right)
\end{multline}
\end{example}
\begin{example}
From Eq. (5.11.4.2) in \cite{prud2}.
\begin{multline}
\sum _{p=0}^{\infty } \sum _{j=0}^p \frac{m^{-j+p} t^{-j} \binom{p}{j} L_n^{-n+p}(x)}{(-j+k)! p!}=\sum _{j=0}^n
   \frac{e^m \left(1+\frac{1}{t}\right)^k (1+t)^{-j} (m-x)^{-j+n} \binom{n}{j}}{(-j+k)! n!}
\end{multline}
\end{example}
\begin{example}
From Eq. (5.2.9.5) in \cite{prud1} and Eq. (4.9.2.12) in \cite{brychkov}.
\begin{equation}
\int_0^a \frac{a^{3/2} e^{b \sqrt{(a-x) x}} I_1\left(b \sqrt{(a-x) x}\right)}{2 \sqrt{a-x} x} \, dx=\sum
   _{k=1}^{\infty } \frac{a \left(a b e^{-a b}\right)^k k^k}{(1+k)!}=-a-\frac{1}{b}+\frac{e^{a b}}{b}
\end{equation}
\end{example}
\begin{example}
From Eq. (2.6.6.4) in \cite{prud1}.
\begin{multline}
\int_0^1 \frac{x^{-1+m} \log ^{\sigma }\left(\frac{1}{x}\right) \log ^k(a x)}{b^2+x^2-2 b x \cos (\gamma )}
   \, dx\\
=\frac{\csc (\gamma ) \Gamma (1+k) \Gamma (1+\sigma ) \log ^k(a) }{b}\sum _{j=0}^{\infty } \sum
   _{p=1}^{\infty } \frac{b^{-p} (-1+m+p)^{-1-j-\sigma } \binom{-1-\sigma }{j} \log ^{-j}(a) \sin (p \gamma
   )}{(-j+k)!}
\end{multline}
where $Re(m)>1,Re(b)>1,Re(\sigma)>-1/2,Re(a)>e^{3\pi}$
\end{example}
\begin{example}
From Eq. (2.6.10.3) in \cite{prud1}.
\begin{multline}
\int_0^1 \frac{(1+x)^{-1+m} \log ^n(1+x) \log ^k(a (1+x))}{k!} \, dx\\
=-\sum _{j=0}^{\infty } \frac{(-1)^n
   m^{-1-j-n} \binom{-1-n}{j} n! \log ^{-j+k}(a)}{(-j+k)!}\\
+\sum _{j=0}^{\infty } \sum _{p=0}^n \frac{(-1)^p 2^m
   m^{-1-j-p} \binom{-1-p}{j} n! \log ^{n-p}(2) \log ^{-j+k}(2 a)}{(-j+k)! (n-p)!}
\end{multline}
where $Re(m)>1,Re(a)>e^{3\pi}$
\end{example}
\section{A few definite integrals}
\begin{example}
\begin{multline}
 \int_0^{\infty } \frac{\log \left(\frac{(1+x) \left(b^2+x\right)}{(b+x)^2}\right) (2+\log (a x) \log (\log
   (a x)))}{2 \sqrt{x} \log (a x)} \, dx\\
=\frac{1}{4} \pi  \left(2 i \left(-1+\sqrt{b}\right)^2 \pi +24 (1+b) \log
   (2)+4 \left(-1+\sqrt{b}\right)^2 \log (\pi )-8 \log (-3 \pi -i \log (a))\right. \\ \left.
+8 \log \left(\frac{1}{4} (-\pi -i
   \log (a))\right)+16 \sqrt{b} (\log (-3 \pi -i \log (a)-i \log (b))\right. \\ \left.
-\log (-2 (\pi +i \log (a)+i \log (b))))-8 b
   \log (-3 \pi -i \log (a)-2 i \log (b))\right. \\ \left.
+8 b \log \left(\frac{1}{4} (-\pi -i \log (a)-2 i \log (b))\right)+8
   \text{log$\Gamma $}\left(-\frac{\pi +i \log (a)}{4 \pi }\right)\right. \\ \left.
-8 \left(\text{log$\Gamma
   $}\left(-\frac{3}{4}-\frac{i \log (a)}{4 \pi }\right)+2 \sqrt{b} \text{log$\Gamma $}\left(-\frac{\pi +i \log
   (a)+i \log (b)}{4 \pi }\right)\right.\right. \\ \left.\left.
-2 \sqrt{b} \text{log$\Gamma $}\left(-\frac{3 \pi +i \log (a)+i \log (b)}{4 \pi
   }\right)-b \text{log$\Gamma $}\left(-\frac{\pi +i \log (a)+2 i \log (b)}{4 \pi }\right)\right.\right. \\ \left.\left.
+b \text{log$\Gamma
   $}\left(-\frac{3 \pi +i \log (a)+2 i \log (b)}{4 \pi }\right)\right)\right)
\end{multline}
\end{example}
\begin{example}
\begin{multline}
 \int_0^{\infty } \frac{\log \left(\frac{(1+x) (16+x)}{(4+x)^2}\right) (2+\log (-x) \log (\log (-x)))}{2
   \sqrt{x} \log (-x)} \, dx\\
=\frac{1}{2} \pi  \left(i \pi +\log (16)-16 \log \left(\frac{(\pi +i \log (4))
   \left(\Gamma \left(-\frac{\pi +i \log (4)}{2 \pi }\right) \Gamma \left(-\frac{i \log (4)}{4 \pi
   }\right)\right)}{(2 \pi +i \log (4)) \left(\Gamma \left(-\frac{i \log (4)}{2 \pi }\right) \Gamma
   \left(-\frac{1}{2}-\frac{i \log (4)}{4 \pi }\right)\right)}\right)\right)
\end{multline}
\end{example}
\begin{example}
\begin{multline}
\int_0^{\infty } \frac{\log \left(\frac{(1+x) \left(b^2+x\right)}{(b+x)^2}\right) (4+\log (a x) \log (\log
   (a x)))}{4 x^{3/4} \log (a x)} \, dx\\
=\frac{\pi }{\sqrt{2}} \left(\log (16)+\left(-1+\sqrt[4]{b}\right)^2 \left(i \pi
   +\log \left(4 \pi ^2\right)\right)-(2-2 i) \log \left(\frac{\Gamma \left(\frac{\pi -i \log (a)}{8 \pi
   }\right)}{\Gamma \left(\frac{5}{8}-\frac{i \log (a)}{8 \pi }\right)}\right)\right. \\ \left.
-(2+2 i) \log \left(\frac{\Gamma
   \left(\frac{3}{8}-\frac{i \log (a)}{8 \pi }\right)}{\Gamma \left(\frac{7}{8}-\frac{i \log (a)}{8 \pi
   }\right)}\right)+(4+4 i) \sqrt[4]{b} \left(-i \log \left(\frac{\Gamma \left(\frac{\pi -i \log (a)-i \log
   (b)}{8 \pi }\right)}{2 \Gamma \left(\frac{5 \pi -i \log (a)-i \log (b)}{8 \pi }\right)}\right)\right.\right. \\ \left.\left.
+\log
   \left(\frac{\Gamma \left(\frac{3 \pi -i \log (a)-i \log (b)}{8 \pi }\right)}{2 \Gamma \left(\frac{7 \pi -i
   \log (a)-i \log (b)}{8 \pi }\right)}\right)\right)-(2-2 i) \sqrt{b} \left(\log \left(\frac{\Gamma
   \left(\frac{\pi -i \log (a)-2 i \log (b)}{8 \pi }\right)}{2 \Gamma \left(\frac{5 \pi -i \log (a)-2 i \log
   (b)}{8 \pi }\right)}\right)\right.\right. \\ \left.\left.
+i \log \left(\frac{\Gamma \left(\frac{3 \pi -i \log (a)-2 i \log (b)}{8 \pi
   }\right)}{2 \Gamma \left(\frac{7 \pi -i \log (a)-2 i \log (b)}{8 \pi}\right)}\right)\right)\right)
\end{multline}
\end{example}
\begin{example}
\begin{multline}
\int_0^{\infty } \frac{\log (\log (x))}{(b+x) \left(b^2+x\right)} \, dx\\
=\frac{i \left(\log (b) (\pi -2 i
   \log (2 \pi ))+4 \pi  \log \left(\frac{-\pi -2 i \log (b)}{-\pi -i \log (b)}\right)+4 \pi  \log
   \left(\frac{\Gamma \left(-\frac{1}{2}-\frac{i \log (b)}{\pi }\right)}{\Gamma \left(-\frac{\pi +i \log (b)}{2
   \pi }\right)}\right)\right)}{2 (-1+b) b}
\end{multline}
\end{example}
\begin{example}
\begin{multline}
\int_0^{\infty } \frac{\log (\log (x))}{\left(e^{i b \pi }+x\right) \left(e^{2 i b \pi }+x\right)} \,
   dx\\
=\frac{1}{4} e^{-\frac{3 i b \pi }{2} } \pi  \csc \left(\frac{b \pi }{2}\right) \left(i b \pi +4 \log
   \left(2+\frac{1}{-1+b}\right)+2 b \log (2 \pi )+4 \log \left(\frac{\Gamma \left(-\frac{1}{2}+b\right)}{\Gamma
   \left(\frac{1}{2} (-1+b)\right)}\right)\right)
\end{multline}
\end{example}
\begin{example}
\begin{equation}
\int_0^{\infty } \frac{x \log (\log (x))}{1+x^2+x^4} \, dx=\frac{\pi  }{6 \sqrt{3}}\log
   \left(\frac{\left(-\frac{27}{25}\right)^3 \pi ^8}{\left(\Gamma \left(-\frac{5}{3}\right) \Gamma
   \left(\frac{1}{3}\right)\right)^6}\right)
\end{equation}
\end{example}
\begin{example}
\begin{multline}
\int_0^{\infty } \frac{x^{-1-m} \left(1+m \log (x)+x^{2 m} (-1+m \log (x))\right) \log \left(\frac{(1+x)
   \left(b^2+x\right)}{(b+x)^2}\right)}{\log ^2(x)} \, dx\\
=2 \left(\coth ^{-1}\left(e^{i m \pi }\right)-\tanh
   ^{-1}\left(e^{i m \pi }\right)\right)+b^{-2 m} e^{-i m \pi } \left(-2 b^m \Phi \left(e^{-2 i m \pi
   },1,\frac{\pi -i \log (b)}{2 \pi }\right)\right. \\ \left.
+\Phi \left(e^{-2 i m \pi },1,\frac{1}{2}-\frac{i \log (b)}{\pi
   }\right)+b^{3 m} e^{2 i m \pi } \left(2 \Phi \left(e^{2 i m \pi },1,\frac{\pi -i \log (b)}{2 \pi }\right)\right.\right. \\ \left.\left.
-b^m
   \Phi \left(e^{2 i m \pi },1,\frac{1}{2}-\frac{i \log (b)}{\pi }\right)\right)\right)
\end{multline}
\end{example}
\begin{example}
From Eq. (2.6.15.13) in \cite{prud1}
\begin{multline}
\int_0^{\infty } x^{-1+m} \log ^{-1+k}(a x) (k+m \log (a x)) \log \left(\frac{c^2+(x+\beta
   )^2}{c^2+(b+x)^2}\right) \, dx\\
=(2 i \pi )^{1+k} \left(\left(b^2+c^2\right)^{m/2} e^{i m \left(\pi -\tan
   ^{-1}\left(\frac{c}{b}\right)\right)} \left(e^{2 i m \tan ^{-1}\left(\frac{c}{b}\right)}\right.\right. \\ \left.\left.
 \Phi \left(e^{2 i m \pi
   },-k,\frac{\pi +\tan ^{-1}\left(\frac{c}{b}\right)-i \log (a)-\frac{1}{2} i \log \left(b^2+c^2\right)}{2 \pi
   }\right)\right.\right. \\ \left.\left.
+\Phi \left(e^{2 i m \pi },-k,-\frac{-2 \pi +2 \tan ^{-1}\left(\frac{c}{b}\right)+2 i \log (a)+i \log
   \left(b^2+c^2\right)}{4 \pi }\right)\right)\right. \\ \left.
-e^{i m \left(\pi -\tan ^{-1}\left(\frac{c}{\beta }\right)\right)}
   \left(c^2+\beta ^2\right)^{m/2} \left(e^{2 i m \tan ^{-1}\left(\frac{c}{\beta }\right)}\right.\right. \\ \left.\left. \Phi \left(e^{2 i m \pi
   },-k,\frac{\pi +\tan ^{-1}\left(\frac{c}{\beta }\right)-i \log (a)-\frac{1}{2} i \log \left(c^2+\beta ^2\right)}{2
   \pi }\right)\right.\right. \\ \left.\left.
+\Phi \left(e^{2 i m \pi },-k,-\frac{-2 \pi +2 \tan ^{-1}\left(\frac{c}{\beta }\right)+2 i \log (a)+i
   \log \left(c^2+\beta ^2\right)}{4 \pi }\right)\right)\right)
\end{multline}
where $Re(a)<0$
\end{example}
\begin{example}
\begin{multline}
\int_0^{\infty } \frac{\left(k+\frac{1}{2} \log (a x)\right) \log ^{-1+k}(a x) \log \left(\frac{c^2+(x+\beta
   )^2}{c^2+(b+x)^2}\right)}{\sqrt{x}} \, dx\\
=i^{2+k} 2^{1+2 k} \pi ^{1+k} \left(\sqrt[4]{b^2+c^2} e^{-\frac{1}{2} i
   \tan ^{-1}\left(\frac{c}{b}\right)} \left(e^{i \tan ^{-1}\left(\frac{c}{b}\right)}\right.\right. \\ \left.\left.
    \left(\zeta \left(-k,\frac{\pi
   +\tan ^{-1}\left(\frac{c}{b}\right)-i \log (a)-\frac{1}{2} i \log \left(b^2+c^2\right)}{4 \pi }\right)\right.\right.\right. \\ \left.\left.\left.
-\zeta
   \left(-k,\frac{6 \pi +2 \tan ^{-1}\left(\frac{c}{b}\right)-2 i \log (a)-i \log \left(b^2+c^2\right)}{8 \pi
   }\right)\right)\right.\right. \\ \left.\left.
-\zeta \left(-k,-\frac{-6 \pi +2 \tan ^{-1}\left(\frac{c}{b}\right)+2 i \log (a)+i \log
   \left(b^2+c^2\right)}{8 \pi }\right)\right.\right. \\ \left.\left.
+\zeta \left(-k,-\frac{-2 \pi +2 \tan ^{-1}\left(\frac{c}{b}\right)+2 i \log
   (a)+i \log \left(b^2+c^2\right)}{8 \pi }\right)\right)\right. \\ \left.
+e^{-\frac{1}{2} i \tan ^{-1}\left(\frac{c}{\beta }\right)}
   \sqrt[4]{c^2+\beta ^2} \left(e^{i \tan ^{-1}\left(\frac{c}{\beta }\right)} \left(-\zeta \left(-k,\frac{\pi +\tan
   ^{-1}\left(\frac{c}{\beta }\right)-i \log (a)-\frac{1}{2} i \log \left(c^2+\beta ^2\right)}{4 \pi }\right)\right.\right.\right. \\ \left.\left.\left.
+\zeta
   \left(-k,\frac{6 \pi +2 \tan ^{-1}\left(\frac{c}{\beta }\right)-2 i \log (a)-i \log \left(c^2+\beta ^2\right)}{8
   \pi }\right)\right)\right.\right. \\ \left.\left.
+\zeta \left(-k,-\frac{-6 \pi +2 \tan ^{-1}\left(\frac{c}{\beta }\right)+2 i \log (a)+i \log
   \left(c^2+\beta ^2\right)}{8 \pi }\right)\right.\right. \\ \left.\left.
-\zeta \left(-k,-\frac{-2 \pi +2 \tan ^{-1}\left(\frac{c}{\beta
   }\right)+2 i \log (a)+i \log \left(c^2+\beta ^2\right)}{8 \pi }\right)\right)\right)
\end{multline}
\end{example}
\begin{example}
\begin{multline}
\int_0^{\infty } \frac{(-2+\log (a x)) \log \left(\frac{c^2+(x+\beta )^2}{c^2+(b+x)^2}\right)}{\sqrt{x} \log
   ^2(a x)} \, dx\\
=\frac{i }{\sqrt{2}}\left(-\frac{\left(b+i c+\sqrt{b^2+c^2}\right) \psi ^{(0)}\left(\frac{2 \pi +2 \tan
   ^{-1}\left(\frac{c}{b}\right)-2 i \log (a)-i \log \left(b^2+c^2\right)}{8 \pi }\right)}{\sqrt[4]{b^2+c^2}
   \sqrt{1+\frac{b}{\sqrt{b^2+c^2}}}}\right. \\ \left.
+\frac{\left(b+i c+\sqrt{b^2+c^2}\right) \psi ^{(0)}\left(\frac{6 \pi +2 \tan
   ^{-1}\left(\frac{c}{b}\right)-2 i \log (a)-i \log \left(b^2+c^2\right)}{8 \pi }\right)}{\sqrt[4]{b^2+c^2}
   \sqrt{1+\frac{b}{\sqrt{b^2+c^2}}}}\right. \\ \left.
+\frac{\sqrt[4]{b^2+c^2} \left(b+i c+\sqrt{b^2+c^2}\right) \psi
   ^{(0)}\left(-\frac{-6 \pi +2 \tan ^{-1}\left(\frac{c}{b}\right)+2 i \log (a)+i \log \left(b^2+c^2\right)}{8 \pi
   }\right)}{(b+i c) \sqrt{1+\frac{b}{\sqrt{b^2+c^2}}}}\right. \\ \left.
-\frac{\sqrt[4]{b^2+c^2} \left(b+i c+\sqrt{b^2+c^2}\right)
   \psi ^{(0)}\left(-\frac{-2 \pi +2 \tan ^{-1}\left(\frac{c}{b}\right)+2 i \log (a)+i \log \left(b^2+c^2\right)}{8
   \pi }\right)}{(b+i c) \sqrt{1+\frac{b}{\sqrt{b^2+c^2}}}}\right. \\ \left.
+\frac{\left(i c+\beta +\sqrt{c^2+\beta ^2}\right) \psi
   ^{(0)}\left(\frac{2 \pi +2 \tan ^{-1}\left(\frac{c}{\beta }\right)-2 i \log (a)-i \log \left(c^2+\beta
   ^2\right)}{8 \pi }\right)}{\sqrt[4]{c^2+\beta ^2} \sqrt{1+\frac{\beta }{\sqrt{c^2+\beta ^2}}}}\right. \\ \left.
-\frac{\left(ic+\beta +\sqrt{c^2+\beta ^2}\right) \psi ^{(0)}\left(\frac{6 \pi +2 \tan ^{-1}\left(\frac{c}{\beta }\right)-2 i
   \log (a)-i \log \left(c^2+\beta ^2\right)}{8 \pi }\right)}{\sqrt[4]{c^2+\beta ^2} \sqrt{1+\frac{\beta
   }{\sqrt{c^2+\beta ^2}}}}\right. \\ \left.
-\frac{\sqrt[4]{c^2+\beta ^2} \left(i c+\beta +\sqrt{c^2+\beta ^2}\right) \psi
   ^{(0)}\left(-\frac{-6 \pi +2 \tan ^{-1}\left(\frac{c}{\beta }\right)+2 i \log (a)+i \log \left(c^2+\beta
   ^2\right)}{8 \pi }\right)}{(i c+\beta ) \sqrt{1+\frac{\beta }{\sqrt{c^2+\beta ^2}}}}\right. \\ \left.
+\frac{\sqrt[4]{c^2+\beta ^2}
   \left(c-i \left(\beta +\sqrt{c^2+\beta ^2}\right)\right) \psi ^{(0)}\left(-\frac{-2 \pi +2 \tan
   ^{-1}\left(\frac{c}{\beta }\right)+2 i \log (a)+i \log \left(c^2+\beta ^2\right)}{8 \pi }\right)}{(c-i \beta )
   \sqrt{1+\frac{\beta }{\sqrt{c^2+\beta ^2}}}}\right)
\end{multline}
\end{example}
\begin{example}
\begin{multline}
\int_0^{\infty } \frac{(-2+\log (i x)) \log \left(\frac{\frac{5}{4}+x+x^2}{1+(1+x)^2}\right)}{\sqrt{x} \log
   ^2(i x)} \, dx\\
=\frac{e^{-\frac{1}{2} i \tan ^{-1}(2)}}{10^{3/4}} \left(5 i \sqrt[4]{2} \psi ^{(0)}\left(\frac{2 \pi +\tan
   ^{-1}\left(\frac{4}{3}\right)-i \log \left(\frac{5}{4}\right)}{8 \pi }\right)\right. \\ \left.
-5 i \sqrt[4]{2} \psi
   ^{(0)}\left(\frac{6 \pi +\tan ^{-1}\left(\frac{4}{3}\right)-i \log \left(\frac{5}{4}\right)}{8 \pi
   }\right)+\sqrt{5} \left(-2 \sqrt[4]{-4-3 i} H_{-\frac{\pi +2 i \log (2)}{16 \pi }}-2 i \sqrt[4]{4+3 i}
   H_{-\frac{11}{16}-\frac{i \log (2)}{8 \pi }}\right.\right. \\ \left.\left.
+2 \sqrt[4]{-4-3 i} H_{-\frac{9}{16}-\frac{i \log (2)}{8 \pi }}+2 i
   \sqrt[4]{4+3 i} H_{-\frac{3}{16}-\frac{i \log (2)}{8 \pi }}\right.\right. \\ \left.\left.
-(2-i) \sqrt[4]{2} \left(\psi ^{(0)}\left(-\frac{-4
   \pi +\tan ^{-1}\left(\frac{4}{3}\right)+i \log \left(\frac{5}{4}\right)}{8 \pi }\right)-\psi ^{(0)}\left(\frac{15
   \pi +\tan ^{-1}\left(\frac{24}{7}\right)-i \log \left(\frac{25}{16}\right)}{16 \pi
   }\right)\right)\right)\right)
\end{multline}
\end{example}
\begin{example}
\begin{multline}
\int_0^{\infty } \frac{\log \left(\frac{1+(1+x)^2}{1+x^2}\right) (2+\log (-x) \log (\log (-x)))}{\sqrt{x}
   \log (-x)} \, dx\\
=\sqrt[8]{-1} \pi  \left(\left((-1-i) \sqrt[8]{-1}+\frac{1+i}{\sqrt[4]{2}}+5 i \sqrt[4]{2}\right)
   \pi -(2+2 i) 2^{3/4} \tan ^{-1}\left(\frac{\log (4)}{9 \pi }\right)\right. \\ \left.
-4 i \sqrt[4]{2} \tan ^{-1}\left(\frac{\log
   (4)}{7 \pi }\right)-4 i \sqrt[4]{2} \tan ^{-1}\left(\frac{\log (4)}{\pi }\right)+(2+2 i) 2^{3/4} \tan
   ^{-1}\left(\frac{\log (4)}{\pi }\right)\right. \\ \left.
+\sqrt[8]{-1} \log \left(\frac{4096}{2401}\right)+(4+4 i) (-1)^{5/8} \log
   (2)+(2-2 i) 2^{3/4} \log (2)\right. \\ \left.
+(-1)^{5/8} \log \left(\frac{4096}{625}\right)+\sqrt[4]{2} \log (16)-\sqrt[8]{-1}
   \log \left(\frac{4096}{81}\right)+\sqrt[8]{-1} \log (4096)-(-1)^{5/8} \log (4096)\right. \\ \left.
+(2+2 i) (-1)^{5/8} \log (\pi
   )+2 \sqrt[4]{2} \log (\pi )+(1-i) 2^{3/4} \log (\pi )+2 \sqrt[4]{2} \log \left(\pi ^2+4 \log ^2(2)\right)\right. \\ \left.
+(1-i)
   2^{3/4} \log \left(\pi ^2+4 \log ^2(2)\right)-2 \sqrt[4]{2} \log \left(49 \pi ^2+4 \log ^2(2)\right)-(1-i)
   2^{3/4} \log \left(81 \pi ^2+4 \log ^2(2)\right)\right. \\ \left.
-4 \sqrt[8]{-1} \text{log$\Gamma $}\left(-\frac{7}{8}\right)-4
   (-1)^{5/8} \text{log$\Gamma $}\left(-\frac{5}{8}\right)+4 \sqrt[8]{-1} \text{log$\Gamma
   $}\left(-\frac{3}{8}\right)+4 (-1)^{5/8} \text{log$\Gamma $}\left(-\frac{1}{8}\right)\right. \\ \left.
+4 \sqrt[4]{2}
   \text{log$\Gamma $}\left(\frac{\pi -2 i \log (2)}{16 \pi }\right)+(2-2 i) 2^{3/4} \text{log$\Gamma
   $}\left(-\frac{\pi +2 i \log (2)}{16 \pi }\right)\right. \\ \left.
-(2-2 i) 2^{3/4} \text{log$\Gamma $}\left(-\frac{9}{16}-\frac{i
   \log (2)}{8 \pi }\right)-4 \sqrt[4]{2} \text{log$\Gamma $}\left(-\frac{7}{16}-\frac{i \log (2)}{8 \pi
   }\right)\right)
\end{multline}
\end{example}
\begin{example}

\end{example}
\section{Extended Abel theorem forms}
\begin{example}
\begin{multline}
\sum _{m=0}^{\infty } \frac{(a+b m)^m \left(e^{-b y} y\right)^m}{(m+3)!}
\\=\frac{e^{3 b y}}{y^3} \left(\frac{1}{(-a+3 b)^3}+\frac{(-a+b) e^{-b y} y-\frac{1}{2} (a-2 b)^2 e^{-2 b y} y^2}{(a-2 b)^2
   (a-b)}\right. \\ \left.
+\frac{e^{(a-3 b) y} \left((a-2 b)^2 (a-b)-(2 a-5 b) (a-3 b) (a-b) b y+(a-3 b)^2 (a-2 b) b^2 y^2\right)}{(a-3 b)^3 (a-2 b)^2 (a-b)}\right)
\end{multline}
\end{example}
\begin{example}
\begin{multline}
\sum _{m=0}^{\infty } \frac{(a+b m)^m \left(e^{-b y} y\right)^m}{(m+4)!}
=-\frac{e^{4 b y}}{(a-4 b)^4 y^4}-\frac{e^{3 b y}}{(a-3 b)^3 y^3}-\frac{e^{2 b y}}{2 (a-2 b)^2 y^2}-\frac{e^{b y}}{6 (a-b)
   y}\\
-\frac{e^{a y}}{(a-4 b)^4
   (a-3 b)^3 (a-2 b)^2 (a-b) y^4} \left(-a^6+14 a^5 b-80 a^4 b^2+238 a^3 b^3-387 a^2 b^4\right. \\ \left.
+324 a b^5-108 b^6+3 a^6 b y-48 a^5 b^2 y+310 a^4 b^3 y-1029 a^3 b^4 y+1840 a^2 b^5 y-1668 a b^6 y\right. \\ \left.
+592 b^7 y-3 a^6 b^2 y^2+54 a^5 b^3
   y^2-395 a^4 b^4 y^2+1494 a^3 b^5 y^2-3054 a^2 b^6 y^2+3152 a b^7 y^2\right. \\ \left.
-1248 b^8 y^2+a^6 b^3 y^3-20 a^5 b^4 y^3+165 a^4 b^5 y^3-718 a^3 b^6 y^3+1736 a^2 b^7 y^3\right. \\ \left.
-2208 a b^8 y^3+1152 b^9 y^3\right)
\end{multline}
\end{example}
\begin{example}
\begin{multline}
\sum _{m=0}^{\infty } \frac{(a+b m)^m \left(e^{-b y} y\right)^m}{(m+5)!}\\
=\frac{e^{5 b y}}{(-a+5 b)^5 y^5}-\frac{e^{4 b y}}{(-a+4 b)^4 y^4}+\frac{e^{3 b y}}{2 (-a+3 b)^3 y^3}-\frac{e^{2 b y}}{6 (-a+2
   b)^2 y^2}+\frac{e^{b y}}{24 (-a+b) y}\\
-\frac{e^{(a-5 b) y+5 b y}}{(-a+b) (-a+2 b)^2 (-a+3 b)^3 (-a+4 b)^4 (-a+5 b)^5 y^5} \left(a^{10}-30 a^9 b+400 a^8 b^2\right. \\ \left.
-3118 a^7 b^3+15715 a^6 b^4-53428 a^5 b^5+123852 a^4 b^6-192832 a^3 b^7+192384 a^2 b^8\right. \\ \left.
-110592 a b^9+27648
   b^{10}-4 a^{10} b y+130 a^9 b^2 y-1870 a^8 b^3 y+15657 a^7 b^4 y-84371 a^6 b^5 y\right. \\ \left.
+305216 a^5 b^6 y-749108 a^4 b^7 y+1228581 a^3 b^8 y-1284363 a^2 b^9 y+769392 a b^{10} y\right. \\ \left.
-199260 b^{11} y+6 a^{10} b^2
   y^2-210 a^9 b^3 y^2+3255 a^8 b^4 y^2-29373 a^7 b^5 y^2+170548 a^6 b^6 y^2\right. \\ \left.
-664156 a^5 b^7 y^2+1751627 a^4 b^8 y^2-3078241 a^3 b^9 y^2+3434084 a^2 b^{10} y^2-2183140 a b^{11} y^2\right. \\ \left.
+595600 b^{12} y^2-4 a^{10}
   b^3 y^3+150 a^9 b^4 y^3-2500 a^8 b^5 y^3+24352 a^7 b^6 y^3-153258 a^6 b^7 y^3\right. \\ \left.
+649686 a^5 b^8 y^3-1873208 a^4 b^9 y^3+3612612 a^3 b^{10} y^3-4434230 a^2 b^{11} y^3+3100400 a b^{12} y^3\right. \\ \left.
-924000 b^{13}
   y^3+a^{10} b^4 y^4-40 a^9 b^5 y^4+715 a^8 b^6 y^4-7518 a^7 b^7 y^4+51471 a^6 b^8 y^4\right. \\ \left.
-239628 a^5 b^9 y^4+767837 a^4 b^{10} y^4-1670990 a^3 b^{11} y^4+2361800 a^2 b^{12} y^4-1956000 a b^{13} y^4\right. \\ \left.
+720000
   b^{14} y^4\right)
\end{multline}
\end{example}
\begin{example}
\begin{multline}
\frac{e^{a x^p}}{1-b p x^p}-1\\
=\sum _{n=1}^{\infty } \frac{(a+b n)^{n/p}
   \left(e^{-b x^p} x\right)^n \left((-1)^n+(-1)^{n+p}+2 \csc \left(\frac{n \pi
   }{p}\right) \sin \left(\frac{n \pi  \left(-1+2 \left\lceil
   \frac{p}{2}\right\rceil \right)}{p}\right)\right)}{2 n \Gamma
   \left(\frac{n}{p}\right)}
\end{multline}
where $p$ is any positive integer.
\end{example}
\begin{example}
\begin{multline}
\frac{2 b e^{a x^p} p x^p}{-1+b^2 p^2 x^{2 p}}
=\sum _{n=1}^{\infty } \frac{\left(-(a+b n)^{n/p} \left(e^{-b
   x^p} x\right)^n+(a-b n)^{n/p} \left(e^{b x^p} x\right)^n\right) }{2 n
   \Gamma \left(\frac{n}{p}\right)}\\
\left((-1)^n+(-1)^{n+p}+2 \csc \left(\frac{n \pi
   }{p}\right) \sin \left(\frac{n \pi  \left(-1+2 \left\lceil \frac{p}{2}\right\rceil \right)}{p}\right)\right)
\end{multline}
\end{example}
\begin{example}
\begin{multline}
\sum _{n=1}^{\infty } \frac{(a+b n)^{n/p}
   \left(x^p\right)^{-\frac{c}{p}} \left(-e^{b c \left(-1+x^p\right)}
   \left(e^{-b}\right)^n+\left(e^{-b x^p} x\right)^n
   \left(x^p\right)^{c/p}\right) }{n (c+n) \Gamma
   \left(\frac{n}{p}\right)}\\\left((-1)^n+(-1)^{n+p}+2 \csc \left(\frac{n
   \pi }{p}\right) \sin \left(\frac{n \pi  \left(-1+2 \left\lceil
   \frac{p}{2}\right\rceil \right)}{p}\right)\right)
\\
=\frac{e^{-b c} \left(x^p\right)^{-\frac{c}{p}}}{c p}
   \left(\left(e^{b c x^p}-e^{b c} \left(x^p\right)^{c/p}\right) \left(-1+e^{i
   p \pi }+4 \left\lceil \frac{p}{2}\right\rceil \right)\right. \\ \left.+2 c e^{b c
   \left(1+x^p\right)} \left(E_{1-\frac{c}{p}}(-a+b c)-\left(x^p\right)^{c/p}
   E_{1-\frac{c}{p}}\left((-a+b c) x^p\right)\right)\right)
\end{multline}
\end{example}
\begin{example}
\begin{multline}
\int_{-1}^1 \frac{e^{a x^2}}{1-2 b x^2} \, dx\\
=\sum _{n=1}^{\infty }
   \frac{\left(1+(-1)^n\right) \left(1+e^{i n \pi }\right) (b
   n)^{-\frac{1}{2}-\frac{n}{2}} (a+b n)^{n/2} \left(\Gamma
   \left(\frac{1+n}{2}\right)-e^{b n} \left(e^{-b}\right)^n \Gamma
   \left(\frac{1+n}{2},b n\right)\right)}{2 n \Gamma
   \left(\frac{n}{2}\right)}+2
\end{multline}
\end{example}
\begin{example}
\begin{multline}
\int_{-1}^1 \frac{e^{\left(\frac{x}{\pi }\right)^2}}{6-x^2} \,
   dx\\
=\frac{1}{6}
   \left(2+\sum _{n=1}^{\infty } \frac{\sqrt{3} \left(1+e^{i n \pi }\right)^2
   n^{-\frac{3}{2}-\frac{n}{2}} \pi ^{-n} \left(12+n \pi ^2\right)^{n/2}
   \left(\Gamma \left(\frac{1+n}{2}\right)-\Gamma
   \left(\frac{1+n}{2},\frac{n}{12}\right)\right)}{\Gamma
   \left(\frac{n}{2}\right)}\right)
\end{multline}
\end{example}
\begin{example}
\begin{multline}
\frac{4 \tanh ^{-1}(b x) \sinh (a x)}{x}\\
=\sum _{n=1}^{\infty } \frac{1}{n^2 \Gamma (n)}\left(e^{(a-b n) x} (a-b n)^{1+n} \left(e^{b x}
   x\right)^n E_{-n}((a-b n) x)\right. \\ \left.
+(a+b n)^{1+n} \left(e^{-((a+b n) x)} \left(-e^{b x} x\right)^n E_{-n}(-((a+b n)
   x))-e^{(a+b n) x} \left(e^{-b x} x\right)^n E_{-n}((a+b n) x)\right)\right. \\ \left.
-e^{-a x+b n x} (a-b n)^{1+n} \left(-e^{-b x}
   x\right)^n E_{-n}(-a x+b n x)\right)
\end{multline}
\end{example}
\begin{example}
\begin{multline}
\frac{4 \tanh ^{-1}\left(3 b x^3\right) \sinh \left(a x^3\right)}{3 x^3}\\
=\sum _{n=1}^{\infty }
   \frac{e^{-\left((a+b n) x^3\right)} \left(1+2 \cos \left(\frac{2 n \pi }{3}\right)\right)}{n^2 \Gamma
   \left(\frac{n}{3}\right)} \left((a-b
   n)^{1+\frac{n}{3}} \left(e^{2 a x^3} \left(e^{b x^3} x\right)^n E_{-\frac{n}{3}}\left((a-b n) x^3\right)\right.\right. \\ \left.\left.
-e^{2 b n
   x^3} \left(-e^{-b x^3} x\right)^n E_{-\frac{n}{3}}\left((-a+b n) x^3\right)\right)+(a+b n)^{1+\frac{n}{3}}
   \left(\left(-e^{b x^3} x\right)^n E_{-\frac{n}{3}}\left(-\left((a+b n) x^3\right)\right)\right.\right. \\ \left.\left.
-e^{2 (a+b n) x^3}
   \left(e^{-b x^3} x\right)^n E_{-\frac{n}{3}}\left((a+b n) x^3\right)\right)\right)
\end{multline}
\end{example}
\begin{example}
\begin{multline}
\frac{2 \left(-e^{a (-x)^p} (-x)^{-p} \tanh ^{-1}\left(b p (-x)^p\right)+e^{a x^p} x^{-p} \tanh ^{-1}\left(b p
   x^p\right)\right)}{p}\\
=\sum _{n=1}^{\infty } \frac{1}{2
   n^2 \Gamma \left(\frac{n}{p}\right)}\left(-e^{(a-b n) (-x)^p} (a-b n)^{\frac{n+p}{p}} \left(-e^{b
   (-x)^p} x\right)^n E_{-\frac{n}{p}}\left((a-b n) (-x)^p\right)\right. \\ \left.
+e^{(a-b n) x^p} (a-b n)^{\frac{n+p}{p}} \left(e^{b
   x^p} x\right)^n E_{-\frac{n}{p}}\left((a-b n) x^p\right)\right. \\ \left.
+(a+b n)^{\frac{n+p}{p}} \left(e^{(a+b n) (-x)^p}
   \left(-e^{-b (-x)^p} x\right)^n E_{-\frac{n}{p}}\left((a+b n) (-x)^p\right)\right.\right. \\ \left.\left.
-e^{(a+b n) x^p} \left(e^{-b x^p}
   x\right)^n E_{-\frac{n}{p}}\left((a+b n) x^p\right)\right)\right)\\
 \left((-1)^n+(-1)^{n+p}+2 \csc \left(\frac{n \pi
   }{p}\right) \sin \left(\frac{n \pi  \left(-1+2 \left\lceil \frac{p}{2}\right\rceil \right)}{p}\right)\right)
\end{multline}
\end{example}
\begin{example}
\begin{equation}
-\frac{e^{a x} (1+x)}{-1+b x (1+x)}=\sum _{n=0}^{\infty } \left(\frac{e^{-b x} x}{1+x}\right)^n L_n(-a-b
   n)
\end{equation}
\end{example}
\begin{example}
\begin{equation}
-\frac{e^{a x} (-1+x)}{1+b (-1+x) x}=\sum _{n=0}^{\infty } \left(\frac{e^{-b x} x}{-1+x}\right)^n L_n(a+b
   n)
\end{equation}
\end{example}
\begin{example}
\begin{multline}
-\frac{e^{(a+b) x} (-1+x)}{1+(-1-b+c) x+b x^2}\\
=\sum _{n=0}^{\infty } \frac{(-1)^{(-1+c) n} \left(e^{-b x}
   (-1+x)^{-c} x\right)^n \Gamma (1+c n) \, _1F_1(-n;1+(-1+c) n;a+b+b n)}{n! \Gamma (1+(-1+c) n)}
\end{multline}
\end{example}
\begin{example}
\begin{multline}
\frac{(1+x) \log \left(1+a x+x^2\right)}{1+x-b x}
=\sum _{n=0}^{\infty }
   \left(\frac{\left(a-\sqrt{-4+a^2}\right) \left(x (1+x)^{-b}\right)^n \Gamma (1+b n) }{2 \Gamma (n) \Gamma (2-n+b
   n)}\right. \\ \left.
\, _3F_2\left(1,1,1-n;2,2-n+b
   n;\frac{1}{2} \left(a-\sqrt{-4+a^2}\right)\right)\right. \\ \left.
+\frac{\left(a+\sqrt{-4+a^2}\right) \left(x (1+x)^{-b}\right)^n \Gamma (1+b n) }{2 \Gamma (n) \Gamma (2-n+b n)}\right. \\ \left.\, _3F_2\left(1,1,1-n;2,2-n+b
   n;\frac{1}{2} \left(a+\sqrt{-4+a^2}\right)\right)\right)
\end{multline}
where $Im(b)\neq 0$
\end{example}
\begin{example}
\begin{multline}
\frac{(1+x) \log \left(1+x+x^2+x^3\right)}{1+x-a x}\\
=\sum _{n=1}^{\infty } \left(-\frac{i \left(x
   (1+x)^{-a}\right)^n \Gamma (1+a n) \left(\, _3F_2(1,1,1-n;2,2-n+a n;-i)-\, _3F_2(1,1,1-n;2,2-n+a
   n;i)\right)}{\Gamma (n) \Gamma (2-n+a n)}\right. \\ \left.
   +\frac{a \left(x (1+x)^{-a}\right)^n \Gamma (a n) \left(-H_{(-1+a)
   n}+H_{a n}\right)}{\Gamma (n) \Gamma (1+(-1+a) n)}\right)
\end{multline}
where $Im(a)\neq 0$
\end{example}
\begin{example}
\begin{multline}
-\frac{(-1+x) \log (1+(a-x) x)}{1+(-1+b) x}\\
=\sum _{n=1}^{\infty } \left(-\frac{(-1)^n
   \left(a+\sqrt{4+a^2}\right) \left((1-x)^{-b} x\right)^n \Gamma (1+b n) }{2 \Gamma (n) \Gamma (2-n+b n)}\right. \\ \left.
   \, _3F_2\left(1,1,1-n;2,2-n+b n;\frac{1}{2}
   \left(-a-\sqrt{4+a^2}\right)\right)\right. \\ \left.
+\frac{(-1)^n \left(\sqrt{2-i a} \sqrt{2+i
   a}-a\right) \left((1-x)^{-b} x\right)^n \Gamma (1+b n)}{2 \Gamma (n) \Gamma (2-n+b n)}\right)\\ \, _3F_2\left(1,1,1-n;2,2-n+b n;\frac{1}{2}
   \left(-a+\sqrt{4+a^2}\right)\right)
\end{multline}
where $Im(b)\neq 0$
\end{example}
\begin{example}
\begin{multline}
\frac{(-1+x) \log \left(1-x^2+2 x \sinh (a)\right)}{-1+x-b x}\\
=\sum _{n=1}^{\infty } -\frac{(-1)^n e^{-a}
   \left((1-x)^{-b} x\right)^n \Gamma (1+b n) }{\Gamma (n) \Gamma (2-n+b n)}\\
\left(-\, _3F_2\left(1,1,1-n;2,2-n+b n;e^{-a}\right)+e^{2 a} \,
   _3F_2\left(1,1,1-n;2,2-n+b n;-e^a\right)\right)
\end{multline}
where $Im(b)\neq 0$
\end{example}
\begin{example}
\begin{multline}
\frac{e^{a x} (1-x)^{1+d}}{1+(-1-b+c) x+b x^2}\\
=\sum _{n=0}^{\infty }
   \frac{(-1)^n \left(e^{-b x} (1-x)^{-c} x\right)^n \Gamma (1+d+c n) \,
   _1F_1(-n;1+d+(-1+c) n;a+b n)}{n! \Gamma (1+d+(-1+c) n)}
\end{multline}
\end{example}
\begin{example}
\begin{multline}
-\frac{(-1+x) x^m (-1+x z)}{-1+x-b x+x (1-a+(-1+a+b) x) z}\\
=\sum
   _{p=0}^{\infty } \frac{(-1)^{m+p} \left((1-x)^{-b} x (1-x
   z)^{-a}\right)^{m+p} \Gamma (1+b (m+p))}{\Gamma (1+p) \Gamma (1-p+b (m+p))}\\
 \, _2F_1(-p,-a (m+p);1+m+(-1+b)
   (m+p);z)
\end{multline}
\end{example}
\begin{example}
\begin{multline}
\frac{e^{(a+b m) x} \cosh (x)}{1-b x}=\sum _{p=0}^{\infty }
   \frac{\left((-1+a+b (m+p))^p+(1+a+b (m+p))^p\right) \left(e^{-b x}
   x\right)^p}{2 \Gamma (1+p)}
\end{multline}
\end{example}
\begin{example}
From Eq. (1.9.54) in \cite{erd_t1}
\begin{multline}
\int_0^{\infty } e^{-i m x} \left((a \pi -i x)^k+e^{2 i m x} (a \pi +i x)^k\right) \text{sech}(x) \sin
   \left(\frac{x^2}{\pi }\right) \, dx\\
=-\sqrt{2} e^{\frac{m \pi }{2}} \pi ^{1+k} \Phi \left(-e^{m \pi
   },-k,\frac{1}{2}+a\right)\\
+\sum _{p=0}^{\infty } \sum _{j=0}^{4 p} \frac{2 (-1)^p e^{\frac{m \pi }{2}} m^{-j+4 p}
   \pi ^{\frac{1}{2}-j+k+2 p} (1-j+k)_j (1-j+4 p)_{-\frac{1}{2}+j-2 p} }{\Gamma (1+j)}\\
\Phi \left(-e^{m \pi
   },j-k,\frac{1}{2}+a\right)
\end{multline}
\end{example}
\begin{example}
From Eq. (1.9.58) in \cite{erd_t1}
\begin{multline}
\int_0^{\infty } \frac{e^{-i m x} \left((-i x+\log (a))^k+e^{2 i m x} (i x+\log (a))^k\right)
   \text{sech}\left(\sqrt{\frac{\pi }{2}} x\right) \left(\cos \left(\frac{x^2}{2}\right)+\sin \left(\frac{x^2}{2}\right)\right)}{2
   k!} \, dx\\
=\sum _{p=0}^{\infty } \sum _{j=0}^{4 p} \frac{\sqrt{\pi } (-1)^p 2^{\frac{1}{2}-j+\frac{j-k}{2}+k-2 p} e^{m
   \sqrt{\frac{\pi }{2}}} m^{-j+4 p} \pi ^{\frac{1}{2} (-j+k)} \binom{4 p}{j} }{(-j+k)! (2 p)!}\\
\Phi \left(-e^{m \sqrt{2 \pi
   }},j-k,\frac{\sqrt{\frac{\pi }{2}}+\log (a)}{\sqrt{2 \pi }}\right)\\
+\sum _{p=0}^{\infty } \sum _{j=0}^{2 (1+2
   p)} \frac{\sqrt{\pi } (-1)^p 2^{-\frac{1}{2}-j+\frac{j-k}{2}+k-2 p} e^{m \sqrt{\frac{\pi }{2}}} m^{-j+2 (1+2 p)} \pi
   ^{\frac{1}{2} (-j+k)} \binom{2 (1+2 p)}{j} }{(-j+k)! (1+2 p)!}\\
\Phi \left(-e^{m \sqrt{2 \pi }},j-k,\frac{\sqrt{\frac{\pi }{2}}+\log (a)}{\sqrt{2\pi }}\right)
\end{multline}
\end{example}
\begin{example}
From (1.9.54) and  (1.9.55) in \cite{erd_t1}
\begin{multline}
\int_0^{\infty } e^{-\frac{i x (m \pi +x)}{\pi }} \left((a \pi -i x)^k+e^{2 i m x} (a \pi +i x)^k\right) \text{sech}(x) \,
   dx\\
=(1+i) \sqrt{2} e^{\frac{m \pi }{2}} \pi ^{1+k} \Phi \left(-e^{m \pi },-k,\frac{1}{2}+a\right)\\
-\sum _{j=0}^{\infty } \sum
   _{p=0}^{2 j} \frac{2 i e^{\frac{1}{2} (i j+m) \pi } m^{2 j-p} \pi ^{\frac{1}{2}+j+k-p} \Phi \left(-e^{m \pi
   },-k+p,\frac{1}{2}+a\right) (1+2 j-p)_{-\frac{1}{2}-j+p} (1+k-p)_p}{\Gamma (1+p)}
\end{multline}
\end{example}
\begin{example}
From (1.14.5) in \cite{erd_t1}
\begin{multline}
\int_0^{\infty } \frac{e^{-i m x} \, _1F_2\left(\alpha ;\frac{1}{2},\beta ;-\frac{x^2}{4}\right) \left((-i x+\log
   (a))^k+e^{2 i m x} (i x+\log (a))^k\right)}{2 k!} \, dx\\
=\sum _{j=0}^{\infty } \sum _{p=0}^{\infty } \frac{(-1)^j
   m^{-1+2 j-p+2 \alpha } \pi  \binom{-1+2 j+2 \alpha }{p} \binom{-1-\alpha +\beta }{j} \Gamma (\beta ) \log
   ^{k-p}(a)}{(k-p)! \Gamma (\alpha ) \Gamma (-\alpha +\beta )}
\end{multline}
\end{example}
\begin{example}
\begin{multline}
\int_0^{\infty } i e^{-i m x} x^{-v} \left((-i x+\log (a))^k-e^{2 i m x} (i x+\log (a))^k\right) \, dx\\
=\sum
   _{j=0}^{\infty } 2 (1-j+k)_j m^{-1-j+v} \binom{-1+v}{j} \cos \left(\frac{\pi  v}{2}\right) \Gamma (1-v) \log
   ^{-j+k}(a)
\end{multline}
where $Re(a)<0,Re(m)>0$
\end{example}
\begin{example}
From Eq.(23) in \cite{ram1}
\begin{multline}
\int_0^{\infty } \frac{1}{2} e^{-2 i m \pi  x} \cos \left(\pi  x^2\right) \left((-2 i \pi  x+\log (a))^k+e^{4 i m
   \pi  x} (2 i \pi  x+\log (a))^k\right) \text{sech}(\pi  x) \, dx\\
=2^{-\frac{1}{2}+k} e^{m \pi } \pi ^k \Phi \left(-e^{2
   m \pi },-k,\frac{\pi +\log (a)}{2 \pi }\right)\\
+\sum _{j=0}^{\infty } \sum _{p=0}^{\infty } \frac{(-1)^p 2^{2-j+k+4 p}
   e^{m \pi } m^{2-j+4 p} \pi ^{\frac{1}{2}-j+k+2 p} (1-j+k)_j (3-j+4 p)_{-\frac{3}{2}+j-2 p} }{\Gamma (1+j)}\\
\Phi \left(-e^{2 m \pi
   },j-k,\frac{\pi +\log (a)}{2 \pi }\right)
\end{multline}
\end{example}
\begin{example}
\begin{multline}
\int_0^{\infty } \frac{\left((a-2 i x)^k+(a+2 i x)^k\right) \left(\cos (2 m \pi  x) \cos \left(\pi 
   x^2\right)\right)}{\cosh (\pi  x)} \, dx\\
=2^{-\frac{1}{2}+k} e^{-m \pi } \left(\Phi \left(-e^{-2 m \pi
   },-k,\frac{1+a}{2}\right)+e^{2 m \pi } \Phi \left(-e^{2 m \pi },-k,\frac{1+a}{2}\right)\right)\\
+\sum _{j=0}^{\infty }
   \sum _{p=0}^{\infty } \frac{(-1)^p 2^{2-j+k+4 p} e^{-m \pi } (-m)^{-j} m^{2-j} \pi ^{\frac{1}{2}-j+2 p} (1-j+k)_j (3-j+4 p)_{-\frac{3}{2}+j-2 p}}{\Gamma (1+j)}\\
\left((-m)^{4
   p} m^j \Phi \left(-e^{-2 m \pi },j-k,\frac{1+a}{2}\right)+e^{2 m \pi } (-m)^j m^{4 p} \Phi \left(-e^{2 m \pi
   },j-k,\frac{1+a}{2}\right)\right) 
\end{multline}
\end{example}
\begin{example}
\begin{multline}
\int_0^{\infty } \frac{\cos \left(\pi  x^2\right) }{\left(25+4 x^2\right) }\frac{\cos (\pi  x)}{\cosh (\pi  x)} \,
   dx\\
=\frac{e^{-\frac{5 \pi }{2}} \left(e^{\pi } \left(-2+e^{\pi }+e^{2 \pi }-2 e^{3 \pi }+2 e^{4 \pi } \log
   \left(1+e^{-\pi }\right)\right)+2 \log \left(1+e^{\pi }\right)\right)}{40 \sqrt{2}}\\
+\sum _{j=0}^{\infty } \sum
   _{p=0}^{\infty } \frac{(-1)^{-j+p} e^{-\frac{\pi }{2}} \pi ^{\frac{1}{2}-j+2 p} \Gamma \left(\frac{3}{2}+2 p\right)
   (-j)_j \left((-1)^j e^{\pi } \Phi \left(-e^{\pi },1+j,3\right)+(-1)^{4 p} \right)}{20
   \Gamma (1+j) \Gamma (3-j+4 p)}\\
\Phi \left(-e^{-\pi },1+j,3\right)
\end{multline}
\end{example}
\begin{example}
\begin{multline}
\int_0^{\infty } \frac{\cos (2 m \pi  x) \cos \left(\pi  x^2\right) \text{sech}(\pi  x)}{a^2+4 x^2} \,
   dx\\
=\frac{e^{-m \pi } \left(\Phi \left(-e^{-2 m \pi },1,\frac{1+a}{2}\right)+e^{2 m \pi } \Phi \left(-e^{2 m \pi
   },1,\frac{1+a}{2}\right)\right)}{2 a 2 \sqrt{2}}\\
+\sum _{j=0}^{\infty } \sum _{p=0}^{\infty } \frac{(-1)^{-j+p} 2^{1-j+4
   p} e^{-m \pi } m^{2-2 j} \pi ^{\frac{1}{2}-j+2 p} (-j)_j (3-j+4 p)_{-\frac{3}{2}+j-2 p}}{2 a \Gamma (1+j)}\\
\left( m^{j+4 p} \Phi \left(-e^{-2 m \pi
   },1+j,\frac{1+a}{2}\right)+(-1)^j e^{2 m \pi } m^{j+4 p} \Phi \left(-e^{2 m \pi },1+j,\frac{1+a}{2}\right)\right)
\end{multline}
\end{example}
\begin{example}
From Eq.(23) and (24) in \cite{ram1}
\begin{multline}
\int_0^{\infty } e^{-i \pi  (2 m-x) x} \left((a-2 i x)^k+e^{4 i m \pi  x} (a+2 i x)^k\right) \text{sech}(\pi 
   x) \, dx\\
=-(-1)^{3/4} 2^{1+k} e^{m \pi } \Phi \left(-e^{2 m \pi },-k,\frac{1+a}{2}\right)\\
+\sum _{n=0}^{\infty }
   \sum _{p=0}^{\infty } \frac{2^{1+k} e^{m \pi } i (-1)^n i^n 2^{-p} m^{2 n-p} \pi ^{n-p} \binom{2 n}{p} \Phi
   \left(-e^{2 m \pi },-k+p,\frac{1+a}{2}\right) (1+k-p)_p}{n!}
\end{multline}
\end{example}
\begin{example}
\begin{multline}
\int_0^{\infty } \frac{e^{i \pi  x^2} x \text{sech}(\pi  x) \sin (2 m \pi
    x)}{a^2+x^2} \, dx\\
=\frac{\sqrt[4]{-1} e^{-m \pi } \left(-\Phi \left(-e^{-2 m
   \pi },1,\frac{1}{2} (1+2 a)\right)+e^{2 m \pi } \Phi \left(-e^{2 m \pi
   },1,\frac{1}{2} (1+2 a)\right)\right)}{2 i}\\
+\sum _{n=0}^{\infty } \sum
   _{p=0}^{\infty } \frac{(-1)^n i^{n-1} 2^{-1-p} e^{-m \pi }
   (-p)_p}{n!}\\
\left(-m^2\right)^{-p} \pi ^{n-p} \binom{2 n}{p} \left((-m)^{2 n} m^p \Phi
   \left(-e^{-2 m \pi },1+p,\frac{1}{2} (1+2 a)\right)\right. \\ \left.
-e^{2 m \pi } (-m)^p m^{2
   n} \Phi \left(-e^{2 m \pi },1+p,\frac{1}{2} (1+2 a)\right)\right)
\end{multline}
\end{example}
\begin{example}
\begin{multline}
\int_0^{\infty } \frac{e^{i \pi  x^2} x \text{sech}(\pi  x) \sinh
   \left(\frac{\pi  x}{2}\right)}{4+x^2} \, dx
=\frac{1}{6} i \left((-4-4 i)-3
   \sqrt[4]{-1} \log \left(3-2 \sqrt{2}\right)\right)\\
+\sum _{n=0}^{\infty } \sum
   _{p=0}^{\infty } \frac{i (-1)^{\frac{1}{4}+n} i^{3 n+p} 2^{-1-4 n+p} \pi
   ^{n-p} \binom{2 n}{p} \left(-i (-1)^p \Phi
   \left(-i,1+p,\frac{5}{2}\right)+\Phi \left(i,1+p,\frac{5}{2}\right)\right)
   (-p)_p}{n!}
\end{multline}
\end{example}
\begin{example}
\begin{multline}
\int_0^{\infty } \frac{e^{i \pi  x^2} \cos (2 m \pi  x) \text{sech}(\pi  x)}{a^2+4 x^2} \, dx\\
=-\frac{(-1)^{3/4}
   e^{-m \pi } \left(\Phi \left(-e^{-2 m \pi },1,\frac{1+a}{2}\right)+e^{2 m \pi } \Phi \left(-e^{2 m \pi
   },1,\frac{1+a}{2}\right)\right)}{4 a}\\
+\sum _{n=0}^{\infty } \sum _{p=0}^{\infty } \frac{i (-1)^n i^n 2^{-2-p} e^{-m \pi } (-p)_p}{a n!}\\
\left(-m^2\right)^{-p} \pi ^{n-p} \binom{2 n}{p} \left((-m)^{2 n} m^p \Phi \left(-e^{-2 m \pi
   },1+p,\frac{1+a}{2}\right)\right. \\ \left.
+e^{2 m \pi } (-m)^p m^{2 n} \Phi \left(-e^{2 m \pi },1+p,\frac{1+a}{2}\right)\right)
\end{multline}
\end{example}
\begin{example}
\begin{multline}
\int_0^{\infty } \frac{e^{i \pi  x^2} \left(a^2-4 x^2\right) \cos (2 m \pi  x) \text{sech}(\pi  x)}{\left(a^2+4
   x^2\right)^2} \, dx\\
=-\frac{1}{8} (-1)^{3/4} e^{-m \pi } \left(\Phi \left(-e^{-2 m \pi },2,\frac{1+a}{2}\right)+e^{2
   m \pi } \Phi \left(-e^{2 m \pi },2,\frac{1+a}{2}\right)\right)\\
+\sum _{n=0}^{\infty } \sum _{p=0}^{\infty } \frac{i
   (-1)^n i^n 2^{-3-p} e^{-m \pi }  (-1-p)_p}{n!}\\
\left(-m^2\right)^{-p} \pi ^{n-p} \binom{2 n}{p} \left((-m)^{2 n} m^p \Phi
   \left(-e^{-2 m \pi },2+p,\frac{1+a}{2}\right)\right. \\ \left.
+e^{2 m \pi } (-m)^p m^{2 n} \Phi \left(-e^{2 m \pi
   },2+p,\frac{1+a}{2}\right)\right)
\end{multline}
\end{example}
\begin{example}
From Eq. (6.4.15) in \cite{erd_t1}
\begin{multline}
\int_0^{\infty } x^{-1+m} \log ^k(a x) \log (1+b x) \, dx\\
=-\sum _{j=0}^{\infty } b^{-m} e^{i m \pi } m^{-1-j} (2
   i \pi )^{1-j+k} \binom{-1}{j} (1-j+k)_j\\
 \Phi \left(e^{2 i m \pi },j-k,\frac{\pi -i \log (a)-i \log
   \left(\frac{1}{b}\right)}{2 \pi }\right)
\end{multline}
where $Re(m)<0,Im(b)\leq 0$
\end{example}
\begin{example}
\begin{multline}
\int_0^{\infty } \frac{x^{-1+m} \log (1+b x)}{(a \pi +\log (x))^2} \, dx
=\sum _{j=0}^{\infty } -i^{1-j} b^{-m}
   e^{i m \pi } (-1-j) m^{-1-j} (2 \pi )^{-1-j} \binom{-1}{j}\\
 \Phi \left(e^{2 i m \pi },2+j,\frac{\pi -i a \pi +i \log
   (b)}{2 \pi }\right) (-j)_j
\end{multline}
\end{example}
\begin{example}
From Eq. (4.359.2) in \cite{grad}
\begin{multline}
\int_0^1 \frac{e^{m x} \left(x^{p-1}-x^{q-1}\right) \log ^k\left(a e^x\right)}{\log (x)} \, dx=\sum
   _{j=0}^{\infty } \sum _{l=0}^{\infty } \frac{m^{j-l} (1+k-l)_l \binom{j}{l} \log ^{k-l}(a) \log
   \left(\frac{j+p}{j+q}\right)}{j!}
\end{multline}
\end{example}
\begin{example}
\begin{multline}
\int_0^1 \frac{e^{m x} \left(x^p-x^q\right)}{\left(-a^2 \pi ^2+x^2\right) \log (x)} \, dx\\
=\sum _{j=0}^{\infty }
   \sum _{l=0}^{\infty } \frac{(-1)^{-l} \left(-1+(-1)^l\right) (a \pi )^{-1-l} m^{j-l} \binom{j}{l} \log
   \left(\frac{j+p}{j+q}\right) (-l)_l}{2 j!}
\end{multline}
\end{example}
\begin{example}
\begin{multline}
\int_0^1 \frac{e^{m x} \left(x^p-x^q\right) \log (x+\log (a))}{x \sqrt{x+\log (a)} \log (x)} \, dx\\
=\sum
   _{j=0}^{\infty } \sum _{l=0}^{\infty } \frac{m^{j-l} \binom{j}{l} \log ^{-\frac{1}{2}-l}(a) \log
   \left(\frac{j+p}{j+q}\right) \left(-H_{-\frac{1}{2}-l}+\log \left(\frac{\log (a)}{4}\right)\right)
   \left(\frac{1}{2}-l\right)_l}{j!}
\end{multline}
\end{example}
\begin{example}
From Eq. (6.215.2) in \cite{grad}
\begin{multline}
\int_0^1 \frac{x^{-1-m} \log ^k\left(\frac{a}{x}\right) \text{li}(x)}{\sqrt{\log \left(\frac{1}{x}\right)}} \,
   dx=-\sum _{j=0}^{\infty } \sum _{n=0}^{\infty } \frac{2^{1-2 n} m^{-j+n} \sqrt{\pi } (1-j+k)_j \binom{n}{j} (2 n)!
   \log ^{-j+k}(a)}{(1+2 n) (n!)^2}
\end{multline}
where $Re(a)>e^{3\pi}$ in order for the sum to converge.
\end{example}
\begin{example}
\begin{multline}
\int_0^1 \frac{x^{-1-m} \text{li}(x)}{\sqrt{\log \left(\frac{1}{x}\right)} \left(a^2 \pi ^2-\log ^2(x)\right)}
   \, dx\\
=-\sum _{j=0}^{\infty } \sum _{n=0}^{\infty } \frac{2^{-2 n} (-a)^{-j} a^{-2-j} \left((-a)^j+a^j\right)
   m^{-j+n} \pi ^{-\frac{3}{2}-j} \binom{n}{j} (2 n)! (-j)_j}{(1+2 n) (n!)^2}
\end{multline}
\end{example}
\begin{example}
From Eq. (1.3.19) in \cite{erd_t1}
\begin{multline}
\int_0^{\infty } \frac{1}{2} e^{-i m x} \left((b-i x)^{-v}+(b+i
   x)^{-v}\right) \left((-i x+\log (a))^k+e^{2 i m x} (i x+\log (a))^k\right) \,
   dx\\
=\sum _{j=0}^{\infty } \frac{e^{-b m} m^{-1-j+v} \pi  (1-j+k)_j
   \binom{-1+v}{j} \log ^{-j+k}\left(a e^{-b}\right)}{\Gamma (v)}
\end{multline}
\end{example}
\begin{example}
From Eq. (1.4.11) in \cite{erd_t1}
\begin{multline}
\int_0^{\infty } \frac{1}{2} e^{-i m x-b x^2} \left((-i x+\log (a))^k+e^{2
   i m x} (i x+\log (a))^k\right) \, dx\\
=\sum _{j=0}^{\infty } \sum _{n=0}^{\infty
   } \frac{(-1)^n 2^{-1-2 n} b^{-\frac{1}{2}-n} m^{-j+2 n} \sqrt{\pi } \binom{2
   n}{j} \log ^{-j+k}(a) (1-j+k)_j}{n!}
\end{multline}
\end{example}
\begin{example}
\begin{multline}
\int_0^{\infty } \frac{e^{-b x^2} x \sin (m x)}{x^2+\log ^2(a)} \,
   dx\\
=-\sum _{j=0}^{\infty } \sum _{n=0}^{\infty } \frac{(-1)^n 2^{-2-2 n}
   b^{-\frac{1}{2}-n} \left(-m^2\right)^{-j} \left((-m)^{2 n} m^j-(-m)^j m^{2
   n}\right) \sqrt{\pi } \binom{2 n}{j} \log ^{-1-j}(a) (-j)_j}{n!}
\end{multline}
\end{example}
\begin{example}
\begin{multline}
\int_0^{\infty } \frac{e^{-b x^2} \cos (m x)}{x^2+\log ^2(a)} \, dx\\
=\sum
   _{j=0}^{\infty } \sum _{n=0}^{\infty } \frac{(-1)^n 2^{-2-2 n}
   b^{-\frac{1}{2}-n} \left(-m^2\right)^{-j} \left((-m)^{2 n} m^j+(-m)^j m^{2
   n}\right) \sqrt{\pi } \binom{2 n}{j} \log ^{-2-j}(a) (-j)_j}{n!}
\end{multline}
\end{example}
\begin{example}
From Eq. (4.1.4.3) in \cite{brychkov}
\begin{multline}
\int_0^b e^{-i m x} \left((-i x+\log (a \sin (b-x)))^k-e^{2 i m x} (i x+\log (a \sin (b-x)))^k\right) \sin
   ^{-2+m}(b-x) \, dx\\
=-\sum _{j=0}^{\infty } 2 i (-1+m)^{-1-j} \binom{-1}{j} \log ^{-j+k}(a \sin (b)) (1-j+k)_j \sin
   ^m(b)
\end{multline}
where $Re(a)<-e^{3\pi},Re(m)>1$
\end{example}
\begin{example}
From Eq. (6.16.4.1) in \cite{brychkov}
\begin{multline}
\sum _{k=0}^{\infty } \frac{2^k \left(\frac{e^w w}{z}\right)^k P_k^{(-k+\rho ,-k+\sigma )}(1+(1+k) z)}{(b+k)
   (k+1) (-\rho -\sigma )_k}\\
=e^{-b w} w^{-b} \int_0^w -\frac{e^{(-1+b) t} t^{-2+b} (1+t) z (1+\rho +\sigma )
   \left(-1+\, _1F_1\left(-1-\rho ;-1-\rho -\sigma ;-\frac{2 t}{z}\right)\right)}{2 (1+\rho )} \, dt
\end{multline}
where $Im(z)\neq 0$
\end{example}
\begin{example}
From Eq. (6.631.4(11)) in \cite{grad}
\begin{multline}
\int_0^{\infty } e^{-m x^2} x^{1+v} J_v(b x) \log ^k\left(a e^{-x^2}\right) \, dx\\
=\sum _{j=0}^{\infty } \sum
   _{n=0}^{\infty } \frac{(-1)^n 2^{-1-2 n-v} b^{2 n+v} m^{-1-j-n-v} (1-j+k)_j \binom{-1-n-v}{j} \log
   ^{-j+k}(a)}{n!}
\end{multline}
\end{example}
\begin{example}
\begin{multline}
\int_0^{\infty } \frac{e^{-(m+p) x^2} \left(e^{m x^2}-e^{p x^2}\right) x^{3+v} J_v(b x)}{a^2 \pi
   ^2-x^4} \, dx\\
=\sum _{j=0}^{\infty } \sum _{n=0}^{\infty } \frac{(-1)^n 2^{-2-2 n-v} (-a)^{-j} a^{-1-j}
   \left((-a)^j-a^j\right) b^{2 n+v} m^{-1-j-n-v} p^{-1-j-n-v} }{n!}\\ \times
\left(m^{1+j+n+v}-p^{1+j+n+v}\right) \pi ^{-1-j}
   \binom{-1-n-v}{j} (-j)_j
\end{multline}
\end{example}
\begin{example}
From Eq. (6.673.3) in \cite{grad}, Eq. (5.2.9.4)  and Eq. (5.2.9.5) in \cite{prud1} 
\begin{multline}
\int_0^{\frac{\pi }{2}} (\cos (x) I_0(a \cos (x))+I_1(a \cos (x))) \, dx=\sum _{k=1}^{\infty } k^k
   \left(\frac{2^{-k} e^{-\frac{k}{2}}}{k!}+\frac{\left(a e^{-a}\right)^k}{\Gamma
   (k+2)}\right)=\frac{-1+e^a}{a}
\end{multline}
\end{example}
\begin{example}
From Eq. (6.3.4.1) in \cite{prud3}
\begin{multline}
\sum _{j=0}^{\infty } \sum _{l=0}^{\infty } \frac{2^j e^{-2 j m \beta } m^{j-l} \binom{j}{l} E_j(\alpha +j \beta
   ) (\log (a)-2 j \beta )^{k-l}}{j! (k-l)!}\\
=\sum _{h=0}^{\infty } \frac{2^{1+k} e^{m+m (-1+2 \alpha )} (-\beta )^h
   (1-2 m \beta )^{-1-h} \binom{-1}{h} \Phi \left(-e^{2 m},h-k,\frac{1}{2} (2 \alpha +\log
   (a))\right)}{(-h+k)!}
\end{multline}
where $Re(a)>e^{3\pi}$
\end{example}
\begin{example}
From Eq. (6.3.2.1) in \cite{prud3}
\begin{multline}
\sum _{j=0}^{\infty } \sum _{p=0}^{\infty } \frac{e^{-2 j m \beta } (2 m)^{j-p} (a \pi -j \beta )^{k-p}
   B_j(\alpha +j \beta ) \binom{j}{p}}{j! (k-p)!}\\
=-\sum _{j=0}^{\infty } \frac{e^{2 m \alpha } (-\beta )^j (1-2 m \beta
   )^{-1-j} \binom{-1}{j} \left(2 m (-1-j+k)! \right)}{(-1-j+k)! (-j+k)!}\\
\Phi \left(e^{2 m},j-k,a \pi +\alpha \right)+(-j+k)! \Phi \left(e^{2
   m},1+j-k,a \pi +\alpha \right)
\end{multline}
where $Re(a)>7\pi$
\end{example}
\begin{example}
From pp.328 in \cite{sylvester}
\begin{multline}
\int_0^{\infty } e^{-2 i m x-x^2} \left((-2 i x+\log (a))^k+e^{4 i m x} (2 i x+\log (a))^k\right) \, dx\\
=\sum
   _{j=0}^{\infty } \sum _{n=0}^{\infty } \frac{(-1)^n m^{-j+2 n} \sqrt{\pi } (1-j+k)_j \binom{2 n}{j} \log
   ^{-j+k}(a)}{n!}
\end{multline}
where $Re(a)>e^{3\pi}$
\end{example}
\begin{example}
\begin{multline}
\int_0^{\infty } \frac{e^{-x^2} \left(2 \gamma +\psi ^{(0)}\left(-\frac{-9 \pi -2 i x}{9 \pi }\right)+\psi
   ^{(0)}\left(-\frac{-9 \pi +2 i x}{9 \pi }\right)\right) \sin (2 m x)}{x} \, dx\\
=-\sum _{j=0}^{\infty } \sum
   _{n=0}^{\infty } \frac{(-1)^n 9^{-1-j} \left(-m^2\right)^{-j} \left((-m)^{2 n} m^j-(-m)^j m^{2 n}\right) \pi
   ^{-\frac{1}{2}-j} \binom{2 n}{j} (-j)_j \zeta (2+j)}{n!}
\end{multline}
where $Im(m)\neq 0$
\end{example}
\begin{example}
\begin{multline}
\int_0^{\infty } \frac{e^{-x^2} \cos (2 m x) \left(-9+2 x \coth \left(\frac{2 x}{9}\right)\right)}{72 x^2} \,
   dx\\
=\sum _{j=0}^{\infty } \sum _{n=0}^{\infty } \frac{(-1)^n 9^{-2-j} \left(-m^2\right)^{-j} \left((-m)^{2 n}
   m^j+(-m)^j m^{2 n}\right) \pi ^{-\frac{3}{2}-j} \binom{2 n}{j} (-j)_j \zeta (2+j)}{4 n!}
\end{multline}
where $Im(m)\neq 0$
\end{example}
\begin{example}
From pp. 341 in \cite{sylvester}
\begin{multline}
\int_0^{\infty } e^{e^{-2 i b x} m-x^2} \left(\left(e^{-2 i b x}+\log (a)\right)^k+e^{2 i m \sin (2 b x)}
   \left(e^{2 i b x}+\log (a)\right)^k\right) \, dx\\
=\sum _{j=0}^{\infty } \sum _{n=0}^{\infty } \frac{e^{-b^2 n}
   m^{-j+n} \sqrt{\pi } (1-j+k)_j \binom{n}{j} \log ^{-j+k}(a)}{n!}
\end{multline}
\end{example}
\begin{example}
\begin{multline}
\int_0^{\infty } e^{e^{-2 i b x} m-2 i b x-x^2} \left(e^{4 i b x} H_{\frac{e^{-2 i b x}}{9 \pi
   }}+e^{2 i m \sin (2 b x)} H_{\frac{e^{2 i b x}}{9 \pi }}\right) \, dx\\
=\sum _{j=0}^{\infty } \sum
   _{n=0}^{\infty } \frac{9^{-1-j} e^{-b^2 n} m^{-j+n} \pi ^{-\frac{1}{2}-j} \binom{n}{j} (-j)_j \zeta
   (2+j)}{n!}
\end{multline}
\end{example}
\begin{example}
Eq. (456.2) in \cite{bdh}
\begin{multline}
\int_0^{\infty } \frac{e^{-i m \tan ^{-1}\left(\frac{x}{q}\right)} }{\left(s^2+x^2\right) \left(q^2+x^2\right)^{m/2}} \left(\left(-i \tan
   ^{-1}\left(\frac{x}{q}\right)+\log \left(\frac{a}{\sqrt{q^2+x^2}}\right)\right)^k\right. \\ \left.
+e^{2 i m \tan
   ^{-1}\left(\frac{x}{q}\right)} \left(i \tan ^{-1}\left(\frac{x}{q}\right)+\log
   \left(\frac{a}{\sqrt{q^2+x^2}}\right)\right)^k\right)
dx\\
=\frac{\pi  (q+s)^{-m} \log ^k\left(\frac{a}{q+s}\right)}{s}
\end{multline}
where $Re(a)< 0$
\end{example}
\begin{example}
\begin{multline}
\int_0^{\infty } \frac{e^{-i m \tan ^{-1}\left(\frac{x}{q}\right)} (-1+\pi  x \coth (\pi  x)) }{x^2 \left(q^2+x^2\right)^{m/2}}\\\left(\left(-i
   \tan ^{-1}\left(\frac{x}{q}\right)+\log \left(\frac{a}{\sqrt{q^2+x^2}}\right)\right)^k\right. \\ \left.
+e^{2 i m \tan
   ^{-1}\left(\frac{x}{q}\right)} \left(i \tan ^{-1}\left(\frac{x}{q}\right)+\log
   \left(\frac{a}{\sqrt{q^2+x^2}}\right)\right)^k\right)
 \, dx\\
=\sum _{s=1}^{\infty }
   \frac{2 \pi  \log ^k\left(\frac{a}{q+s}\right)}{s (q+s)^m}
\end{multline}
where $Re(a)< 0,Re(m)>0$
\end{example}
\begin{example}
From Eq. (456.2) in \cite{bdh}
\begin{multline}
\int_0^1 x^{-m x} \log ^k\left(a x^{-x}\right) \, dx=\sum _{n=0}^{\infty } \sum _{j=0}^{n-1} m^{-1-j+n} n^{-n}
   \binom{-1+n}{j} (1-j+k)_j \log ^{-j+k}(a)
\end{multline}
\end{example}
\begin{example}
\begin{multline}
\int_0^1 \frac{x^{-2-m x} \left(-1+x^{2 m x}\right) (-1+\pi  x \cot (\pi  x \log (x)) \log (x))}{\log ^2(x)} \,
   dx\\
=\sum _{n=0}^{\infty } \sum _{j=0}^{n-1} e^{-i j \pi } \left(1+e^{i j \pi }\right) (-m)^{-j} m^{-1-j}
   \left((-m)^n m^j+(-m)^j m^n\right) n^{-n} \binom{-1+n}{j} (-j)_j \zeta (2+j)
\end{multline}
\end{example}
\begin{example}
\begin{multline}
\int_0^1 x^{-m x} \left(1+x^{2 m x}\right) \psi ^{(1)}(1-x \log (x)) \, dx\\
=\sum _{n=0}^{\infty } \sum
   _{j=0}^{n-1} (1+j) (-m)^{-j} m^{-1-j} \left(-(-m)^n m^j+(-m)^j m^n\right) n^{-n} \binom{-1+n}{j} (-j)_j \zeta
   (2+j)
\end{multline}
\end{example}
\begin{example}
\begin{multline}
\int_0^1 \frac{x^{-m x} \log (a-x \log (x))}{\sqrt{a-x \log (x)}} \, dx\\
=\sum _{n=0}^{\infty } \sum _{j=0}^{n-1}
   \left(a^{-\frac{1}{2}-j} m^{-1-j+n} n^{-n} \binom{-1+n}{j} \log (a) \left(\frac{1}{2}-j\right)_j\right. \\ \left.
+a^{-\frac{1}{2}-j}
   m^{-1-j+n} n^{-n} \binom{-1+n}{j} \left(\frac{1}{2}-j\right)_j \left(\psi ^{(0)}\left(\frac{1}{2}\right)-\psi
   ^{(0)}\left(\frac{1}{2}-j\right)\right)\right)
\end{multline}
\end{example}
\begin{example}
\begin{multline}
\int_0^1 \frac{x^{-x} \log (1-x \log (x))}{\sqrt{1-x \log (x)}} \, dx\\
=\sum _{n=0}^{\infty } \sum _{j=0}^{-1+n}
   n^{-n} \binom{-1+n}{j} \left(\frac{1}{2}-j\right)_j \left(\psi ^{(0)}\left(\frac{1}{2}\right)-\psi
   ^{(0)}\left(\frac{1}{2}-j\right)\right)
\end{multline}
\end{example}
\begin{example}
\begin{multline}
\sum _{n=1}^{\infty } \sum _{j=0}^{\infty } \sum _{p=0}^{2 j} \frac{(-1)^n
   m^{2 j-p} \left(n^2\right)^{-1-j} \binom{-1}{j} \binom{2 j}{p} \log
   ^{k-p}(a)}{(k-p)!}\\
=\sum _{j=0}^{\infty } -\frac{m^{-2-j} \binom{-2}{j} \log
   ^{-j+k}(a)}{2 (-j+k)!}\\
   -\sum _{j=0}^{\infty } \frac{2^{-j+k} e^{m \pi }
   m^{-1-j} \pi ^{1-j+k} \binom{-1}{j} \Phi \left(e^{2 m \pi },j-k,\frac{\pi
   +\log (a)}{2 \pi }\right)}{(-j+k)!}
\end{multline}
where $Re(a)>e^{3\pi},Re(m)<0$
\end{example}
\begin{example}
\begin{multline}
\sum _{n=1}^{\infty } \sum _{j=0}^{\infty } \sum _{p=0}^{2 j} \frac{(-1)^n
   \left(-\frac{i}{2}\right)^{2 j-p} \left(n^2\right)^{-1-j} \binom{-1}{j}
   \binom{2 j}{p} \log ^{k-p}(a)}{(k-p)!}\\
=-\sum _{j=0}^{\infty }
   \frac{(-i)^{-2-j} 2^{1+j} \binom{-2}{j} \log ^{-j+k}(a)}{(-j+k)!}\\
+\sum
   _{j=0}^{\infty } \frac{i (-i)^{-1-j} 2^{1+k} \pi ^{1-j+k} \binom{-1}{j}
   \left(2^{-j+k} \zeta \left(j-k,\frac{\pi +\log (a)}{4 \pi }\right)-2^{-j+k}
   \zeta \left(j-k,\frac{1}{2} \left(1+\frac{\pi +\log (a)}{2 \pi
   }\right)\right)\right)}{(-j+k)!}
\end{multline}
\end{example}
\begin{example}
\begin{multline}
\sum _{n=1}^{\infty } \sum _{j=0}^{\infty } \sum _{p=0}^{2 j} \frac{(-1)^n
   \left(\frac{i}{2}\right)^{2 j-p} \left(n^2\right)^{-1-j} \binom{-1}{j}
   \binom{2 j}{p} \log ^{k-p}(a)}{(k-p)!}\\
=-\sum _{j=0}^{\infty } \frac{i^{-2-j}
   2^{1+j} \binom{-2}{j} \log ^{-j+k}(a)}{(-j+k)!}\\
-\sum _{j=0}^{\infty } \frac{i
   i^{-1-j} 2^{1+k} \pi ^{1-j+k} \binom{-1}{j} \left(2^{-j+k} \zeta
   \left(j-k,\frac{\pi +\log (a)}{4 \pi }\right)-2^{-j+k} \zeta
   \left(j-k,\frac{1}{2} \left(1+\frac{\pi +\log (a)}{2 \pi
   }\right)\right)\right)}{(-j+k)!}
\end{multline}
\end{example}
\begin{example}
From Eq. (43.16) in \cite{bdh}
\begin{multline}
\int_0^{\frac{\pi }{2}} e^{-i m \cos (x)} \cos (2 q x) \left((-i \cos (x)+\log (a))^k+e^{2 i m \cos (x)} (i
   \cos (x)+\log (a))^k\right) \, dx\\
=\frac{\log ^k(a) \sin (\pi  q)}{q}+\sum _{n=1}^{\infty } \sum _{j=0}^{2 n}
   \frac{\left(-\frac{1}{4}\right)^n m^{-j+2 n} \binom{2 n}{j} \log ^{-j+k}(a) (1-j+k)_j \sin (\pi  q)}{q (1-q)_n
   (1+q)_n}
\end{multline}
\end{example}
\begin{example}
\begin{multline}
\int_0^{\frac{\pi }{2}} \frac{\cos (x) \cos (2 q x) \sin (m \cos (x))}{\cos ^2(x)+\log ^2(a)} \, dx\\
=-\sum
   _{n=1}^{\infty } \sum _{j=0}^{2 n} \frac{(-1)^n 2^{-2-2 n} \left(-m^2\right)^{-j} \left((-m)^{2 n} m^j-(-m)^j m^{2
   n}\right) \binom{2 n}{j} \log ^{-1-j}(a) (-j)_j \sin (\pi  q)}{q (1-q)_n (1+q)_n}
\end{multline}
\end{example}
\begin{example}
\begin{multline}
\int_0^{\frac{\pi }{2}} \frac{\cos (2 q x) \cos (m \cos (x)) \log (a)}{\cos ^2(x)+\log ^2(a)} \, dx
=\frac{\sin
   (\pi  q)}{2 q \log (a)}\\+\sum _{n=1}^{\infty } \sum _{j=0}^{2 n} \frac{(-1)^n 4^{-1-n} \left(-m^2\right)^{-j}
   \left((-m)^{2 n} m^j+(-m)^j m^{2 n}\right) \binom{2 n}{j} \log ^{-1-j}(a) (-j)_j \sin (\pi  q)}{q (1-q)_n
   (1+q)_n}
\end{multline}
\end{example}
\begin{example}
From Eq. (81.1) in \cite{bdh}
\begin{multline}
\sum _{j=0}^{\infty } m^{-j-p} \binom{-p}{j} (1-j+k)_j \Gamma (p) \log ^{-j+k}(a)\\
=\frac{m^{-k-p} \pi  \csc
   ((k+p) \pi ) \, _1F_1(-k;1-k-p;-m \log (a)) (-1)^k}{\Gamma (1-k-p)}\\
-\frac{\pi  \csc ((k+p) \pi ) \Gamma (p) \,
   _1F_1(p;1+k+p;-m \log (a)) (-\log (a))^p \log ^k(a)}{\Gamma (-k) \Gamma (1+k+p)}
\end{multline}
where $Re(a)>e^{3\pi},Re(m)>1$
\end{example}
\begin{example}
From Eq. (110.13) in \cite{bdh}
\begin{multline}
\int_0^1 \frac{x^m \log ^{-1+p}(x) \log ^k(a x)}{1-x^2} \, dx\\
=\sum _{j=0}^{\infty } \sum _{n=0}^{\infty }
   (-1)^{-1+p} (1+m+2 n)^{-j-p} (1-j+k)_j \binom{-p}{j} \Gamma (p) \log ^{-j+k}(a)
\end{multline}
where $Re(a)>e^{3\pi},Re(m)<0$
\end{example}
\begin{example}
\begin{multline}
\int_0^1 \frac{x^{-m} \left(-1+x^{2 m}\right) \log ^{-1+p}(x)}{\left(-1+x^2\right) \log (a x)} \, dx\\
=\sum
   _{j=0}^{\infty } \sum _{n=0}^{\infty } (-1)^p \left(\left(\frac{1}{1+m+2 n}\right)^{j+p}-\left(\frac{1}{1-m+2
   n}\right)^{j+p}\right) \binom{-p}{j} \Gamma (p) \log ^{-1-j}(a) (-j)_j
\end{multline}
\end{example}
\begin{example}
From Eq. (113.11) in \cite{bdh} errata.
\begin{multline}
\int_0^1 \frac{x^{-1+m} (p x-\cos (\lambda )) \log ^{-1+r}(x) \log ^k(a x)}{1+p^2 x^2-2 p x \cos (\lambda )} \,
   dx\\
=\sum _{j=0}^{\infty } \sum _{n=1}^{\infty } (-1)^r (-1+m+n)^{-j-r} p^{-1+n} (1-j+k)_j \binom{-r}{j} \cos
   (\lambda  n) \Gamma (r) \log ^{-j+k}(a)
\end{multline}
\end{example}
\begin{example}
From Eq. (256.9) in \cite{bdh}, Lobatto, N.V.
\begin{multline}
\int_0^{\infty } e^{-m x^2} \left(a^2-x^2\right)^k \log
   \left(q^2+x^2\right) \, dx\\
=\sum _{n=1}^{\infty } \sum _{j=0}^{\infty }
   -\left(-\frac{1}{4}\right)^n m^{-\frac{1}{2}-j-n} (1+n)^{-1+n} \sqrt{\pi }
   q^{-2 n} (1-j+k)_j \binom{-\frac{1}{2}-n}{j} a^{2 (-j+k)}\\
+\sum _{j=0}^{\infty
   } m^{-\frac{1}{2}-j} \sqrt{\pi } (1-j+k)_j \binom{-\frac{1}{2}}{j} a^{2
   (-j+k)} \log (q)
\end{multline}
\end{example}
\begin{example}
\begin{multline}
\int_0^{\infty } \frac{e^{-m x^2} \left(-4 z+x \cot \left(\frac{x}{4
   z}\right)\right) \log \left(q^2+x^2\right)}{x^2} \, dx\\
=\sum _{n=1}^{\infty }
   \sum _{j=0}^{\infty } (-1)^n 2^{-1-4 j-2 n} m^{-\frac{1}{2}-j-n} (1+n)^{-1+n}
   \pi ^{-\frac{3}{2}-2 j} q^{-2 n} z^{-1-2 j} \binom{-\frac{1}{2}-n}{j} (-j)_j
   \zeta (2 (1+j))\\
-\sum _{j=0}^{\infty } 2^{-1-4 j} m^{-\frac{1}{2}-j} \pi
   ^{-\frac{3}{2}-2 j} z^{-1-2 j} \binom{-\frac{1}{2}}{j} \log (q) (-j)_j \zeta
   (2 (1+j))
\end{multline}
\end{example}
\begin{example}
\begin{multline}
\int_0^{\infty } \frac{e^{-\pi  x^2} \left(-4 \pi +x \cot \left(\frac{x}{4
   \pi }\right)\right) \log \left(\pi ^2+x^2\right)}{x^2} \, dx\\
=\sum
   _{n=1}^{\infty } \sum _{j=0}^{\infty } (-1)^n 2^{-1-4 j-2 n} (1+n)^{-1+n} \pi
   ^{-3-5 j-3 n} \binom{-\frac{1}{2}-n}{j} (-j)_j \zeta (2 (1+j))\\
-\sum
   _{j=0}^{\infty } 2^{-1-4 j} \pi ^{-3-5 j} \binom{-\frac{1}{2}}{j} \log (\pi )
   (-j)_j \zeta (2 (1+j))
\end{multline}
\end{example}
\begin{example}
\begin{multline}
\int_0^{\infty } \frac{e^{-x^2} \left(-12+\pi  x \cot \left(\frac{\pi 
   x}{12}\right)\right) \log \left(4 \pi ^2+x^2\right)}{x^2} \, dx\\
=\sum
   _{n=1}^{\infty } \sum _{j=0}^{\infty } (-1)^n 2^{-1-4 j-4 n} 3^{-1-2 j}
   (1+n)^{-1+n} \pi ^{\frac{1}{2}-2 n} \binom{-\frac{1}{2}-n}{j} (-j)_j \zeta (2
   (1+j))\\
-\sum _{j=0}^{\infty } 2^{-1-4 j} 3^{-1-2 j} \sqrt{\pi }
   \binom{-\frac{1}{2}}{j} \log (2 \pi ) (-j)_j \zeta (2 (1+j))
\end{multline}
\end{example}
\begin{example}
From Eq. (266.1) in \cite{bdh}
\begin{multline}
\int_0^{\infty } \frac{e^{-m x^2} \left(\log (a)-x^2\right)^k}{1+r^2-2 r \cos (x)} \, dx=-\sum _{j=0}^{\infty }
   \frac{m^{-\frac{1}{2}-j} \sqrt{\pi } (1-j+k)_j \binom{-\frac{1}{2}}{j} \log ^{-j+k}(a)}{2 \left(-1+r^2\right)}\\
+\sum
   _{j=0}^{\infty } \sum _{n=1}^{\infty } \sum _{p=0}^{\infty } \sum _{l=0}^p \frac{4^{-p} m^{-\frac{1}{2}-j-l-p}
   \left(-n^2\right)^p \sqrt{\pi } r^n (1-j+k-l)_{j+l} \binom{-\frac{1}{2}}{j} \binom{-p}{l} \log
   ^{-j+k-l}(a)}{\left(1-r^2\right) p!}
\end{multline}
\end{example}
\begin{example}
\begin{multline}
\int_0^{\infty } \frac{e^{-m x^2}}{\left(-a^2+x^4\right) \left(1+r^2-2 r \cos (x)\right)} \, dx\\
=\sum
   _{j=0}^{\infty } \frac{\left(-a^2\right)^{-j} \left((-a)^j+a^j\right) m^{-\frac{1}{2}-j} \sqrt{\pi }
   \binom{-\frac{1}{2}}{j} (-j)_j}{4 a^2 \left(-1+r^2\right)}\\
+\sum _{j=0}^{\infty } \sum _{n=1}^{\infty } \sum
   _{p=0}^{\infty } \sum _{l=0}^p \frac{(-1)^p 2^{-1-2 p} \left(-a^2\right)^{-j-l} \left((-a)^{j+l}+a^{j+l}\right)
   m^{-\frac{1}{2}-j-l-p} n^{2 p} \sqrt{\pi } r^n \binom{-\frac{1}{2}}{j} \binom{-p}{l} (-j-l)_{j+l}}{a^2
   \left(-1+r^2\right) p!}
\end{multline}
\end{example}
\begin{example}
\begin{multline}
\int_0^{\infty } \frac{e^{-5 \pi  x^2} \left(\frac{1}{6}-\cos (x)\right)}{\left(729 \pi ^3-9 \pi  x^4\right)
   \left(\frac{37}{36}-\frac{\cos (x)}{3}\right)^2} \, dx\\
=\sum _{j=0}^{\infty } \frac{1}{49} (-2) (-1)^{-j} 3^{-3-4 j}
   5^{-\frac{5}{2}-j} \pi ^{-3-3 j} \left((-9 \pi )^j+(9 \pi )^j\right) \binom{-\frac{1}{2}}{j} (-j)_j\\
+\sum
   _{j=0}^{\infty } \sum _{n=1}^{\infty } \sum _{p=0}^{\infty } \sum _{l=0}^p \frac{(-1)^{-j-l+p} 2^{3-n-2 p} 3^{-1-4
   j-4 l-n} 5^{-\frac{5}{2}-j-l-p} \left(-\frac{1}{18}-\frac{35 n}{36}\right) n^{2 p} \pi ^{-3-j+2 (-j-l)-l-p}
  }{49 p!}\\
 \left((-9 \pi )^{j+l}+(9 \pi )^{j+l}\right) \binom{-\frac{1}{2}}{j} \binom{-p}{l} (-j-l)_{j+l}
\end{multline}
\end{example}
\begin{example}
\begin{multline}
\int_0^{\infty } \frac{e^{-m x^2} x^2}{\left(-a^2+x^4\right) \left(1+r^2-2 r \cos (x)\right)} \, dx\\
=-\sum
   _{j=0}^{\infty }\frac{a \left(-a^2\right)^{-1-j} \left((-a)^j-a^j\right) m^{-\frac{1}{2}-j} \sqrt{\pi }
   \binom{-\frac{1}{2}}{j} (-j)_j}{4 \left(-1+r^2\right)}\\
-\sum _{j=0}^{\infty } \sum _{n=1}^{\infty } \sum_{p=0}^{\infty } \sum _{l=0}^p \frac{(-1)^p 2^{-1-2 p} a \left(-a^2\right)^{-1-j-l} \left((-a)^{j+l}-a^{j+l}\right)
   m^{-\frac{1}{2}-j-l-p} n^{2 p} \sqrt{\pi } r^n \binom{-\frac{1}{2}}{j} \binom{-p}{l}
   (-j-l)_{j+l}}{\left(-1+r^2\right) p!}
\end{multline}
\end{example}
\begin{example}
From Eq. (289.12) in \cite{bdh}
\begin{multline}
\int_0^{\frac{\pi }{4}} \log ^k(a (\csc (x) (\cos (x)-\sin (x)))) \csc ^2(x) \log (\tan (x)) (\csc (x) (\cos
   (x)-\sin (x)))^{-1+m} \, dx\\
=\sum _{j=0}^{\infty } e^{i m \pi } m^{-1-j} (2 i \pi )^{1-j+k} (1-j+k)_j \binom{-1}{j}
   \Phi \left(e^{2 i m \pi },j-k,\frac{\pi -i \log (a)}{2 \pi }\right)
\end{multline}
where $Re(a)<0,Re(m)<0$.
\end{example}
\begin{example}
\begin{multline}
\int_0^1 \frac{e^{-i m t} \left(1-t^2\right)^{-\frac{1}{2}+v} \left((-i t+\log (a))^k+e^{2 i m t} (i t+\log
   (a))^k\right)}{\sqrt{\pi } \Gamma \left(\frac{1}{2}+v\right)} \, dt\\
=\sum _{j=0}^{\infty } \sum _{p=0}^{\infty }
   \frac{(-1)^j 2^{-2 j} m^{2 j-p} (1+k-p)_p \binom{2 j}{p} \log ^{k-p}(a)}{j! \Gamma (1+j+v)}
\end{multline}
\end{example}
\begin{example}
\begin{multline}
\int_0^1 \frac{\left(1-t^2\right)^{-\frac{1}{2}+v} \left(2 \gamma +\psi ^{(0)}\left(-\frac{-\pi -i t}{\pi
   }\right)+\psi ^{(0)}\left(-\frac{-\pi +i t}{\pi }\right)\right) \sin (m t)}{t \Gamma \left(\frac{1}{2} (1+2
   v)\right)} \, dt\\
=\sum _{j=0}^{\infty } \sum _{p=0}^{\infty } \frac{(-1)^j 2^{-1-2 j} \left(1-(-1)^{-p}\right) m^{2
   j-p} \pi ^{-\frac{1}{2}-p} \binom{2 j}{p} (-p)_p \zeta (2+p)}{j! \Gamma (1+j+v)}
\end{multline}
\end{example}
\begin{example}
From Eq. (2.4.17) in \cite{erd_t1}
\begin{multline}
\int_0^{\infty } e^{-i m x-x \alpha } \left(1-e^{-x \beta }\right)^{-1+v} \left((-i x+\log (a))^k-e^{2 i m x}
   (i x+\log (a))^k\right) \, dx\\
=\sum _{j=0}^{\infty } \sum _{l=0}^{\infty } \sum _{p=0}^j \frac{(-1)^{j+p} \beta
   ^{-1-p} \left((-1)^l i^l (-i m-\alpha +\beta )^{-l+p}-i^l (i m-\alpha +\beta )^{-l+p}\right) }{(j+v) j!}\\\times
\binom{p}{l} \log
   ^{k-l}(a) (1+k-l)_l S_j^{(p)}
\end{multline}
where $Re(a)>e^{3\pi},Re(\alpha) >Re(\beta) < Re(v) < Re(m)$
\end{example}
\begin{example}
\begin{multline}
\int_0^{\infty } \exp (-i x (m-i \alpha )) \left(1-e^{-x \beta }\right)^{-1+v} \left((a \pi -x)^k-e^{2 i m x}
   (a \pi +x)^k\right) \, dx\\
=\sum _{j=0}^{\infty } \sum _{l=0}^{\infty } \sum _{p=0}^j \frac{(-1)^{j+p} (a \pi )^{k-l}
   \beta ^{-1-p} \left(\frac{(-1)^l}{(-i m-\alpha +\beta )^{l-p}}-\frac{1}{(i m-\alpha +\beta )^{l-p}}\right)
   \binom{p}{l} (1+k-l)_l S_j^{(p)}}{(j+v) j!}
\end{multline}
where $Re(a)>e^{3\pi},Re(\alpha) >Re(\beta) < Re(v) < Re(m)$
\end{example}
\begin{example}
\begin{multline}
\int_0^{\infty } \frac{e^{x (-\alpha +\beta )} \left(1-e^{-x \beta }\right)^v \left(-4+x \cot
   \left(\frac{x}{4}\right)\right) \sin (m x)}{\left(-1+e^{x \beta }\right) x^2} \, dx\\
=\sum _{j=0}^{\infty } \sum
   _{l=0}^{\infty } \sum _{p=0}^j \frac{i (-1)^{j+p} 2^{-3-2 l} \left(1+(-1)^l\right) \pi ^{-2-l} \beta ^{-1-p} (-i
   m-\alpha +\beta )^{-l} (i m-\alpha +\beta )^{-l} }{(j+v) j!}\\ \times
\left(-(-i m-\alpha +\beta )^p (i m-\alpha +\beta )^l+(-i m-\alpha+\beta )^l (i m-\alpha +\beta )^p\right) \binom{p}{l} (-l)_l S_j^{(p)} \zeta (2+l)
\end{multline}
\end{example}
\begin{example}
\begin{multline}
\int_0^{\infty } e^{-m x} x^v J_v(x) (\log (a)-x)^k \, dx\\=\sum _{j=0}^{\infty } \sum _{p=0}^{\infty } \frac{2^v
   m^{-1-2 j-p-2 v} (1+k-p)_p \binom{-1-2 j-2 v}{p} \binom{-\frac{1}{2}-v}{j} \Gamma \left(\frac{1}{2}+v\right) \log
   ^{k-p}(a)}{\sqrt{\pi }}
\end{multline}
\end{example}
\begin{example}
\begin{multline}
\int_0^{\infty } \frac{e^{-m x} x^{-2+v} J_v(x) \left(-4 z+x \cot \left(\frac{x}{4 z}\right)\right)}{z} \,
   dx\\
=\sum _{j=0}^{\infty } \sum _{p=0}^{\infty } -(-1)^{-p} 2^{-2-2 p+v} \left(1+(-1)^p\right) m^{-1-2 j-p-2 v} \pi
   ^{-\frac{5}{2}-p} z^{-2-p} \binom{-1-2 j-2 v}{p}\\
 \binom{-\frac{1}{2}-v}{j} \Gamma \left(\frac{1}{2} (1+2 v)\right)
   (-p)_p \zeta (2+p)
\end{multline}
\end{example}
\begin{example}
\begin{multline}
\int_0^{\infty } \frac{e^{-2 x} J_0(x) \left(-4 \pi +x \cot \left(\frac{x}{4 \pi }\right)\right)}{x^2} \,
   dx\\
=-\sum _{j=0}^{\infty } \sum _{p=0}^{\infty } 2^{-2-2 j-3 p} e^{-\frac{1}{2} i p \pi } \pi ^{-3-2 p}
   \binom{-\frac{1}{2}}{j} \binom{-1-2 j}{p} \cos \left(\frac{p \pi }{2}\right) \zeta (2+p) (-p)_p
\end{multline}
\end{example}
\begin{example}
\begin{multline}
\int_0^{\infty } \frac{e^{-2 x} \left(-4 \pi +x \cot \left(\frac{x}{4 \pi }\right)\right) \sin (x)}{x^2} \,
   dx\\
=-\sum _{j=0}^{\infty } \sum _{p=0}^{\infty } (-1)^{-p} 2^{-4-2 j-3 p} \left(1+(-1)^p\right) \pi ^{-3-2 p}
   \binom{-1}{j} \binom{-2-2 j}{p} (-p)_p \zeta (2+p)
\end{multline}
\end{example}
\begin{example}
\begin{multline}
\int_0^{\infty } \frac{e^{-2 x} J_{\frac{1}{3}}(x) \left(-4 \pi +x \cot \left(\frac{x}{4 \pi
   }\right)\right)}{x^{5/3}} \, dx\\
=-\sum _{j=0}^{\infty } \sum _{p=0}^{\infty } (-1)^{-p} 2^{-\frac{10}{3}-2 j-3 p}
   \left(1+(-1)^p\right) \pi ^{-\frac{7}{2}-2 p} \binom{-\frac{5}{6}}{j} \binom{-\frac{5}{3}-2 j}{p} \Gamma
   \left(\frac{5}{6}\right) (-p)_p \zeta (2+p)
\end{multline}
\end{example}
\begin{example}
From Eq. (7.13.14) in \cite{titch}. The theta function.
\begin{multline}
\int_0^{\infty } e^{-m x} x^n \vartheta _3\left(0,e^{-\pi ^2 x}\right) (\log (a)-x)^k \, dx,\\
=\sum
   _{p=0}^{\infty } m^{-1-n-p} (1+k-p)_p \binom{-1-n}{p} n! \log ^{k-p}(a)\\
+\sum _{j=1}^{\infty } \sum
   _{p=0}^{\infty } 2 \left(m+j^2 \pi ^2\right)^{-1-n-p} (1+k-p)_p \binom{-1-n}{p} n! \log
   ^{k-p}(a)
\end{multline}
\end{example}
\begin{example}
\begin{multline}
\int_0^{\infty } \frac{e^{-x} \left(-4 \pi +x \cot \left(\frac{x}{4 \pi }\right)\right) \vartheta
   _3\left(0,e^{-\pi ^2 x}\right)}{x^{3/2}} \, dx\\
=-\sum _{p=0}^{\infty } (-1)^{-p} 2^{-3-2 p} \left(1+(-1)^p\right) \pi
   ^{-\frac{5}{2}-2 p} \binom{-\frac{3}{2}}{p} (-p)_p \zeta (2+p)\\
-\sum _{j=1}^{\infty } \sum _{p=0}^{\infty } (-1)^{-p}
   2^{-2-2 p} \left(1+(-1)^p\right) \pi ^{-\frac{5}{2}-2 p} \left(1+j^2 \pi ^2\right)^{-\frac{3}{2}-p}
   \binom{-\frac{3}{2}}{p} (-p)_p \zeta (2+p)
\end{multline}
\end{example}
\begin{example}
\begin{multline}
\int_0^{\infty } \frac{e^{-x} \left(-4 \pi +x \cot \left(\frac{x}{4 \pi }\right)\right) \vartheta
   _3\left(0,e^{-\pi ^2 x}\right)}{x^{4/3}} \, dx\\
=-\sum _{p=0}^{\infty } (-1)^{-p} 2^{-2-2 p} \left(1+(-1)^p\right) \pi
   ^{-3-2 p} \binom{-\frac{5}{3}}{p} \frac{2}{3}! (-p)_p \zeta (2+p)\\
-\sum _{j=1}^{\infty } \sum _{p=0}^{\infty }
   (-1)^{-p} 2^{-1-2 p} \left(1+(-1)^p\right) \pi ^{-3-2 p} \left(1+j^2 \pi ^2\right)^{-\frac{5}{3}-p}
   \binom{-\frac{5}{3}}{p} \frac{2}{3}! (-p)_p \zeta (2+p)
\end{multline}
\end{example}
\begin{example}
\begin{multline}
\int_0^{\infty } e^{-m x} x^{-2+n} \left(-4+x z \cot \left(\frac{x z}{4}\right)\right) \vartheta
   _3\left(0,e^{-\pi ^2 x}\right) \, dx\\
=-\sum _{p=0}^{\infty } (-1)^{-p} 2^{-2-2 p} \left(1+(-1)^p\right) m^{-1-n-p}
   \left(\frac{z}{\pi }\right)^{2+p} \binom{-1-n}{p} n! (-p)_p \zeta (2+p)\\
-\sum _{j=1}^{\infty } \sum _{p=0}^{\infty }
   (-1)^{-p} 2^{-1-2 p} \left(1+(-1)^p\right) \left(m+j^2 \pi ^2\right)^{-1-n-p} \left(\frac{z}{\pi }\right)^{2+p}
   \binom{-1-n}{p} n! (-p)_p \zeta (2+p)
\end{multline}
\end{example}
\begin{example}
\begin{multline}
\int_0^{\infty } e^{-m x} \vartheta _4\left(0,e^{-\pi ^2 x}\right) (\log (a)-x)^k \, dx\\
=-\frac{\Gamma (1+k,-m
   \log (a))}{(-1)^{k+1} a^m m^{k+1}}-\sum _{n=1}^{\infty } \frac{2 (-1)^n \Gamma \left(1+k,-\left(\left(m+n^2 \pi
   ^2\right) \log (a)\right)\right)}{(-1)^{k+1} a^{m+(n \pi )^2} \left(m+n^2 \pi ^2\right)^{k+1}}
\end{multline}
\end{example}
\begin{example}
\begin{multline}
\int_0^{\infty } e^{-2 x} \vartheta _4\left(0,e^{-\pi ^2 x}\right) \log \left(-\frac{1}{2} (i \pi )-x\right) \,
   dx\\
=\frac{1}{4} \left(-i \pi  \left(-2+\sqrt{2} \text{csch}\left(\sqrt{2}\right)\right)+2 \text{Ei}(-i \pi
   )-\sqrt{2} \text{csch}\left(\sqrt{2}\right) \log \left(\frac{4}{\pi ^2}\right)\right)\\
-\sum _{n=1}^{\infty } \frac{2(-1)^n e^{\frac{1}{2} i n^2 \pi ^3} \Gamma \left(0,\frac{1}{2} i \pi  \left(2+n^2 \pi ^2\right)\right)}{2+n^2 \pi^2}
\end{multline}
\end{example}
\begin{example}
\begin{multline}
\int_0^{\infty } e^{-m x} \vartheta _4\left(0,e^{-\pi ^2 x}\right) \log \left(-\frac{i \pi }{a}-x\right) \,
   dx\\
=\frac{e^{\frac{i m \pi }{a}} \Gamma \left(0,\frac{i m \pi }{a}\right)+\sqrt{m} \text{csch}\left(\sqrt{m}\right)
   \left(-i \pi +\log \left(\frac{i \pi }{a}\right)\right)}{m}\\
+\sum _{n=1}^{\infty } \frac{2 (-1)^n e^{\frac{i \pi 
   \left(m+n^2 \pi ^2\right)}{a}} \Gamma \left(0,\frac{i \pi  \left(m+n^2 \pi ^2\right)}{a}\right)}{m+n^2 \pi
   ^2}
\end{multline}
where $Re(a)>0$
\end{example}
\begin{example}
\begin{multline}
\int_0^{\infty } e^{-m x} \vartheta _4\left(0,e^{-\pi ^2 x}\right) \log \left(-\frac{i \pi }{a}-x\right) \,
   dx\\
=\frac{e^{\frac{i m \pi }{a}} \Gamma \left(0,\frac{i m \pi }{a}\right)+\sqrt{m} \text{csch}\left(\sqrt{m}\right)
   \left(i \pi +\log \left(\frac{i \pi }{a}\right)\right)}{m}\\
+\sum _{n=1}^{\infty } \frac{2 (-1)^n e^{\frac{i \pi 
   \left(m+n^2 \pi ^2\right)}{a}} \Gamma \left(0,\frac{i \pi  \left(m+n^2 \pi ^2\right)}{a}\right)}{m+n^2 \pi
   ^2}
\end{multline}
where $Re(a)<0$
\end{example}
\begin{example}
\begin{multline}
\int_0^{\infty } e^{-m x} \vartheta _4\left(0,e^{-\pi ^2 x}\right) \log \left(\frac{\pi ^2}{a^2}+x^2\right) \,
   dx\\
=\sum _{n=1}^{\infty } \frac{2 (-1)^n e^{-\frac{i \pi  \left(m+n^2 \pi ^2\right)}{a}} \left(\Gamma
   \left(0,-\frac{i \pi  \left(m+n^2 \pi ^2\right)}{a}\right)+e^{\frac{2 i \pi  \left(m+n^2 \pi ^2\right)}{a}} \Gamma
   \left(0,\frac{i \pi  \left(m+n^2 \pi ^2\right)}{a}\right)\right)}{m+n^2 \pi ^2}\\
+\frac{e^{-\frac{i m \pi }{a}} \Gamma
   \left(0,-\frac{i m \pi }{a}\right)+e^{\frac{i m \pi }{a}} \Gamma \left(0,\frac{i m \pi }{a}\right)+\sqrt{m}
   \text{csch}\left(\sqrt{m}\right) \log \left(\frac{\pi ^2}{a^2}\right)}{m}
\end{multline}
\end{example}
\begin{example}
\begin{multline}
\int_0^{\infty } e^{-i x} \vartheta _4\left(0,e^{-\pi ^2 x}\right) \log \left(1+x^2\right) \, dx\\
=\sum
   _{n=1}^{\infty } \frac{2 (-1)^n e^{-1-i n^2 \pi ^2} \left(e^2 \Gamma \left(0,1-i n^2 \pi ^2\right)+e^{2 i n^2 \pi
   ^2} \Gamma \left(0,-1+i n^2 \pi ^2\right)\right)}{i+n^2 \pi ^2}\\
-i (e E_1(1)+\text{Subfactorial}[-1])
\end{multline}
\end{example}
\begin{example}
\begin{multline}
\int_0^{\infty } e^{-x} \vartheta _4\left(0,e^{-\pi ^2 x}\right) \log \left(1+x^2\right) \, dx\\
=\sum
   _{n=1}^{\infty } \frac{2 (-1)^n e^{-i \left(1+n^2 \pi ^2\right)} \left(\Gamma \left(0,-i \left(1+n^2 \pi
   ^2\right)\right)+e^{2 i \left(1+n^2 \pi ^2\right)} \Gamma \left(0,i \left(1+n^2 \pi ^2\right)\right)\right)}{1+n^2\pi ^2}\\
+e^{-i} \Gamma (0,-i)+e^i \Gamma (0,i)
\end{multline}
\end{example}
\begin{example}
\begin{multline}
\int_0^{\infty } e^{-m x} \vartheta _4\left(0,e^{-\pi ^2 x}\right) \tan ^{-1}\left(\frac{a x}{\pi }\right) \,
   dx\\
=\frac{i e^{-\frac{i m \pi }{a}} \left(-\Gamma \left(0,-\frac{i m \pi }{a}\right)+e^{\frac{2 i m \pi }{a}} \Gamma
   \left(0,\frac{i m \pi }{a}\right)\right)}{2 m}\\
+\sum _{n=1}^{\infty } \frac{i (-1)^n e^{-\frac{i \pi  \left(m+n^2 \pi
   ^2\right)}{a}} \left(-\Gamma \left(0,-\frac{i \pi  \left(m+n^2 \pi ^2\right)}{a}\right)+e^{\frac{2 i \pi 
   \left(m+n^2 \pi ^2\right)}{a}} \Gamma \left(0,\frac{i \pi  \left(m+n^2 \pi ^2\right)}{a}\right)\right)}{m+n^2 \pi^2}
\end{multline}
\end{example}
\begin{example}
\begin{multline}
\int_0^{\infty } e^{-2 i x} \tan ^{-1}(x) \vartheta _4\left(0,e^{-\pi ^2 x}\right) \, dx\\
=\frac{1}{4} e^2
   \left(\frac{\Gamma (0,-2)}{e^4}-\Gamma (0,2)\right)\\
+\sum _{n=1}^{\infty } \frac{i (-1)^n e^{-i \left(2 i+n^2 \pi
   ^2\right)} \left(-\Gamma \left(0,-i \left(2 i+n^2 \pi ^2\right)\right)+e^{2 i \left(2 i+n^2 \pi ^2\right)} \Gamma
   \left(0,i \left(2 i+n^2 \pi ^2\right)\right)\right)}{2 i+n^2 \pi ^2}
\end{multline}
\end{example}
\begin{example}
\begin{multline}
\int_0^{\infty } e^{-x} \tan ^{-1}(2 x) \vartheta _4\left(0,e^{-\pi ^2 x}\right) \, dx\\
=\frac{1}{2} i
   e^{-\frac{i}{2}} \left(-\Gamma \left(0,-\frac{i}{2}\right)+e^i \Gamma \left(0,\frac{i}{2}\right)\right)\\
+\sum
   _{n=1}^{\infty } \frac{i (-1)^n e^{-\frac{1}{2} i \left(1+n^2 \pi ^2\right)} \left(-\Gamma \left(0,-\frac{1}{2} i
   \left(1+n^2 \pi ^2\right)\right)+e^{i \left(1+n^2 \pi ^2\right)} \Gamma \left(0,\frac{1}{2} i \left(1+n^2 \pi
   ^2\right)\right)\right)}{1+n^2 \pi ^2}
\end{multline}
\end{example}
\begin{example}
\begin{multline}
\int_0^{\infty } \frac{e^{-m x} \pi  \vartheta _4\left(0,e^{-\pi ^2 x}\right)}{\pi ^2+a^2 x^2} \, dx\\
=\sum
   _{n=1}^{\infty } -\frac{i (-1)^n e^{-\frac{i \pi  \left(m+n^2 \pi ^2\right)}{a}} \left(\Gamma \left(0,-\frac{i \pi 
   \left(m+n^2 \pi ^2\right)}{a}\right)-e^{\frac{2 i \pi  \left(m+n^2 \pi ^2\right)}{a}} \Gamma \left(0,\frac{i \pi 
   \left(m+n^2 \pi ^2\right)}{a}\right)\right)}{a}\\
-\frac{i \left(e^{-\frac{i m \pi }{a}} \Gamma \left(0,-\frac{i m \pi
   }{a}\right)-e^{\frac{i m \pi }{a}} \Gamma \left(0,\frac{i m \pi }{a}\right)\right)}{2 a}
\end{multline}
\end{example}
\begin{example}
\begin{multline}
\int_0^{\infty } \frac{e^{-m x} x^2 \vartheta _4\left(0,e^{-\pi ^2 x}\right)}{\left(\pi ^2+a^2 x^2\right)^2} \,
   dx
=-\sum _{n=1}^{\infty } \frac{(-1)^n e^{-\frac{i \pi  \left(m+n^2 \pi ^2\right)}{a}} }{2 a^4 \pi }\\\times
\left(\left(i a+m \pi +n^2
   \pi ^3\right) \Gamma \left(0,-\frac{i \pi  \left(m+n^2 \pi ^2\right)}{a}\right)\right. \\ \left.+e^{\frac{2 i \pi  \left(m+n^2 \pi
   ^2\right)}{a}} \left(-i a+m \pi +n^2 \pi ^3\right) \Gamma \left(0,\frac{i \pi  \left(m+n^2 \pi
   ^2\right)}{a}\right)\right)\\
+\frac{e^{-\frac{i m \pi }{a}} \left((-i a-m \pi ) \Gamma \left(0,-\frac{i m
   \pi }{a}\right)+i e^{\frac{2 i m \pi }{a}} (a+i m \pi ) \Gamma \left(0,\frac{i m \pi }{a}\right)\right)}{4 a^4 \pi
   }
\end{multline}
\end{example}
\begin{example}
\begin{multline}
\int_0^{\infty } \frac{e^{-m x} \vartheta _4\left(0,e^{-\pi ^2 x}\right) \log (-x+\log (a))}{-x+\log (a)} \,
   dx\\
=\frac{1}{12} a^{-m} \left(-6 \gamma ^2-\pi ^2-12 m \, _3F_3(1,1,1;2,2,2;m \log (a)) \log (a)\right. \\ \left.
+12 \text{Chi}(m \log
   (a)) (-i \pi -\log (m)+\log (-m \log (a)))\right. \\ \left.
+6 \left(\log ^2(-m \log (a))+2 (i \pi +\log (m)) \log (m \log (a))\right.\right. \\ \left.\left.
-2 \log(-m \log (a)) (\gamma +i \pi +\log (m)+\log (m \log (a)))\right)\right. \\ \left.
+12 (-i \pi -\log (m)+\log (-m \log (a)))
   \text{Shi}(m \log (a))\right)\\
-\sum _{n=1}^{\infty } \frac{1}{6} (-1)^n a^{-m-n^2 \pi ^2} \left(6 \gamma ^2+\pi
   ^2+12 \left(m+n^2 \pi ^2\right)\right. \\ \left. 
\, _3F_3\left(1,1,1;2,2,2;\left(m+n^2 \pi ^2\right) \log (a)\right) \log (a)\right. \\ \left.
+12\text{Chi}\left(\left(m+n^2 \pi ^2\right) \log (a)\right) \left(i \pi +\log \left(m+n^2 \pi ^2\right)\right.\right. \\ \left.\left.
-\log\left(-\left(\left(m+n^2 \pi ^2\right) \log (a)\right)\right)\right)+6 \left(2 \left(\gamma +i \pi +\log \left(m+n^2 \pi ^2\right)\right)\right.\right. \\ \left.\left.
-\log \left(-\left(\left(m+n^2 \pi ^2\right) \log (a)\right)\right)\right) \log
   \left(-\left(\left(m+n^2 \pi ^2\right) \log (a)\right)\right)\right. \\ \left.
+12 \left(-i \pi -\log \left(m+n^2 \pi ^2\right)+\log
   \left(-\left(\left(m+n^2 \pi ^2\right) \log (a)\right)\right)\right) \log \left(\left(m+n^2 \pi ^2\right) \log
   (a)\right)\right. \\ \left.
+12 \left(i \pi +\log \left(m+n^2 \pi ^2\right)-\log \left(-\left(\left(m+n^2 \pi ^2\right) \log
   (a)\right)\right)\right) \text{Shi}\left(\left(m+n^2 \pi ^2\right) \log (a)\right)\right)
\end{multline}
\end{example}
\begin{example}
From Eq. (1) in \cite{glasser}
\begin{multline}
\int_0^{\infty } \frac{e^{-m x (i \pi +x)} \cosh (x) \left(\left(-i \pi  x+\log (a)-x^2\right)^k+e^{2 i m \pi  x} \left(i \pi  x+\log (a)-x^2\right)^k\right)}{1+\alpha ^4+2 \alpha ^2 \cosh (2 x)} \,
   dx\\
=\frac{e^{-\frac{1}{4} \left(m \pi ^2\right)-m \log ^2(\alpha )} \pi  \left(\log (a)+\log \left(e^{-\frac{\pi ^2}{4}-\log ^2(\alpha )}\right)\right)^k}{2 \alpha  \left(1+\alpha ^2\right)}
\end{multline}
\end{example}
\begin{example}
\begin{multline}
\int_0^{\infty } \frac{e^{-m x (i \pi +x)} \cosh (x) \left(e^{2 i m \pi  x} \Phi \left(-1,1,1-i \pi  x+x^2\right)+\Phi \left(-1,1,1+i \pi  x+x^2\right)\right)}{1+\alpha ^4+2 \alpha ^2 \cosh (2 x)} \,
   dx\\
=\frac{e^{-\frac{1}{4} \left(m \pi ^2\right)-m \log ^2(\alpha )} \pi  \Phi \left(-1,1,1-\log \left(e^{-\frac{\pi ^2}{4}-\log ^2(\alpha )}\right)\right)}{2 \alpha  \left(1+\alpha ^2\right)}
\end{multline}
\end{example}
\begin{example}
\begin{multline}
\int_0^{\infty } \frac{e^{-m x (i \pi +x)} \cosh (x) \left(e^{2 i m \pi  x} \psi ^{(1)}\left(1-i \pi  x+x^2\right)+\psi ^{(1)}\left(1+i \pi  x+x^2\right)\right)}{1+\alpha ^4+2 \alpha ^2 \cosh (2 x)} \,
   dx\\
=\frac{e^{-\frac{1}{4} m \left(\pi ^2+4 \log ^2(\alpha )\right)} \pi  \psi ^{(1)}\left(1-\log \left(e^{-\frac{\pi ^2}{4}-\log ^2(\alpha )}\right)\right)}{2 \alpha  \left(1+\alpha ^2\right)}
\end{multline}
\end{example}
\begin{example}
\begin{equation}
\int_0^{\infty } \frac{\cosh (x) \log \left(a^2+\left(2 a+\pi ^2\right) x^2+x^4\right)}{1+\alpha ^4+2 \alpha ^2 \cosh (2 x)} \, dx=\frac{\pi  \log \left(a+\frac{\pi ^2}{4}+\log ^2(\alpha )\right)}{2
   \alpha  \left(1+\alpha ^2\right)}
\end{equation}
\end{example}
\begin{example}
\begin{multline}
\int_0^{\infty } \frac{e^{-m x^2} \cosh (x) \left(\left(-a^2+(\pi -x) x^2 (\pi +x)\right) \cos (m \pi  x)+2 \pi  x^3 \sin (m \pi  x)\right)}{\left(-a^2+\left(-2 i a+\pi ^2\right) x^2+x^4\right)
   \left(-a^2+\left(2 i a+\pi ^2\right) x^2+x^4\right) \left(1+\alpha ^4+2 \alpha ^2 \cosh (2 x)\right)} \, dx\\
=-\frac{4 e^{-\frac{1}{4} m \left(\pi ^2+4 \log ^2(\alpha )\right)} \pi }{\alpha  \left(1+\alpha
   ^2\right) \left(16 a^2+\pi ^4+8 \pi ^2 \log ^2(\alpha )+16 \log ^4(\alpha )\right)}
\end{multline}
\end{example}
\begin{example}
\begin{multline}
\int_0^{\infty } \frac{e^{-m x (i \pi +x)} \left(-1+a-i \pi  x+x^2+e^{2 i m \pi  x} \left(-1+a+i \pi  x+x^2\right)\right) \cosh (x)}{\left(-1+a-i \pi  x+x^2\right) \left(-1+a+i \pi  x+x^2\right)
   \left(1+\alpha ^4+2 \alpha ^2 \cosh (2 x)\right)} \, dx\\
=\frac{2 e^{-\frac{1}{4} m \left(\pi ^2+4 \log ^2(\alpha )\right)} \pi }{\left(\alpha +\alpha ^3\right) \left(-4+4 a+\pi ^2+4 \log ^2(\alpha
   )\right)}
\end{multline}
\end{example}
\begin{example}
\begin{multline}
\int_0^{\infty } \frac{e^{-m x (i \pi +x)} }{\left((1+a)^2+\left(2+2 a+\pi ^2\right) x^2+x^4\right) \left(1+\alpha ^4+2 \alpha ^2 \cosh (2 x)\right)}\\\times
\left(-(-a+(-i \pi -x) x)^{1+k} \left(1+a-i \pi  x+x^2\right)\right. \\ \left.
-e^{2 i m \pi  x} (-a+i (\pi +i x) x)^{1+k} \left(1+a+i \pi  x+x^2\right)\right) \cosh
   (x) \, dx\\
=-\frac{2^{-1-2 k} e^{-\frac{1}{4} m \left(\pi ^2+4 \log ^2(\alpha )\right)} \pi 
   \left(-4 a-\pi ^2-4 \log ^2(\alpha )\right)^{1+k}}{\left(\alpha +\alpha ^3\right) \left(4+4 a+\pi ^2+4 \log ^2(\alpha )\right)}
\end{multline}
\end{example}
\section{Application of an Adams theorem}
In this section we apply the theorem given in Eq. (6.706) in \cite{adams} to the formulae cited.
\begin{example}
From Eq. (7.704) in \cite{adams}
\begin{multline}
\int_0^1 \frac{(-1)^k a^{-m+t^{-\beta }} (1-t)^{-1+\gamma } t^{-1+\alpha +k \beta } \Gamma
   \left(1+k,\left(-m+t^{-\beta }\right) \log (a)\right)}{\left(-1+m t^{\beta }\right)^{k+1}} \, dt\\
=-\sum
   _{n=0}^{\infty } \sum _{j=0}^n \frac{m^{-j+n} (1-j+k)_j \binom{n}{j} (-1+\gamma )! \log ^{-j+k}(a)}{(\alpha +n
   \beta ) (1+\alpha +n \beta )_{-1+\gamma }}
  \end{multline}
   where $Im(m)\geq 0$
\end{example}
\begin{example}
\begin{multline}
\int_0^1 (2 i)^k e^{i m t^{\beta }} (1-t)^{-1+\gamma } t^{-1+\alpha +k \beta } \Phi \left(-e^{2 i m t^{\beta
   }},-k,\frac{1}{2}-\frac{1}{2} i t^{-\beta } \log (a)\right) \, dt\\
=\sum _{n=0}^{\infty } \sum _{j=0}^n \frac{(4 i)^n
   m^{-j+n} \Gamma (\gamma ) (1-j+k)_j \binom{n}{j} \log ^{-j+k}(a) \left(\zeta \left(-n,\frac{1}{4}\right)-\zeta
   \left(-n,\frac{3}{4}\right)\right)}{(\alpha +n \beta ) n! (1+\alpha +n \beta )_{-1+\gamma }}
\end{multline}
\end{example}
\begin{example}
\begin{multline}
\int_0^1 t^{a-1} (1-t)^{k-1} \left(t^{b/2} \sqrt{x}\right)^{-v} J_v\left(2 t^{b/2} \sqrt{x}\right) \, dt\\
=\sum
   _{n=0}^{\infty } \frac{(-1)^n \left((k-1)! x^n\right)}{(n! \Gamma (v+n+1)) \left((a+b n) (1+a+b
   n)_{-1+k}\right)}
\end{multline}
\end{example}
\begin{example}
\begin{multline}
\int_0^1 t^{a-1} (1-t)^{k-1} \, _2F_1\left(\alpha ,\beta ;\gamma ;t^b x\right) \, dt
=\sum _{n=0}^{\infty }
   \frac{\left((\alpha )_n (\beta )_n\right) \left((k-1)! x^n\right)}{\left((\gamma )_n n!\right) \left((a+b n) (1+a+b
   n)_{-1+k}\right)}
\end{multline}
\end{example}
\begin{example}
\begin{multline}
\int_0^1 \frac{t^{a-1} (1-t)^{k-1} \, _1F_2\left(\frac{1}{2};1-v,1+v;-4 t^b x\right)}{\Gamma (1-v) \Gamma (1+v)}
   \, dt\\
=\sum _{n=0}^{\infty } \frac{\left((-1)^n (2 n)!\right) \left((k-1)! x^n\right)}{\left((n!)^2 \Gamma (n+v+1)
   \Gamma (n-v+1)\right) \left((a+b n) (1+a+b n)_{-1+k}\right)}
\end{multline}
\end{example}
\begin{example}
From Eq. (6.7.1.1) in \cite{brychkov}
\begin{multline}
\int_0^1 t^{a-1} (1-t)^{k-1} e^{\frac{t^b x z^2}{2 \left(1+2 t^b x\right)}} \left(1+2 t^b x\right)^{\frac{1}{2}
   (-1-v)} D_v\left(\frac{z}{\sqrt{1+2 t^b x}}\right) \, dt\\
=\sum _{n=0}^{\infty } \frac{D_{2 n+v}(z) \left((k-1)!
   x^n\right)}{n! \left((a+b n) (1+a+b n)_{-1+k}\right)}
\end{multline}
\end{example}
\begin{example}
From Eq. (6.12.1.1) in \cite{brychkov}
\begin{multline}
\int_0^1 t^{a-1} (1-t)^{k-1} e^{t^b x z} \cos \left(t^b x \sqrt{1-z^2}\right) \, dt=\sum _{n=0}^{\infty }
   \frac{T_n(z) \left((k-1)! x^n\right)}{n! \left((a+b n) (1+a+b n)_{-1+k}\right)}
\end{multline}
\end{example}
\begin{example}
From Eq. (6.14.1.4) in \cite{brychkov}
\begin{multline}
\int_0^1 t^{a-1} (1-t)^{k-1} e^{t^b x} \left(t^b x z\right)^{-\frac{\lambda }{2}} J_{\lambda }\left(2 \sqrt{t^b x
   z}\right) \Gamma (1+\lambda ) \, dt\\
=\sum _{n=0}^{\infty } \frac{L_n^{\lambda }(z) \left((k-1)! x^n\right)}{(1+\lambda
   )_n \left((a+b n) (1+a+b n)_{-1+k}\right)}
\end{multline}
\end{example}
\begin{example}
From Eq. (6.14.1.9) in \cite{brychkov}
\begin{multline}
\int_0^1 t^{a-1} (1-t)^{k-1} e^{2 t^b x} I_0\left(2 \sqrt{t^b x \left(t^b x-z\right)}\right) \, dt=\sum
   _{n=0}^{\infty } \frac{L_n^n(z) \left((k-1)! x^n\right)}{n! \left((a+b n) (1+a+b n)_{-1+k}\right)}
\end{multline}
\end{example}
\begin{example}
The Bell number $B_n$.
\begin{multline}
\int_0^1 t^{a-1} (1-t)^{k-1} e^{t^b x} B_j\left(t^b x\right) \, dt=\sum _{n=0}^{\infty } \frac{n^j \left((k-1)!
   x^n\right)}{n! \left((a+b n) (1+a+b n)_{-1+k}\right)}
\end{multline}
\end{example}
\begin{example}
\begin{multline}
\int_0^1 \frac{t^{a-1} (1-t)^{k-1} \,
   _2F_3\left(\frac{1}{2}+\frac{u}{2}+\frac{v}{2},1+\frac{u}{2}+\frac{v}{2};1+u,1+v,1+u+v;-t^b x\right)}{\Gamma (1+u)
   \Gamma (1+v)} \, dt\\
=\sum _{n=0}^{\infty } \frac{\left(\left(-\frac{1}{4}\right)^n (1+n+u+v)_n\right) \left((k-1)!
   x^n\right)}{(n! \Gamma (1+n+u) \Gamma (1+n+v)) \left((a+b n) (1+a+b n)_{-1+k}\right)}
\end{multline}
\end{example}
\begin{example}
From Eq. (10.8.3) in \cite{dlmf}
\begin{multline}
\int_0^1 t^{a-1} (1-t)^{k-1} 2^{u+v} \left(t^b x\right)^{\frac{1}{2} (-u-v)} J_u\left(\sqrt{t^b x}\right)
   J_v\left(\sqrt{t^b x}\right) \, dt\\
=\sum _{n=0}^{\infty } \frac{\left(\left(-\frac{1}{4}\right)^n (1+n+u+v)_n\right)
   \left((k-1)! x^n\right)}{(n! \Gamma (1+n+u) \Gamma (1+n+v)) \left((a+b n) (1+a+b n)_{-1+k}\right)}
\end{multline}
\end{example}
\begin{example}
From Eq. (5.2.16.6) in \cite{prud1}. Complete elliptic integral of the first kind $K(m)$.
\begin{multline}
\frac{2 }{\pi }\int_0^1 t^{a-1} (1-t)^{k-1} K\left(16 t^b x\right) \, dt=\sum _{n=0}^{\infty } \frac{((2 n)!)^2
   \left((k-1)! x^n\right)}{(n!)^4 \left((a+b n) (1+a+b n)_{-1+k}\right)}
\end{multline}
\end{example}
\begin{example}
From Eq. (5.2.16.7) in \cite{prud1}. The complete elliptic integral E(m).
\begin{multline}
-\frac{2 }{\pi }\int_0^1 t^{a-1} (1-t)^{k-1} E\left(16 t^b x\right) \, dt=\sum _{n=0}^{\infty } \frac{((2 n)!)^2
   \left((k-1)! x^n\right)}{\left((-1+2 n) (n!)^4\right) \left((a+b n) (1+a+b n)_{-1+k}\right)}
\end{multline}
\end{example}
\begin{example}
From Eq. (5.3.1.5(i)) in\cite{prud2}.
\begin{multline}
\frac{1}{2}\int_0^1  t^{a-1} (1-t)^{k-1} \left(-\pi  \sqrt{t^b x} \cot \left(\pi  \sqrt{t^b x}\right)\right) \,
   dt=\sum _{n=0}^{\infty } \frac{\zeta (2 n) \left((k-1)! x^n\right)}{(a+b n) (1+a+b n)_{-1+k}}
\end{multline}
\end{example}
\begin{example}
From Eq. (5.3.1.5(ii)) in\cite{prud2}. 
\begin{multline}
-\frac{1}{2} \pi  \int_0^1 t^{a-1} (1-t)^{k-1} \left(\sqrt{t^b x} \coth \left(\pi  \sqrt{t^b x}\right)\right) \,
   dt\\
=(k-1)! \sum _{n=0}^{\infty } \frac{\left((-1)^n \zeta (2 n)\right) x^n}{(a+b n) (1+a+b n)_{-1+k}}
\end{multline}
\end{example}
\begin{example}
From Eq. (5.3.1.11) in \cite{prud2}. 
\begin{multline}
\int_0^1 t^{a-1} (1-t)^{k-1} \zeta \left(\alpha ,1-t^b x\right) \, dt=\sum _{n=0}^{\infty } \frac{\left((\alpha
   )_n \zeta (n+\alpha )\right) \left((k-1)! x^n\right)}{n! \left((a+b n) (1+a+b n)_{-1+k}\right)}
\end{multline}
\end{example}
\begin{example}
From Eq. (5.3.1.13) in \cite{prud2}.
\begin{multline}
\int_0^1 \frac{t^{a-1} (1-t)^{k-1} \left((-1)^j \psi ^{(-1+j)}\left(1-t^b x\right)\right)}{(-1+j)!} \, dt=\sum
   _{n=0}^{\infty } \frac{(\binom{-1+j+n}{n} \zeta (j+n)) \left((k-1)! x^n\right)}{(a+b n) (1+a+b n)_{-1+k}}
\end{multline}
\end{example}
\begin{example}
From Eq. (5.7.6.1) in \cite{prud2} 
\begin{multline}
\int_0^1 t^{a-1} (1-t)^{k-1} \left(z^{v/2} \left(-2 t^b x+z\right)^{-\frac{v}{2}} J_v\left(\sqrt{-2 t^b x
   z+z^2}\right)\right) \, dt\\
=\sum _{n=0}^{\infty } \frac{J_{n+v}(z) \left((k-1)! x^n\right)}{n! \left((a+b n) (1+a+b
   n)_{-1+k}\right)}
\end{multline}
\end{example}
\begin{example}
From Eq. (5.9.1.1) in \cite{prud2}. 
\begin{multline}
\int_0^1 \frac{t^{a-1} (1-t)^{k-1} \left(e^{\frac{1}{2} \left(\frac{t^{-b}}{x}+t^b x\right) z} \pi 
   \text{erfc}\left(\sqrt{\frac{t^{-b} z}{2 x}}+\frac{\sqrt{t^b x z}}{\sqrt{2}}\right)\right)}{2 \sqrt{t^b x}} \,
   dt\\
=\sum _{n=0}^{\infty } \frac{\left((-1)^n K_{\frac{1}{2}+n}(z)\right) \left((k-1)! x^n\right)}{(a+b n) (1+a+b
   n)_{-1+k}}
\end{multline}
\end{example}
\begin{example}
From Eq. (6.703) in \cite{adams}
\begin{equation}
\int_0^1 \frac{t^{-1+m} \log ^k\left(e^a t\right)}{1+t^b} \, dt=-\sum _{n=0}^{\infty } \frac{(-1)^n \Gamma
   (1+k,-a (m+b n))}{e^{a (m+b n)} (-m-b n)^{k+1}}
\end{equation}
\end{example}
\begin{example}
\begin{multline}
\int_0^1 \frac{t^{-1+m}}{\left(1+t^b\right) \left(a^2-\log ^2(t)\right)} \, dt\\
=\sum _{n=0}^{\infty } \frac{(-1)^n
   e^{-a (m+b n)} \left(-\Gamma (0,-a (m+b n))+e^{2 a (m+b n)} \Gamma (0,a (m+b n))\right)}{2 a}
\end{multline}
\end{example}
\begin{example}
\begin{multline}
\int_0^1 \frac{t^{-1+m} \log (t)}{\left(1+t^b\right) \left(-a^2+\log ^2(t)\right)} \, dt\\
=-\frac{1}{2} \sum _{n=0}^{\infty }
   (-1)^n e^{-a (m+b n)} \left(\Gamma (0,-a (m+b n))+e^{2 a (m+b n)} \Gamma (0,a (m+b n))\right)
\end{multline}
\end{example}
\begin{example}
\begin{multline}
\int_0^1 \frac{t^{-1+m}}{\left(1+t^b\right) (a+\log (t))^2} \, dt=\sum _{n=0}^{\infty } (-1)^n
   \left(-\frac{1}{a}-e^{-a (m+b n)} (m+b n) \Gamma (0,-a (m+b n))\right)
\end{multline}
\end{example}
\begin{example}
\begin{multline}
\int_0^1 \frac{t^{-1+m} \log \left(a^2+\log ^2(t)\right)}{1+t^b} \, dt\\
=\sum _{n=0}^{\infty } \frac{(-1)^n e^{-i a
   (m+b n)} \left(E_1(-i a (m+b n))+e^{2 i a (m+b n)} E_1(i a (m+b n))+2 e^{i a (m+b n)} \log (a)\right)}{m+b n}
\end{multline}
\end{example}
\begin{example}
\begin{multline}
\int_0^1 \frac{t^{-1+m} \tanh ^{-1}\left(\frac{\log (t)}{a}\right)}{1+t^b} \, dt\\
=\sum _{n=0}^{\infty }
   \frac{(-1)^n  \left(E_1(-a (m+b n))-e^{a (m+b n)} \left(e^{a (m+b n)} E_1(a (m+b n))+\log (-a)-\log(a)\right)\right)}{2 (m+b n)e^{a (m+b n)}}\\
+\frac{i \pi  \left(-\psi ^{(0)}\left(\frac{m}{2 b}\right)+\psi ^{(0)}\left(\frac{b+m}{2
   b}\right)\right)}{4 b}
\end{multline}
\end{example}
\begin{example}
From Eq. (1) in \cite{brychkov_1}
The modified Bessel function and the Gegenbauer polynomial
\begin{multline}
\int_0^1 t^{a-1} (1-t)^{k-1} 2^v e^{t^b x} \left(t^b x \alpha \right)^{-v} I_v\left(t^b x \alpha
   \right) \, dt\\
=(k-1)!\sum _{n=0}^{\infty } \frac{2^{-n} (-\alpha )^n C_n^{(-n-v)}\left(\frac{1}{\alpha
   }\right)x^n}{\Gamma (1+n+v) \left((a+b n) (1+a+b n)_{-1+k}\right)}
\end{multline}
\end{example}
\begin{example}
From Eq. (133.1) pp.256 in \cite{rainville}
\begin{multline}
\int_0^1 \frac{t^{a-1} (1-t)^{k-1} \left(I_{\alpha }\left(\sqrt{2} \sqrt{-1+x} \sqrt{t^b x}\right) I_{\beta
   }\left(\sqrt{2} \sqrt{t^b x} \sqrt{1+x}\right) \Gamma (1+\alpha ) \Gamma (1+\beta )\right)}{\left(\sqrt{-1+x}
   \sqrt{t^b x}\right)^{\alpha } \left(\sqrt{t^b x} \sqrt{1+x}\right)^{\beta }} \, dt\\
=\frac{(k-1)! }{2^{\frac{\alpha }{2}+\frac{\beta }{2}}}\sum _{n=0}^{\infty
   } \frac{P_n^{(\alpha ,\beta )}(x) x^n}{\left((1+\alpha )_n (1+\beta )_n\right) \left((a+b n) (1+a+b
   n)_{-1+k}\right)}
\end{multline}
\end{example}
\begin{example}
From Eq. (148.1) in \cite{rainville}. Bateman-Brafman integral.
\begin{multline}
\int_0^1 \frac{t^{a-1} (1-t)^{k-1} \, _2F_1\left(\frac{1}{2},\frac{1}{2}+\frac{z}{2};1;-\frac{4 t^b x}{\left(1-t^b
   x\right)^2}\right)}{1-t^b x} \, dt\\
=(k-1)! \sum _{n=0}^{\infty } \frac{\,
   _3F_2\left(-n,1+n,\frac{1}{2}+\frac{z}{2};1,1;1\right) x^n}{(a+b n) (1+a+b n)_{-1+k}}
\end{multline}
\end{example}
\begin{example}
From Eq. (6.704) in \cite{adams}.
\begin{multline}
\sum _{j=0}^{\infty } \frac{(-1)^j }{a+j}\binom{k-1}{j} \,
   _2F_1\left(1,\frac{a+j}{b};\frac{a+b+j}{b};x\right)\\=(k-1)!\sum _{n=0}^{\infty
   } \frac{x^n }{(a+b n) (1+a+b n)_{k-1}}
\end{multline}
where $Re(k)>2$.
\end{example}
\begin{example}
\begin{multline}
\sum _{j=0}^{\infty } \sum _{l=0}^{\infty } \sum _{p=l}^{\infty } \frac{(-1)^j m^{l-p} \binom{l}{p}
   \binom{-1+\gamma }{j} (1)_l \left(\frac{j+\alpha }{\beta }\right)_l}{(j+\alpha ) l! (k-p)! \left(\frac{j+\alpha
   +\beta }{\beta }\right)_l \log ^p(a)}\\
=\sum _{n=0}^{\infty } \sum _{p=0}^n \frac{m^{n-p} \binom{n}{p} (\gamma-1
   )!}{(\alpha +n \beta ) (k-n)! (1+\alpha +n \beta )_{\gamma-1 } \log ^n(a)}
\end{multline}
where $Re(a)>e^{11\pi}$
\end{example}
\begin{example}
The elliptic integral of the first kind
\begin{equation}
\sum _{j=0}^{\infty } \sum _{p=0}^{\infty } \frac{(-1)^{j+p} k^{2 p}
   x^{1+2 j+2 p} \binom{-\frac{1}{2}}{j} \binom{-\frac{1}{2}}{p}}{1+2 j+2
   p}=\int_0^x \frac{1}{\sqrt{1-t^2}} \frac{1}{\sqrt{1-(k t)^2}} \, dt
\end{equation}
\end{example}
\begin{example}
\begin{multline}
\sum _{j=0}^{\infty } \frac{(-1)^j \left(e^{-k} k\right)^{2 j} \sqrt{\pi
   } \, _2F_1\left(\frac{1}{2},\frac{1}{2}+j;\frac{3}{2}+j;\sin ^2(\phi
   )\right) \sin ^{1+2 j}(\phi )}{(1+2 j) \Gamma \left(\frac{1}{2}-j\right)
   j!}=F\left(\phi \left|e^{-2 k} k^2\right.\right)
\end{multline}
\end{example}
\begin{example}
\begin{multline}
\sum _{j=0}^{\infty } \frac{(-1)^j \left(e^{-k} k\right)^{2 j} \sqrt{\pi
   } \, _2F_1\left(\frac{1}{2},\frac{1}{2}+j;\frac{3}{2}+j;\sin ^2(\phi
   )\right) \sin ^{1+2 j}(\phi )}{(1+2 j) \Gamma \left(\frac{1}{2}-j\right)
   (j+1)!}\\
=-\frac{2e^{2 k} }{k^2}\int_0^k  e^{-2 t} t \left(-F\left(\phi \left|e^{-2
   t} t^2\right.\right)+t F\left(\phi \left|e^{-2 t} t^2\right.\right)\right)
   \, dt
\end{multline}
\end{example}
\begin{example}
\begin{equation}
\int_0^z t^{-1+s} \left(1-(t \alpha )^{\beta }\right)^{\gamma } \, dt=\frac{z^s }{s}\, _2F_1\left(\frac{s}{\beta
   },-\gamma ;\frac{s+\beta }{\beta };z^{\beta } \alpha ^{\beta }\right)
\end{equation}
\end{example}
\begin{example}
\begin{multline}
\int_0^{\phi } \left(1-(a \sin (t))^c\right)^n \, dt\\
=\phi +\sum _{k=1}^n \frac{(-1)^k \binom{n}{k} \,
   _2F_1\left(\frac{1}{2},\frac{1}{2} (1+c k);\frac{1}{2} (3+c k);\sin ^2(\phi )\right) \sin (\phi ) (a \sin (\phi
   ))^{c k}}{1+c k}
\end{multline}
\end{example}
\begin{example}
\begin{multline}
\int_0^{\phi } e^{(a \sin (t))^c} \, dt
=\phi \\+e \sum _{n=0}^{\infty } \sum _{k=1}^n \frac{(-1)^n \left((-1)^k
   \binom{n}{k} \, _2F_1\left(\frac{1}{2},\frac{1}{2} (1+c k);\frac{1}{2} (3+c k);\sin ^2(\phi )\right) \sin (\phi )
   (a \sin (\phi ))^{c k}\right)}{n! (1+c k)}
\end{multline}
\end{example}
\begin{example}
\begin{equation}
\int_0^z t^{-1+s} \left(1-(t \alpha )^{\beta }\right)^{\gamma } \, dt=\frac{z^s }{s}\, _2F_1\left(\frac{s}{\beta
   },-\gamma ;\frac{s+\beta }{\beta };z^{\beta } \alpha ^{\beta }\right)
\end{equation}
\end{example}
\begin{example}
From Eq.(3.237.6) in \cite{grad}.
\begin{multline}
\int_0^{\infty } \frac{x^{-1+m} \log ^k(a x)}{(x+\gamma )^n (x+\beta )} \, dx\\
=-e^{i m \pi } (2 i \pi )^{1+k}
   \beta ^{-1+m} (-\beta +\gamma )^{-n} \Phi \left(e^{2 i m \pi },-k,\frac{\pi -i \log (a)-i \log (\beta )}{2 \pi
   }\right)\\
+\sum _{j=0}^{n-1} \sum _{p=0}^j \sum _{l=0}^p \frac{(-1)^{j+l-p} e^{i m \pi } (1-m)^{-l+p} (2 i \pi
   )^{1+k-l} \gamma ^{-1-j+m} (-\beta +\gamma )^{j-n} (1+k-l)_l \binom{p}{l} }{j!}\\
\times \Phi \left(e^{2 i m \pi },-k+l,\frac{\pi -i \log (a)-i \log (\gamma )}{2 \pi }\right) S_j^{(p)}
\end{multline}
\end{example}
\begin{example}
From  Eq. (324.8b) in \cite{grobner}
\begin{multline}
\int_0^1 x^m \left(1-x^{\alpha }\right)^n \log ^p(x) \log ^k(a x) \, dx\\
=\sum _{j=0}^{\infty } \sum
   _{l=0}^{\infty } \sum _{v=0}^n (-1)^{p+v} m^{j-l} (1+v \alpha )^{-1-j-p} (1+k-l)_l \binom{j}{l} \binom{n}{v}
   \binom{-1-p}{j} p! \log ^{k-l}(a)
\end{multline}
where $Re(a)>e^{3\pi}$ in order for the sum to converge.
\end{example}
\begin{example}
\begin{multline}
\int_0^1 e^{x^{\alpha }} x^m \log ^p(x) \log ^k(a x) \, dx\\
=e \sum _{j=0}^{\infty } \sum _{l=0}^{\infty } \sum
   _{n=0}^{\infty } \sum _{v=0}^{\infty } \frac{(-1)^n (-1)^{p+v} m^{j-l} (1+v \alpha )^{-1-j-p} (1+k-l)_l
   \binom{j}{l} \binom{n}{v} \binom{-1-p}{j} p! \log ^{k-l}(a)}{n!}
\end{multline}
\end{example}
\begin{example}
\begin{multline}
\int_0^{\infty } \frac{z^{-1+m} J_v(z) \log ^k(a z)}{1+b z^{\alpha }} \, dz\\
=-\sum _{j=0}^{\infty } \frac{i
   (-1)^j 2^{1-2 j+k-v} b^{-\frac{2 j+m+v}{\alpha }} e^{\frac{i \pi  (2 j+m+v)}{\alpha }} \pi ^{1+k}
   \left(\frac{i}{\alpha }\right)^k }{\alpha  j! \Gamma (1+j+v)}\\
\times \Phi \left(e^{\frac{2 i \pi  (2 j+m+v)}{\alpha }},-k,\frac{\alpha  \left(\frac{\pi
   }{\alpha }-i \log (a)-i \log \left(b^{-1/\alpha }\right)\right)}{2 \pi }\right)
\end{multline}
where $Re(a)<0,Re(v)<1,Im(m)\neq 0$
\end{example}
\begin{example}
\begin{multline}
\int_0^{\infty } \frac{z^{-1+m+\alpha } J_v(z) \log ^k(a z)}{\left(1+b z^{\alpha }\right)^2} \, dz\\
=-\sum
   _{j=0}^{\infty } \frac{(-1)^j 2^{-2 j+k-v} b^{-\frac{2 j+m+v+\alpha }{\alpha }} e^{\frac{i \pi  (2 j+m+v)}{\alpha
   }} \pi ^k \left(\frac{i}{\alpha }\right)^k }{\alpha ^2 \Gamma (1+j) \Gamma (1+j+v)}\\
\times \left(k \alpha  \Phi \left(e^{\frac{2 i \pi  (2 j+m+v)}{\alpha
   }},1-k,\frac{\pi -i \alpha  \left(\log (a)+\log \left(b^{-1/\alpha }\right)\right)}{2 \pi }\right)\right. \\ \left.
+2 i \pi  (2
   j+m+v) \Phi \left(e^{\frac{2 i \pi  (2 j+m+v)}{\alpha }},-k,\frac{\pi -i \alpha  \left(\log (a)+\log
   \left(b^{-1/\alpha }\right)\right)}{2 \pi }\right)\right)
\end{multline}
\end{example}
\begin{example}
\begin{multline}
\int_0^{\infty } \frac{z^{-1+m} J_v(z) \log (\log (a z))}{1+b z^{\alpha }} \, dz\\
=\sum _{j=0}^{\infty }
   \frac{2^{-2 j-v} b^{-\frac{2 j+m+v}{\alpha }} e^{\frac{i \pi  (m+v+j (2+\alpha ))}{\alpha }} \pi 
  }{\alpha  \Gamma (1+j) \Gamma (1+j+v)}\\
\times  \left(\left(-i+\cot \left(\frac{\pi  (2 j+m+v)}{\alpha }\right)\right) \log \left(\frac{2 i \pi }{\alpha }\right)\right. \\ \left.+2i \Phi'\left(e^{\frac{2 i \pi  (2 j+m+v)}{\alpha }},0,\frac{\pi -i \alpha  \left(\log
   (a)+\log \left(b^{-1/\alpha }\right)\right)}{2 \pi }\right)\right)
\end{multline}
\end{example}
\begin{example}
\begin{multline}
\int_0^{\infty } \frac{z^{-1+i} J_1(z) \log (\log (z))}{1+z} \, dz\\
=\sum _{j=0}^{\infty } \frac{2^{-1-2 j} e^{i
   ((1+i)+3 j) \pi } \pi  \left((-i+\cot (((1+i)+2 j) \pi )) \log (2 i \pi )+2 i
   \Phi'\left(e^{(-2+4 i j) \pi },0,\frac{1}{2}\right)\right)}{\Gamma (1+j) \Gamma
   (2+j)}
\end{multline}
\end{example}
\begin{example}
From Eq. (2.5.52.5(i)) in \cite{prud1}
\begin{multline}
\int_0^{\infty } e^{-i m x} \left((-i x+\log (a))^k+e^{2 i m x} (i x+\log (a))^k\right) \text{sech}(c x) \sin
   \left(\frac{c^2 x^2}{\pi }\right) \, dx\\
=-\frac{2 \left(-\frac{1}{c}\right)^k e^{-\frac{m \pi }{2 c}} \pi ^{1+k}
   \Phi \left(-e^{-\frac{m \pi }{c}},-k,-\frac{i c \left(\frac{i \pi }{2 c}-i \log (a)\right)}{\pi }\right)}{\sqrt{2}
   c}\\
+\sum _{j=0}^{\infty } \sum _{p=0}^{\infty } \frac{2 (-1)^j 2^{-4 j} \left(-\frac{1}{c}\right)^{k-p} c^{-1-4 j}
   e^{-\frac{m \pi }{2 c}} m^{4 j-p} \pi ^{1+2 j+k-p} (1+k-p)_p \binom{4 j}{p}}{(2 j)!}\\
\times  \Phi \left(-e^{-\frac{m \pi
   }{c}},-k+p,-\frac{i c \left(\frac{i \pi }{2 c}-i \log (a)\right)}{\pi }\right)
\end{multline}
where $Re(m)< Re(c), Re(a) \leq 0$.
\end{example}
\begin{example}
\begin{multline}
\int_0^{\infty } \frac{x \text{sech}(c x) \sin (m x) \sin \left(\frac{c^2 x^2}{\pi }\right)}{-a^2 \pi ^2+x^2}
   \, dx\\
=\frac{e^{-\frac{m \pi }{2 c}} \left(\Phi \left(-e^{-\frac{m \pi }{c}},1,\frac{1}{2}-i a c\right)-e^{\frac{m
   \pi }{c}} \Phi \left(-e^{\frac{m \pi }{c}},1,\frac{1}{2}-i a c\right)\right)}{2 \sqrt{2}}\\
-\sum _{j=0}^{\infty }
   \sum _{p=0}^{\infty } \frac{(-1)^{j-2 p} 2^{-1-4 j} c^{-4 j+p} e^{-\frac{m \pi }{2 c}} m^{4 j-p} \pi ^{2 j-p}
   \binom{4 j}{p} }{(2 j)!}\\
\times \left((-1)^p \Phi \left(-e^{-\frac{m \pi }{c}},1+p,\frac{1}{2}-i a c\right)-(-1)^{4 j} e^{\frac{m
   \pi }{c}} \Phi \left(-e^{\frac{m \pi }{c}},1+p,\frac{1}{2}-i a c\right)\right) (-p)_p
\end{multline}
where $Im(m)\neq 0$
\end{example}
\begin{example}
From Eq. (4.294.1) in \cite{grad}
\begin{multline}
\int_0^1 x^{-2-m} \left(-k x \log ^{-1+k}\left(\frac{a}{x}\right)-m x \log ^k\left(\frac{a}{x}\right)\right. \\ \left.
+x^{2 m} \log ^{-1+k}(a x) (k+(-1+m) \log (a x))\right) \log (1+x) \, dx\\
=e^{i m \pi } (2
   i \pi )^{1+k} \Phi \left(e^{2 i m \pi },-k,\frac{\pi -i \log (a)}{2 \pi }\right)+2 \log (2) \log ^k(a)
\end{multline}
\end{example}
\begin{example}
\begin{multline}
\int_0^1 \frac{\log (1+x) \left(-\frac{2}{\log \left(\frac{a}{x}\right)}+\frac{2}{\log (a x)}-\log \left(\log \left(\frac{a}{x}\right)\right)-\log (\log (a x))\right)}{2 x^{3/2}} \,
   dx\\
=\frac{1}{2} \pi  (-i \pi -2 \log (4 \pi )+4 \log (-3 \pi -i \log (a))-4 \log (-\pi -i \log (a)))\\
+\log (4) \log (\log (a))-2 \pi  \text{log$\Gamma $}\left(-\frac{\pi +i \log (a)}{4 \pi
   }\right)+2 \pi  \text{log$\Gamma $}\left(-\frac{3}{4}-\frac{i \log (a)}{4 \pi }\right)
\end{multline}
\end{example}
\begin{example}
\begin{multline}
\int_0^1 \frac{\log (1+x) \left(\frac{4 \log (x)}{\pi ^2+\log ^2(x)}-\log \left(-\pi ^2-\log ^2(x)\right)\right)}{2 x^{3/2}} \, dx\\
=-\frac{1}{2} \left(i \pi ^2\right)-(2-i) \pi  \log (2)+\log
   (4) \log (\pi )
\end{multline}
\end{example}
\begin{example}
\begin{multline}
\int_0^1 \left(\frac{2 \log (x) \log (1+x)}{x^{3/2} \left(\pi ^2+\log ^2(x)\right)}-\frac{\log (1+x) \log \left(-\pi ^2-\log ^2(x)\right)}{2 x^{3/2}}\right) \, dx\\
=-2 \pi  \log (2)+i
   \left(-\frac{\pi ^2}{2}+\pi  \log (2)\right)+\log (4) \log (\pi )
\end{multline}
\end{example}
\begin{example}
\begin{multline}
\int_0^1 \frac{\log (1+x) \left(\frac{4 \log (x)}{\pi ^2+\log ^2(x)}-\log \left(\pi ^2+\log ^2(x)\right)\right)}{2 x^{3/2}} \, dx=-2 \pi  \log (2)+\log (4) \log (\pi )
\end{multline}
\end{example}
\begin{example}
\begin{multline}
\int_0^1 \frac{\log (1+x) \left(-\frac{2}{\log \left(\frac{a}{x}\right)}+\frac{2}{\log (a x)}-\log \left(\log \left(\frac{a}{x}\right)\right)-\log (\log (a x))\right)}{2 x^{3/2}} \,
   dx\\
=\frac{1}{2} \pi  (-i \pi -2 \log (4 \pi )+4 \log (-3 \pi -i \log (a))-4 \log (-\pi -i \log (a)))+\log (4) \log (\log (a))\\
-2 \pi  \text{log$\Gamma $}\left(-\frac{\pi +i \log (a)}{4 \pi
   }\right)+2 \pi  \text{log$\Gamma $}\left(-\frac{3}{4}-\frac{i \log (a)}{4 \pi }\right)
\end{multline}
\end{example}
\begin{example}
\begin{multline}
\int_0^1 \frac{\log (1+x) \left(-\frac{2}{\log \left(\frac{i}{x}\right)}+\frac{2}{\log (i x)}-\log \left(\log \left(\frac{i}{x}\right)\right)-\log (\log (i x))\right)}{2 x^{3/2}} \,
   dx\\
=\frac{1}{2} \left(-i \pi ^2-\log ^2(4)+\log (16) \log (\pi )+\pi  \left(i \log (4)+\log \left(\frac{625}{16}\right)\right.\right. \\ \left.\left.
-2 \log (\pi )+4 \text{log$\Gamma $}\left(-\frac{5}{8}\right)-4
   \text{log$\Gamma $}\left(-\frac{1}{8}\right)\right)\right)
\end{multline}
\end{example}
\begin{example}
\begin{multline}
\int_0^1 \frac{x^{-2-u-i v} }{\log \left(-\frac{1}{x}\right) \log (-x)}\\
\left(x \left(-u-i v+(u-i v) x^{2 i v}\right) \log (-x)-\frac{x \left(-1+x^{2 i v}\right) \log (-x)}{\log \left(-\frac{1}{x}\right)}\right. \\ \left.
+\frac{x^{2 u} \log \left(-\frac{1}{x}\right) \left(1-x^{2 i v}+\left(1-u+i v+(-1+u+i v) x^{2 i v}\right) \log (-x)\right)}{\log (-x)}\right) \log (1+x) \, dx\\
=e^{-\pi  (i
   u+v)} \left(-e^{2 \pi  v} \log \left(1-e^{2 i \pi  (u+i v)}\right)+\log \left(1-e^{2 \pi  (i u+v)}\right)\right)
\end{multline}
\end{example}
\begin{example}
\begin{multline}
\int_0^1 \frac{\log (1+x) \left(\frac{2}{\log \left(-\frac{1}{x}\right)}-\frac{2}{\log (-x)}+\log \left(\log \left(-\frac{1}{x}\right)\right)+\log (\log (-x))\right)}{2 x^{3/2}} \,
   dx\\
=\frac{i \pi ^2}{2}+(2-i) \pi  \log (2)-\log (4) \log (\pi )
\end{multline}
\end{example}
\begin{example}
\begin{multline}
\int_0^1 \frac{\log (1+x) }{x^{5/2}}\\
\left(-x \log \left(\log \left(-\frac{1}{x}\right)\right)-\frac{1}{2} x \log \left(-\frac{1}{x}\right) \log \left(\log \left(-\frac{1}{x}\right)\right)\right. \\ \left.
+x \left(1-\frac{\log (-x)}{2}\right) \log (\log (-x))\right) \, dx\\
=4 i \pi ^2 \log \left(\frac{A^3}{\sqrt[3]{2} \sqrt[4]{e}}\right)+2 i \pi  \log (2) \log (i \pi )-i \pi ^2 \log (2 i \pi
   )
\end{multline}
\end{example}
\begin{example}
\begin{multline}
\int_0^1 \frac{\log (1+x) }{2 x^{3/2}}\\
\left(\log \left(-\frac{1}{x}\right) \left(2+\left(4+\log \left(-\frac{1}{x}\right)\right) \log \left(\log \left(-\frac{1}{x}\right)\right)\right)\right. \\ \left.
+\log (-x)(-2+(-4+\log (-x)) \log (\log (-x)))\right) \, dx\\
=2 \pi  (\pi  \log (2) \log (i \pi )+7 \zeta (3))
\end{multline}
\end{example}
\begin{example}
From Eq. (3.27.1) and (3.10.1.2) in \cite{brychkov}, 
\begin{multline}
\int _0^{\infty }\int _0^{\infty }\int _0^{\infty }e^{-b z^2} x^{-1+m} y^{\frac{1}{2} (-m+v)} z^{1-m-v} J_v(x \beta ) E_v(y \alpha ) \log ^{-1+k}\left(\frac{a x}{\sqrt{y}
   z}\right)\\
 \left(k+(m-3 v) \log \left(\frac{a x}{\sqrt{y} z}\right)\right)dzdydx\\
=2^m b^{\frac{1}{2} (-2+m+v)} e^{\frac{1}{2} i \pi  (m+v)} (i \pi )^{1+k} \alpha ^{\frac{1}{2}
   (-2+m-v)} \beta ^{-m}\\
 \Phi \left(e^{i \pi  (m+v)},-k,\frac{\pi -2 i \log (2 a)-i \log (b)-i \log (\alpha )+2 i \log (\beta )}{2 \pi }\right)
\end{multline}
\end{example}
\begin{example}
\begin{equation}
\int _0^{\infty }\int _0^{\infty }\int _0^{\infty }\frac{e^{-z^2} y^{5/6} J_{\frac{4}{3}}(x) E_{\frac{4}{3}}(y) \left(3+13 \log \left(\frac{i x}{2 \sqrt{y} z}\right)\right)}{3
   x^{4/3} \log ^2\left(\frac{i x}{2 \sqrt{y} z}\right)}dzdydx=-\frac{i \log (2)}{\sqrt[3]{2}}
\end{equation}
\end{example}
\begin{example}
\begin{equation}
\int _0^{\infty }\int _0^{\infty }\int _0^{\infty }\frac{e^{-z^2} y^{5/6} J_{\frac{4}{3}}(x) E_{\frac{4}{3}}(y) \left(9+13 \log \left(\frac{i x}{2 \sqrt{y} z}\right)\right)}{3
   x^{4/3} \log ^4\left(\frac{i x}{2 \sqrt{y} z}\right)}dzdydx=\frac{3 i \zeta (3)}{4 \sqrt[3]{2} \pi ^2}
\end{equation}
\end{example}
\begin{example}
\begin{multline}
\int _0^{\infty }\int _0^{\infty }\int _0^{\infty }\int _0^{\infty }\int _0^{\infty }\log ^k\left(\frac{a u^2 x
   y^{\nu }}{t^2 z \left(p+q y^{\nu }\right)}\right) e^{-b \left(t^2+u^2+z\right)-x \alpha } t^{-2 m} u^{2 (m-v)}
   x^{-1+m-v} y^{-1+\left(\frac{1}{2}+m\right) \nu }\\
 \left(p+q y^{\nu }\right)^{-m+2 v} z^{-m} K_v(x \alpha
   )dudtdzdydx\\
=-\frac{i 2^{2 k-m+v} b^{-2+m+v} e^{\frac{1}{2} i (k+4 m) \pi } p^{-m+2 v} \pi ^{\frac{5}{2}+k}
   \left(\frac{p}{q}\right)^{\frac{1}{2}+m} \alpha ^{-m+v} \Gamma \left(-\frac{1}{2}-2 v\right) }{\nu }\\
\Phi \left(e^{4 i m \pi
   },-k,\frac{2 \pi -i \log (a)-i \log (b)+i \log (2 p)-i \log \left(\frac{p}{q}\right)+i \log (\alpha )}{4 \pi
   }\right)
\end{multline}
where $Re(m)<|Re(v)|$
\end{example}
\begin{example}
\begin{multline}
\sum _{j=0}^{2 n} (-1)^j e^{i \left(m+\frac{j \pi }{1+2 n}\right)} \Phi \left(-e^{2 i \left(m+\frac{j \pi
   }{1+2 n}\right)},-k,\frac{1}{2} (1-i \log (a))\right)\\
=(-1)^n i^{-k} e^{i m (1+2 n)} (i (1+2 n))^k (1+2 n) \Phi
   \left(-e^{2 i m (1+2 n)},-k,\frac{1+2 n-i \log (a)}{2+4 n}\right)
\end{multline}
\end{example}
\begin{example}
\begin{multline}
\Phi (z,s,a)=3^{1-s} z \Phi \left(z^3,s,\frac{a+1}{3}\right)-(-1)^{2/3} \Phi \left(-\sqrt[3]{-1}
   z,s,a\right)+\sqrt[3]{-1} \Phi \left((-1)^{2/3} z,s,a\right)
\end{multline}
\end{example}
\begin{example}
\begin{multline}
\Phi (-z,s,a)\\
=5^{1-s} z^2 \Phi
   \left(-z^5,s,\frac{a+2}{5}\right)+(-1)^{3/5} \Phi \left(\sqrt[5]{-1}
   z,s,a\right)+\sqrt[5]{-1} \Phi \left(-(-1)^{2/5} z,s,a\right)\\
   -(-1)^{2/5}
   \left((-1)^{2/5} \Phi \left((-1)^{3/5} z,s,a\right)+\Phi \left(-(-1)^{4/5}
   z,s,a\right)\right)
\end{multline}
\end{example}
\begin{example}
\begin{multline}
\sqrt[3]{-1} \Phi'\left((-1)^{2/3},0,a\right)-(-1)^{2/3}
   \Phi'\left(-\sqrt[3]{-1},0,a\right)
   =\log \left(\frac{2 \pi  3^{\frac{1}{2}-a}
   \Gamma (a)}{\Gamma \left(\frac{a+1}{3}\right)^3}\right)
\end{multline}
\end{example}
\begin{example}
\begin{multline}
\sqrt[3]{-1}
   \Phi'\left(-(-1)^{2/3},0,a\right)-(-1)^{2/3}
   \Phi'\left(\sqrt[3]{-1},0,a\right)\\
   =3 \log
   \left(\frac{\sqrt{3} \sqrt[3]{\frac{\Gamma
   \left(\frac{a-2}{2}\right)}{\Gamma \left(\frac{a-1}{2}\right)}} \Gamma
   \left(\frac{a+1}{6}\right)}{\left((a-2) \sqrt{a-1}\right)^{2/3} \Gamma
   \left(\frac{a-2}{6}\right)}\right)+\log (2)
\end{multline}
\end{example}
\begin{example}
\begin{multline}
\prod _{j=0}^{2 n} \left(\frac{\cos \left(\frac{1}{2} \left(m+\frac{2 j \pi }{1+2 n}+r\right)\right)+\sin
   \left(\frac{m-r}{2}\right)}{\cos \left(\frac{m+2 m n+2 j \pi +r+2 n r}{2+4 n}\right)-\sin
   \left(\frac{m-r}{2}\right)}\right)^{(-1)^j}\\
   =\left(-\frac{\sin \left(\frac{1}{4} (m (2+4 n)+\pi )\right)}{\tan
   \left(\frac{1}{4} (\pi +2 r+4 n r)\right) \sin \left(m \left(\frac{1}{2}+n\right)-\frac{\pi
   }{4}\right)}\right)^{(-1)^n}
\end{multline}
\end{example}
\begin{example}
\begin{multline}
\prod _{j=0}^{2 n} \exp \left(-2 (-1)^j \sec \left(\frac{j \pi }{1+2 n}+x\right) \sec \left(\frac{j \pi }{1+2
   n}+\frac{x}{b}\right)\right. \\ \left.
    \sin \left(\frac{(-1+b) x}{2 b}\right) \sin \left(\frac{1}{2} \left(\frac{2 j \pi }{1+2
   n}+x+\frac{x}{b}\right)\right)\right)\\
    \left(\frac{\left(-1+\sin \left(\frac{j \pi }{1+2 n}+x\right)\right)
   \left(1+\sin \left(\frac{j \pi }{1+2 n}+\frac{x}{b}\right)\right)}{\left(1+\sin \left(\frac{j \pi }{1+2
   n}+x\right)\right) \left(-1+\sin \left(\frac{j \pi }{1+2
   n}+\frac{x}{b}\right)\right)}\right)^{\frac{(-1)^j}{2}}\\
   =\exp \left(\frac{2 (-1)^n (1+2 n) \left(\cos (x+2 n
   x)-\cos \left(\frac{x+2 n x}{b}\right)\right)}{\cos \left(\frac{(-1+b) (x+2 n x)}{b}\right)+\cos
   \left(\frac{(1+b) (x+2 n x)}{b}\right)}\right)\\
    \left(\frac{(-1+\sin (x+2 n x)) \left(1+\sin \left(\frac{x+2 n
   x}{b}\right)\right)}{(1+\sin (x+2 n x)) \left(-1+\sin \left(\frac{x+2 n
   x}{b}\right)\right)}\right)^{\frac{(-1)^n}{2}}
\end{multline}
\end{example}
\begin{example}
\begin{equation}
\prod _{j=0}^{2 n} \exp \left(\frac{(-1)^{j-n} e^{\frac{i j \pi }{1+2 n}}
   \Phi'\left(-e^{\frac{2 i j \pi }{1+2 n}},0,\frac{1+a}{2}\right)}{1+2
   n}\right)=\frac{\Gamma \left(\frac{1+a+2 n}{4+8 n}\right)}{\sqrt{2+4 n} \Gamma \left(\frac{3+a+6 n}{4+8
   n}\right)}
\end{equation}
\end{example}
\begin{example}
\begin{multline}
\sum _{n=1}^{\infty } -a^{-i n} \left(n^2\right)^{-2-k} \cos (n t) \left((i n)^{2+k} \Gamma (2+k,-i n \log (a))+a^{2 i n} (-i n)^{2+k} \Gamma
   (2+k,i n \log (a))\right)\\
   =-2^{k+1} \pi ^{2+k} \left(\zeta \left(-1-k,\frac{2 \pi -t+\log (a)}{2 \pi }\right)+\zeta \left(-1-k,\frac{t+\log (a)}{2
   \pi }\right)\right)-\frac{\log ^{2+k}(a)}{k+2}
\end{multline}
\end{example}
\begin{example}
\begin{equation}
\sum _{n=1}^{\infty } e^{-i n \pi } \cos (n \pi ) \left(\Gamma (0,-i n \pi )+e^{2 i n \pi } \Gamma (0,i n \pi )\right)=\gamma -\log (2)
\end{equation}
\end{example}
\begin{example}
\begin{multline}
-\sum _{n=1}^{\infty } \frac{\cos (n \pi )}{6 \pi } \left(12 (1+\log (\pi ))+n \pi  \left(\left(6 \gamma ^2+\pi ^2-12 \text{Ci}(n \pi ) (1+\log (\pi
   ))\right) \sin (n \pi )\right.\right. \\ \left.\left.
+3 i e^{-i n \pi } \left(2 i n \pi  \, _3F_3(1,1,1;2,2,2;i n \pi )+(-2+\log (-i n)) \log (-i n)+\log ^2(\pi )\right.\right.\right. \\ \left.\left.\left.
+2 \log (n)(1+\log (\pi ))+2 \gamma  \log (-i n \pi )-2 i (1+\log (\pi )) \text{Si}(n \pi )\right)\right.\right. \\ \left.\left.
-3 i e^{i n \pi } \left(-2 i n \pi  \, _3F_3(1,1,1;2,2,2;-i n \pi )+(-2+\log (i n)) \log (i n)+\log ^2(\pi )\right.\right.\right. \\ \left.\left.\left.
+2 \log (n) (1+\log (\pi ))+2 \gamma  \log (i n \pi )+2 i (1+\log (\pi )) \text{Si}(n \pi
   )\right)\right)\right)\\
=\frac{\gamma  \pi }{12}-\pi  \log (A)+\frac{1+\log (\pi )}{\pi }
\end{multline}
\end{example}
\begin{example}
\begin{equation}
\sum _{n=1}^{\infty } e^{-i n \pi } \cos \left(\frac{n \pi }{2}\right) \left(E_2(-i n \pi )+e^{2 i n \pi } E_2(i n \pi )\right)=-5+\frac{\pi
   ^2}{2}
\end{equation}
\end{example}
\begin{example}
\begin{multline}
\sum _{n=1}^{\infty } \frac{e^{-i n \pi } \cos (n t) \left(-\Gamma (0,-i n \pi )+e^{2 i n \pi } \Gamma (0,i n \pi )\right)}{n}\\
=-i \pi 
   \left(1+\log \left(\frac{\pi -t}{4 \pi ^3}\right)+\log (-\pi +t)+\text{log$\Gamma $}\left(\frac{\pi -t}{2 \pi }\right)+\text{log$\Gamma
   $}\left(\frac{1}{2} \left(-1+\frac{t}{\pi }\right)\right)\right)
\end{multline}
\end{example}
\begin{example}
\begin{multline}
\sum _{n=1}^{\infty } \frac{e^{-i n \pi } \cos \left(\frac{n \pi }{2}\right) \left(E_1(-i n \pi )-e^{2 i n \pi } E_1(i n \pi )\right)}{n}\\
=i \pi 
   \left(1+i \pi -\log (16 \pi )+\text{log$\Gamma $}\left(-\frac{1}{4}\right)+\text{log$\Gamma $}\left(\frac{1}{4}\right)\right)
\end{multline}
\end{example}
\begin{example}
\begin{multline}
\sum _{n=1}^{\infty } \frac{e^{-n} \cos (n t) \left(-i \Gamma \left(\frac{5}{2},-n\right)+e^{2 n} \Gamma
   \left(\frac{5}{2},n\right)\right)}{n^{5/2}}\\
   =\frac{2}{5}-(2-2 i) \pi ^{5/2} \left(\zeta \left(-\frac{3}{2},\frac{i+2 \pi -t}{2 \pi }\right)+\zeta
   \left(-\frac{3}{2},\frac{i+t}{2 \pi }\right)\right)
\end{multline}
\end{example}
\begin{example}
\begin{equation}
\frac{2 \left(1-x^2\right)^{b m} \left(1-x^2\right)}{1-x^2+2 b x^2}
=\sum _{n=0}^{\infty }
   \frac{\left(1+(-1)^n\right) \left(\frac{x}{\left(1-x^2\right)^b}\right)^n \Gamma \left(\frac{1}{2} (n-2 b
   (m+n))\right)}{\Gamma \left(\frac{2+n}{2}\right) \Gamma (-b (m+n))}
\end{equation}
where $Re(b)<0,|Re(x)|<1,m\in\mathbb{C}$
\end{example}
\begin{example}
\begin{multline}
\int_u^v \frac{2 (1+x)^a \left(1-x^2\right)^{1+b m}}{1-x^2+2 b x^2} \, dx\\
=-\sum _{n=0}^{\infty } \frac{\left(1+(-1)^n\right)}{(1+n)
   \Gamma \left(\frac{2+n}{2}\right) \Gamma (-b (m+n))}
   \left(u \left(1-u^2\right)^{b n} \left(u \left(1-u^2\right)^{-b}\right)^n \right. \\ \left.
F_1(1+n;b n,-a+b n;2+n;u,-u)\right. \\ \left.
-v \left(1-v^2\right)^{b n}
   \left(v \left(1-v^2\right)^{-b}\right)^n F_1(1+n;b n,-a+b n;2+n;v,-v)\right) \Gamma \left(\frac{1}{2} (n-2 b (m+n))\right)
\end{multline}
\end{example}
\begin{example}
\begin{equation}
\frac{e^{\frac{a x}{-1+x}} (-1+x)^2 x^m}{1+x (-2+b+x)}=\sum _{n=0}^{\infty } \left(e^{-\frac{b x}{-1+x}}
   x\right)^{m+n} L_n^{-1}(a+b (m+n))
\end{equation}
where $|Re(b)|<1,|Re(x)|<1$
\end{example}
\begin{example}
\begin{equation}
\frac{e^{\frac{a x}{1-x}} (1-x)^{1-v}}{1-(2+b) x+x^2}=\sum _{n=0}^{\infty } \left(e^{\frac{b x}{-1+x}} x\right)^n
   L_n^v(-a-b n)
\end{equation}
where $|Re(x)|<1$
\end{example}
\begin{example}
\begin{equation}
e^{-\frac{b m}{-1+x}} x^{-m} \int_0^x e^{\frac{b m-a t}{-1+t}} (1-t)^{-1-v} t^{-1+m} \, dt=\sum _{n=0}^{\infty }
   \frac{\left(e^{\frac{b x}{-1+x}} x\right)^n L_n^v(-a-b n)}{n+m}
\end{equation}
\end{example}
\begin{example}
\begin{multline}
\int_0^x \exp \left(\frac{a x \beta +t (\alpha -x \alpha -a \beta )}{(-1+t) (-1+x)}\right) (1-t)^{-1-v} t^{-1+a}
   x^{-a} \, dt\\
=\sum _{n=0}^{\infty } \frac{\left(e^{\frac{x \beta }{-1+x}} x\right)^n \text{csch}((a+n) \pi ) \Gamma
   (1+n+v) }{(a+n) \Gamma (1+n) \Gamma (1+v)}\\
\times \, _1F_1(-n;1+v;-\alpha -n \beta ) ((a+n) \pi +\sinh ((a+n) \pi ))
\end{multline}
where $Re(a)>2\pi$
\end{example}
\begin{example}
\begin{equation}
\frac{e^{\frac{a x}{1-x}} (1-x)^{2-c} x^m}{1-(2+b) x+x^2}=\sum _{n=0}^{\infty } \left(e^{\frac{b x}{-1+x}}
   x\right)^{m+n} L_n^{-1+c}(-a-b (m+n))
\end{equation}
\end{example}
\begin{example}
\begin{multline}
\frac{e^{\frac{a x}{1-x}} (1-x)^{2-c} (1+x)^p}{1-(2+b) x+x^2}=\sum _{m=0}^{\infty } \sum _{n=0}^{\infty }
   \left(e^{\frac{b x}{-1+x}} x\right)^{m+n} \binom{p}{m} L_n^{-1+c}(-a-b (m+n))
\end{multline}
\end{example}
\begin{example}
\begin{multline}
\frac{(1-x)^{1+b} x^m \log (1-x z)}{1+(-1+a) x}\\
=\sum _{n=0}^{\infty } \frac{(-1)^{2 m+n} n \left((1-x)^{-a}
   x\right)^{m+n} z \Gamma (1+b+a (m+n)) }{(m+n)! \Gamma (2+b-n+a
   (m+n))}\\
\, _3F_2(1,1,1-n;2,2+b+a m-n+a n;z) (1+n)_m
\end{multline}
\end{example}
\begin{example}
\begin{multline}
\frac{x^m (1+x)^{1+b} \log (1-x z)}{1+x-a x}\\
=-\sum _{n=0}^{\infty } \frac{n \left(x (1+x)^{-a}\right)^{m+n} z
   \Gamma (1+b+a (m+n)) }{(m+n)! \Gamma (2+b-n+a (m+n))}\\\, _3F_2(1,1,1-n;2,2+b+a m-n+a n;-z) (1+n)_m
\end{multline}
\end{example}
\begin{example}
\begin{multline}
-\frac{e^{a x} x^m}{-1+b x+x \log \left(1-\frac{z}{c}\right)}\\
=\sum _{n=0}^{\infty } \frac{\left(e^{-b x} x
   \left(1-\frac{z}{c}\right)^{-x}\right)^{m+n} \left(a+b (m+n)+(m+n) \log \left(1-\frac{z}{c}\right)\right)^n}{\Gamma
   (1+n)}
\end{multline}
\end{example}
\begin{example}
\begin{multline}
\frac{(1-x)^{1+b} \tanh ^{-1}(x z)}{1+(-1+a) x}\\
=-\sum _{n=1}^{\infty } \frac{(-1)^n \left((1-x)^{-a} x\right)^n z
   \Gamma (1+b+a n) }{2 \Gamma (n) \Gamma
   (2+b-n+a n)}\\
\left(\, _3F_2(1,1,1-n;2,2+b-n+a n;-z)+\, _3F_2(1,1,1-n;2,2+b-n+a n;z)\right)
\end{multline}
\end{example}
\begin{example}
\begin{multline}
\frac{(1-x)^b \tanh ^{-1}(u-x z)}{1+\frac{a x}{1-x}}\\
=-\sum _{n=0}^{\infty } \frac{(-1)^n \left((1-x)^{-a}
   x\right)^n \Gamma (1+b+a n) }{2 \left(-1+u^2\right) n! \Gamma (2+b-n+a n)}\left(2 (1+b-n+a n) \tanh ^{-1}(u)\right. \\ \left.
-2 (1+b-n+a n) u^2 \tanh ^{-1}(u)+n z \,_3F_2\left(1,1,1-n;2,2+b-n+a n;\frac{z}{-1+u}\right)\right. \\ \left.
+n u z \, _3F_2\left(1,1,1-n;2,2+b-n+a n;\frac{z}{-1+u}\right)\right. \\ \left.
+nz \, _3F_2\left(1,1,1-n;2,2+b-n+a n;\frac{z}{1+u}\right)-n u z \, _3F_2\left(1,1,1-n;2,2+b-n+an;\frac{z}{1+u}\right)\right)
\end{multline}
\end{example}
\begin{example}
From Eq. (3.311.1(11)) in \cite{grad}
\begin{multline}
\int_{-\infty }^{\infty } \frac{e^{-m x} \log ^k\left(a e^{-x}\right)}{1+e^{-q x}} \, dx=-\frac{i e^{\frac{i m
   \pi }{q}} (2 \pi )^{1+k} \left(\frac{i}{q}\right)^k \Phi \left(e^{\frac{2 i m \pi }{q}},-k,\frac{q \left(\frac{\pi
   }{q}-i \log (a)\right)}{2 \pi }\right)}{q}
\end{multline}
\end{example}
\begin{example}
\begin{equation}
\int_{-\infty }^{\infty } \frac{e^{-u x} \sinh (v x)}{\left(1+e^{-q x}\right) x} \, dx=\frac{1}{2} \log
   \left(\cot \left(\frac{\pi  (u-v)}{2 q}\right) \tan \left(\frac{\pi  (u+v)}{2 q}\right)\right)
\end{equation}
\end{example}
\begin{example}
\begin{equation}
\int_{-\infty }^{\infty } \frac{e^{x/2} \log (\pi  i-x)}{1+e^x} \, dx=\frac{i \pi ^2}{2}+\pi  \log (4)
\end{equation}
\end{example}
\begin{example}
\begin{multline}
\int_{-\infty }^{\infty } (a-x)^k \text{sech}\left(\frac{q x}{2}\right) \, dx=\frac{(4 \pi )^{1+k}
   \left(\frac{i}{q}\right)^k \left(\zeta \left(-k,\frac{\pi -i q a}{4 \pi }\right)-\zeta \left(-k,\frac{3}{4}-\frac{i q
   a}{4 \pi }\right)\right)}{q}
\end{multline}
\end{example}
\begin{example}
\begin{multline}
\int_{-\infty }^{\infty } \log (a-x) \text{sech}\left(\frac{q x}{2}\right) \, dx\\
=\frac{2 \pi }{q} \left(\log
   \left(\frac{4 i \pi }{q}\right)-2 \log (-3 \pi -i a q)+2 \log (-\pi -i a q)\right. \\ \left.
+2 \text{log$\Gamma $}\left(-\frac{\pi +i
   a q}{4 \pi }\right)-2 \text{log$\Gamma $}\left(-\frac{3}{4}-\frac{i a q}{4 \pi }\right)\right)
\end{multline}
\end{example}
\begin{example}
\begin{equation}
\int_{-\infty }^{\infty } \log \left(x^2-a^2 \pi ^2\right) \text{sech}\left(\frac{q x}{2}\right) \, dx=\frac{4
   \pi  \log \left(\frac{4 i \pi  \Gamma \left(\frac{3}{4}-\frac{i a q}{4}\right) \Gamma \left(\frac{3}{4}+\frac{i a
   q}{4}\right)}{q \Gamma \left(\frac{1}{4}-\frac{i a q}{4}\right) \Gamma \left(\frac{1}{4}+\frac{i a
   q}{4}\right)}\right)}{q}
\end{equation}
\end{example}
\begin{example}
\begin{equation}
\int_{-\infty }^{\infty } \log \left(\frac{x-a}{a+x}\right)
   \text{sech}\left(\frac{q x}{2}\right) \, dx=\frac{4 \pi  \log \left(\tan
   \left(\frac{1}{4} (\pi +i a q)\right)\right)}{q}
\end{equation}
\end{example}
\begin{example}
\begin{equation}
\int_{-\infty }^{\infty } x \log (x) \text{sech}\left(\frac{q x}{2}\right) \, dx=\frac{8 i C \pi }{q^2}
\end{equation}
\end{example}
\begin{example}
\begin{equation}
\int_{-\infty }^{\infty } \log \left(\log \left(i e^{-x}\right)\right) \text{sech}(x) \, dx=\frac{1}{2} \pi  (i
   \pi +\log (4))
\end{equation}
\end{example}
\begin{example}
\begin{equation}
\int_{-\infty }^{\infty } \frac{\log \left(\log \left(i e^{-x}\right)\right) \text{sech}(x)}{\log \left(i
   e^{-x}\right)} \, dx=\log (2) \left(2 i \gamma +\pi -i \log \left(2 \pi ^2\right)\right)
\end{equation}
\end{example}
\begin{example}
\begin{equation}
\int_{-\infty }^{\infty } \frac{\log \left(\log \left(i e^{-x}\right)\right) \text{sech}(x)}{\log ^2\left(i
   e^{-x}\right)} \, dx=\frac{1}{12} \pi  \log \left(\frac{16 \exp (2 \gamma )}{\exp (\pi  i) A^{24}}\right)
\end{equation}
\end{example}
\begin{example}
\begin{multline}
\int_{-\infty }^{\infty } \frac{\log \left(\log \left(i e^{-x}\right)\right) \text{sech}(x)}{\log ^3\left(i
   e^{-x}\right)} \, dx=-\frac{3 \pi  \zeta (3)+2 i \left(\log \left(\frac{2}{\pi ^3}\right) \zeta (3)+3 \zeta
   '(3)\right)}{4 \pi ^2}
\end{multline}
\end{example}
\begin{example}
\begin{multline}
\int_{-\infty }^{\infty } \frac{\log \left(\log \left(i e^{-x}\right)\right) \text{sech}(x)}{\log ^5\left(i
   e^{-x}\right)} \, dx=\frac{15 \pi  \zeta (5)+2 i \left((\log (2)-15 \log (\pi )) \zeta (5)+15 \zeta '(5)\right)}{16
   \pi ^4}
\end{multline}
\end{example}
\begin{example}
From Eq. (1.9.58) in \cite{erd_t1}
\begin{multline}
\int_0^{\infty } e^{-i m x} \left((-i x+\log (a))^k+e^{2 i m x} (i x+\log (a))^k\right)\\
   \text{sech}\left(\sqrt{\frac{\pi }{2}} x\right) \left(\cos \left(\frac{x^2}{2}\right)+\sin
   \left(\frac{x^2}{2}\right)\right) \, dx\\
=\sum _{j=0}^{\infty } \sum _{n=0}^{\infty } \frac{2 (-1)^n
   2^{\frac{1}{2}-j+\frac{j-k}{2}+k-2 n} e^{m \sqrt{\frac{\pi }{2}}} m^{-j+4 n} \pi
   ^{\frac{1}{2}+\frac{1}{2} (-j+k)} (1-j+k)_j \binom{4 n}{j} }{(2 n)!}\\
\times \Phi \left(-e^{m \sqrt{2 \pi
   }},j-k,\frac{\sqrt{\frac{\pi }{2}}+\log (a)}{\sqrt{2 \pi }}\right)\\
+\sum _{j=0}^{\infty } \sum
   _{n=0}^{\infty } \frac{2 (-1)^n 2^{-\frac{1}{2}-j+\frac{j-k}{2}+k-2 n} e^{m \sqrt{\frac{\pi }{2}}}
   m^{2-j+4 n} \pi ^{\frac{1}{2}+\frac{1}{2} (-j+k)} (1-j+k)_j \binom{4 n}{j} }{(1+2 n)!}\\
\times \Phi \left(-e^{m \sqrt{2 \pi
   }},j-k,\frac{\sqrt{\frac{\pi }{2}}+\log (a)}{\sqrt{2 \pi }}\right)
\end{multline}
where $Re(a)<-e^{e^{3\pi}},|Re(m)|<1$
\end{example}
\begin{example}
From Eq. (6.1.75) in \cite{hansen}
\begin{multline}
\sum _{j=0}^{\infty } \sum _{p=0}^{\infty } \frac{(-1)^p a^{-2 i p x} m^{-j} (1-j+k)_j \binom{2}{j} }{(1+2 p)^2}\\
\times \left(a^{2 i
   (1+2 p) x} (-m+i (x+2 p x))^{-1+j-k} \Gamma (1-j+k,-((m-i (x+2 p x)) \log (a)))\right. \\ \left.
-(-m-i (x+2 p x))^{-1+j-k} \Gamma(1-j+k,-((m+i (x+2 p x)) \log (a)))\right)\\
=-\frac{i a^{m+i x} \pi  x \left(-2 e^{\frac{m \pi }{2 x}} \pi
   ^k x^{-k} \Phi \left(-e^{\frac{m \pi }{x}},-k,\frac{1}{2}+\frac{x \log (a)}{\pi }\right)+\log ^k(a)\right)}{2
   m^2}
\end{multline}
\end{example}
\begin{example}
\begin{multline}
\sum _{j=0}^{\infty } \sum _{p=0}^{\infty } \frac{(-1)^p m^{-j} \left(-(-m-i (x+2 p x))^{-1+j-k}+(-m+i (x+2 p
   x))^{-1+j-k}\right) \binom{2}{j}}{(1+2 p)^2}\\
=\frac{i e^{\frac{m \pi }{2 x}} \pi ^{1+k} x^{1-k} \Phi \left(-e^{\frac{m
   \pi }{x}},-k,\frac{1}{2}\right)}{\Gamma (1+k) m^2}
\end{multline}
\end{example}
\begin{example}
\begin{multline}
\sum _{j=0}^{\infty } \sum _{p=0}^{\infty } \frac{(-1)^p 2^{-j} (i x)^{-j} \binom{2}{j}}{(1+2 p)^2}\\
\times  \left(-(i (-1+2 p)
   x)^{-1+j} \log (i (-1+2 p) x)+(-i (3+2 p) x)^{-1+j} \log (-i (3+2 p) x)\right)\\
=\frac{i \pi  \log
   \left(\frac{2 e^{\gamma } \pi  \Gamma \left(-\frac{1}{4}\right)^2}{9 x \Gamma \left(-\frac{3}{4}\right)^2}\right)}{8
   x}
\end{multline}
\end{example}
\begin{example}
\begin{multline}
\sum _{j=0}^{\infty } \sum _{p=0}^{\infty } \frac{\left((-1)^p \binom{2}{j}\right) \left(\frac{\log ((3+2 p)
   x)}{(3+2 p)^{1-j}}-\frac{\log (x (1-2 p))}{(1-2 p)^{1-j}}\right)}{(-1)^j 2^j (1+2 p)^2}=\frac{1}{8} \pi  \log
   \left(\frac{2 \pi  \exp (\gamma ) \Gamma \left(-\frac{1}{4}\right)^2}{9 x i \Gamma
   \left(-\frac{3}{4}\right)^2}\right)
\end{multline}
\end{example}
\begin{example}
\begin{multline}
\prod _{j=0}^{\infty } \prod _{p=0}^{\infty } \exp \left(\frac{(-1)^{-j+p} 2^{3-j} \binom{2}{j} \left((3+2
   p)^{-1+j} \log ((3+2 p) x)-(1-2 p)^{-1+j} \log (x-2 p x)\right)}{(1+2 p)^2 \pi }\right)\\
=-\frac{2 i e^{\gamma } \pi 
   \Gamma \left(-\frac{1}{4}\right)^2}{9 x \Gamma \left(-\frac{3}{4}\right)^2}
\end{multline}
\end{example}
\begin{example}
\begin{multline}
\frac{e^{a x} \cos (x \sin (t))}{1-b x}=\sum _{n=0}^{\infty }
   \frac{\left(e^{-b x} x\right)^n \left((a+b n-i \sin (t))^n+(a+b n+i \sin
   (t))^n\right)}{2 n!}
\end{multline}
where $f(x)=e^{ax}\cos(x\sin(t))$ and $\phi(x)=e^{bx}$.
\end{example}
\begin{example}
\begin{multline}
\frac{2 e^{a+(a+b m) x}}{(b x-1) \left(x^2-1\right)}\\
=\sum _{n=0}^{\infty
   } \frac{(-1)^{-n} e^{-b (m+n)} \left(\frac{x}{e^{b x}}\right)^n \left(\Gamma
   (1+n,-a-b (m+n))+(-1)^n e^{2 (a+b (m+n))} \Gamma (1+n,a+b
   (m+n))\right)}{\Gamma (1+n)}
\end{multline}
where $f(x)=\frac{e^{a x} x^m}{1-x^2}$ and $\phi(x)=e^{bx}$.
\end{example}
\begin{example}
\begin{multline}
\frac{e^{b x} \left(-e^{b-b m} \text{Ei}(a+b (-1+m))+e^{2 a+b (-1+m)}
   \text{Ei}(-a+b-b m)\right)}{x}\\+\frac{e^{b x} \left(-e^{2 a+b (-1+m)}
   \text{Ei}((a+b (-1+m)) (-1+x))+e^{b-b m} \text{Ei}((a+b (-1+m))
   (1+x))\right)}{x}\\
=\sum _{n=0}^{\infty } \frac{(-1)^{-n} e^{-b (m+n)}
   \left(\frac{x}{e^{b x}}\right)^n \left(\Gamma (1+n,-a-b (m+n))+(-1)^n e^{2
   (a+b (m+n))} \Gamma (1+n,a+b (m+n))\right)}{(n+1)!}
\end{multline}
\end{example}
\begin{example}
\begin{multline}
\frac{(1-x)^a}{\left(1+x^3 z^3\right)^{b+1}}
=\sum _{n=0}^{\infty } \frac{(1+a-n)_n (-x)^n }{3 n!}\\ \times
\left(\,
   _2F_1(1+b,-n;1+a-n;-z)+\, _2F_1\left(1+b,-n;1+a-n;\frac{1}{2} \left(1-i \sqrt{3}\right) z\right)\right. \\ \left.
+\,
   _2F_1\left(1+b,-n;1+a-n;\frac{1}{2} \left(1+i \sqrt{3}\right) z\right)\right)
\end{multline}
where $|Re(x)|<1,|Re(z)|<1,f(x)=(1-x)^a \left(x^3 z^3+1\right)^{-b-1},\phi(x)=1$.
\end{example}
\begin{example}
\begin{multline}
\frac{(1-x)^a}{\left(1+x^3 z^3\right)^{b+1}}
=\sum _{n=0}^{\infty } \frac{(1+a-n)_n (-x)^n }{3 n!}\\ \times
\left(\,
   _2F_1(1+b,-n;1+a-n;-z)+\, _2F_1\left(1+b,-n;1+a-n;\frac{1}{2} \left(1-i \sqrt{3}\right) z\right)\right. \\ \left.
+\,
   _2F_1\left(1+b,-n;1+a-n;\frac{1}{2} \left(1+i \sqrt{3}\right) z\right)\right)
\end{multline}
where $|Re(x)|<1,|Re(z)|<1$.
\end{example}
\begin{example}
\begin{multline}
\frac{e^{a x} (1+x)^{1+d}}{1-x (-1+b+c+b x)}\\
=\sum _{n=0}^{\infty } \frac{\left(e^{-b x} x (1+x)^{-c}\right)^n
   \Gamma (1+d+c n) \, _1F_1(-n;1+d+(-1+c) n;-a-b n)}{\Gamma (1+n) \Gamma (1+d-n+c n)}
\end{multline}
where $|Re(a)|<1,|Re(b)|<1,|Re(c)|<1,|Re(d)|<1,|Re(x)|<1$.
\end{example}
\begin{example}
\begin{multline}
\frac{e^{a x} x^m (1+x)^{1+d} \left(e^{-b x} x (1+x)^{-c}\right)^{-m}}{1-x (-1+b+c+b x)}\\
=\sum _{n=0}^{\infty }
   \frac{\left(\frac{x}{e^{b x} (1+x)^c}\right)^n (\Gamma (1+d+c (m+n)) \, _1F_1(-n;1+d-n+c (m+n);-a-b (m+n)))}{\Gamma
   (1+n) \Gamma (1+d-n+c (m+n))}
\end{multline}
where $|Re(a)|<1,|Re(b)|<1,|Re(c)|<1,|Re(d)|<1,|Re(m)|<1,|Re(x)|<1$.
\end{example}
\begin{example}
\begin{multline}
\frac{x^m (1+x)^b (1+x z)}{1+x (1+a+a x) z}=\sum _{n=0}^{\infty } \frac{\left(\frac{x (1+x
   z)^a}{1+x}\right)^{m+n} \Gamma (b+m+n) \, _2F_1(-n,a (m+n);b+m;z)}{\Gamma (b+m) \Gamma (1+n)}
\end{multline}
where $|Re(a)|<1,|Re(b)|<1,|Re(m)|<1,|Re(z)|<1,|Re(x)|<1$.
\end{example}
\begin{example}
\begin{equation}
\frac{(1-x z)^{1-c} \left((1-x z)^b\right)^{-m}}{1-(1+b) x z}=\sum _{n=0}^{\infty } \frac{z^n x^n \left((1-x
   z)^b\right)^n \Gamma (c+n+b (m+n))}{\Gamma (1+n) \Gamma (c+b (m+n))}
\end{equation}
where $|Re(b)|<1,|Re(c)|<1,|Re(m)|<1,|Re(z)|<1,|Re(x)|<1$.
\end{example}
\begin{example}
\begin{multline}
\frac{\left(1-(x z)^m\right)^{1-b}}{1-(1+a m) (x z)^m}
=\frac{3+e^{i m \pi }+4 \left\lfloor \frac{1}{2}
   (-1+m)\right\rfloor }{2 m}\\
+\sum _{n=1}^{\infty } \frac{(z x)^n \left(\left(1-(x z)^m\right)^a\right)^n \Gamma
   \left(b+\left(a+\frac{1}{m}\right) n\right) }{2 n \Gamma
   \left(\frac{n}{m}\right) \Gamma (b+a n)}\\ \times \left((-1)^n+(-1)^{m+n}+2 \csc \left(\frac{n \pi }{m}\right) \sin
   \left(\frac{\pi  \left(n+2 n \left\lfloor \frac{1}{2} (-1+m)\right\rfloor \right)}{m}\right)\right)
\end{multline}
where $|Re(z)|<1,|Re(x)|<1,m\in\mathbb{Z_{+}}$.
\end{example}
\begin{example}
\begin{multline}
\int_0^1 \frac{\left(1-(x z)^m\right)^{1-b}}{1-(1+a m) (x z)^m} \, dz=\frac{3+e^{i m \pi }+4 \left\lfloor
   \frac{1}{2} (-1+m)\right\rfloor }{2 m}\\
-\sum _{n=1}^{\infty } \frac{x^n \left(\left(1-x^m\right)^a\right)^n
   \left(-1+x^m\right) \Gamma \left(b+\left(a+\frac{1}{m}\right) n\right) }{2 n (1+n) \Gamma
   \left(\frac{n}{m}\right) \Gamma (b+a n)}\\
\times \, _2F_1\left(1,\frac{1+m+n+a m
   n}{m};\frac{1+m+n}{m};x^m\right)\\
 \left((-1)^n+(-1)^{m+n}+2 \csc \left(\frac{n \pi }{m}\right) \sin \left(\frac{\pi 
   \left(n+2 n \left\lfloor \frac{1}{2} (-1+m)\right\rfloor \right)}{m}\right)\right)
\end{multline}
\end{example}
\begin{example}
\begin{multline}
\left(-1+a e^x\right)^m=(-1+a)^m+(-1)^{m+1} \sum _{n=1}^{\infty } \sum
   _{j=1}^m \frac{a x^n (-1)^{1+j} a^{-1+j} j^n \binom{m}{j}}{n!}
\end{multline}
where $|Re(a)|<1,|Re(x)|<1,m\in\mathbb{Z_{+}}$.
\end{example}
\section{A short note on a definite integral.}
\begin{example}
\begin{multline}
\int_0^{\infty } \frac{\log ^k(a x) x^{-1+m} (1+c x)^u}{\left(1+b x^n\right)^p} \, dx\\
=\sum _{j=0}^{\infty }
   \sum _{r=0}^{p-1} \sum _{l=0}^{\infty } \frac{i (-1)^p i^{k-l} b^{-\frac{j+m}{n}} e^{\frac{i (j+m) \pi }{n}} c^j
   (j+m-n)^{-l+r} n^{-1-k+l-r} (2 \pi )^{1+k-l} \binom{r}{l} \binom{u}{j} }{\Gamma (p)}\\
\times \Phi \left(e^{\frac{2 i (j+m) \pi
   }{n}},-k+l,\frac{\pi -i n \left(\log (a)+\log \left(b^{-1/n}\right)\right)}{2 \pi }\right) (1+k-l)_l
   S_{-1+p}^{(r)}
\end{multline}
where $Re(a)<0,Re(m)>0$.
\end{example}
\begin{example}
\begin{multline}
B(m,p)
=\sum _{j=0}^{\infty } \sum _{r=0}^{p-1} \sum _{l=0}^{\infty } \frac{i (-1)^p i^{-l} e^{i (j+m) \pi }
   (-1+j+m)^{-l+r} (2 \pi )^{1-l} \binom{-m}{j} \binom{r}{l} }{\Gamma (p)}\\ \times \Phi \left(e^{2 i (j+m) \pi },l,\frac{1}{2}\right)
   (1-l)_l S_{-1+p}^{(r)}
\end{multline}
where $Re(m)<-\pi$.
\end{example}
%%
%\begin{example}
%
%\end{example}
%%%
\section{A table of formulae involving the $arcsin(x)$ function  and special functions.}
In this section we derive definite integrals and series involving the arcsin function in terms of special functions. In this work we have taken a combinatorial approach to deriving the $arcsin$ function in terms of the hypergeometric function. This process yielded a new special function forms for this trigonometric function not seen in current literature to the best of our knowledge. This form allowed us to write down definite integrals involving elliptic functions and Catalan's constant. The $arcsin$ function is used in the special integral Owen T-function \cite{nikolay}, pp.137. when dealing with wireless communications theory. 
\begin{lemma}\label{eq:arcsin}
\begin{multline}
\sinh ^{-1}(x)=-\frac{x^3 \, _2F_1\left(\frac{3}{4},\frac{3}{4};\frac{7}{4};\frac{1}{2}
   \left(1-\sqrt{1-x^4}\right)\right)}{3 \sqrt[4]{2} \left(1+\sqrt{1-x^4}\right)^{3/4}}+\frac{\sqrt[4]{2} x \,
   _2F_1\left(\frac{1}{4},\frac{1}{4};\frac{5}{4};\frac{1}{2}-\frac{\sqrt{1-x^4}}{2}\right)}{\sqrt[4]{1+\sqrt{1-x^4}}}
\end{multline}
Use (\ref{eq:lag1}) where $f(x)=\log \left(\sqrt{x^2+1}+x\right)$ and $\phi(x)=1$.
\end{lemma}
\section{Table of results}
The proceeding formulae are derived using the cited equations applied to equation (\ref{eq:arcsin}). 
\begin{example}
\begin{equation}
\pi =\frac{6 \sqrt{2 \left(4+\sqrt{15}\right)} \, _2F_1\left(\frac{1}{4},\frac{1}{4};\frac{5}{4};\frac{1}{2}-\frac{\sqrt{15}}{8}\right)+\,
   _2F_1\left(\frac{3}{4},\frac{3}{4};\frac{7}{4};\frac{1}{2}-\frac{\sqrt{15}}{8}\right)}{\left(2 \left(4+\sqrt{15}\right)\right)^{3/4}}
\end{equation}
\end{example}
\begin{example}
Catalan's constant $C$,
\begin{equation}
C=\int_0^1 \frac{-\frac{x^3 \, _2F_1\left(\frac{3}{4},\frac{3}{4};\frac{7}{4};\frac{1}{2}
   \left(1-\sqrt{1-x^4}\right)\right)}{3 \sqrt[4]{2} \left(1+\sqrt{1-x^4}\right)^{3/4}}+\frac{\sqrt[4]{2} x \,
   _2F_1\left(\frac{1}{4},\frac{1}{4};\frac{5}{4};\frac{1}{2}-\frac{\sqrt{1-x^4}}{2}\right)}{\sqrt[4]{1+\sqrt{1-x^4}}}}{
   \sqrt{1-x^2}} \, dx
\end{equation}
\end{example}
\begin{example}
From Eq. (4.521.1) in \cite{grad}
\begin{multline}
\frac{1}{2} \pi  \log (2)\\
=\int_0^1 \frac{1}{6 x \left(1+\sqrt{1-x^4}\right)^{3/4}}\left(2^{3/4} x^3 \, _2F_1\left(\frac{3}{4},\frac{3}{4};\frac{7}{4};\frac{1}{2}
   \left(1-\sqrt{1-x^4}\right)\right)\right. \\ \left.
   +6 \sqrt[4]{2} x \sqrt{1+\sqrt{1-x^4}} \,
   _2F_1\left(\frac{1}{4},\frac{1}{4};\frac{5}{4};\frac{1}{2}-\frac{\sqrt{1-x^4}}{2}\right)\right) \,dx
\end{multline}
\end{example}
\begin{example}
From Eq. (4.521.3) in \cite{grad}
\begin{multline}
\frac{\pi  \left(\log (4)+\log (1+q)-2 \log \left(1+\sqrt{1+q}\right)\right)}{4 q}\\
=\int_0^1 \frac{x}{6 \left(1+q x^2\right)
   \left(1+\sqrt{1-x^4}\right)^{3/4}} \left(2^{3/4} x^3 \, _2F_1\left(\frac{3}{4},\frac{3}{4};\frac{7}{4};\frac{1}{2}
   \left(1-\sqrt{1-x^4}\right)\right)\right. \\ \left.
   +6 \sqrt[4]{2} x \sqrt{1+\sqrt{1-x^4}} \,
   _2F_1\left(\frac{1}{4},\frac{1}{4};\frac{5}{4};\frac{1}{2}-\frac{\sqrt{1-x^4}}{2}\right)\right) \, dx
\end{multline}
\end{example}
\begin{example}
From Eq. (4.521.4) in \cite{grad}
\begin{multline}
-\frac{\pi  \left(\log \left(4-4 p^2\right)-2 \log \left(1+\sqrt{1-p^2}\right)\right)}{4 p^2}\\
=\int_0^1 \frac{x}{6 \left(1-p^2 x^2\right)
   \left(1+\sqrt{1-x^4}\right)^{3/4}} \left(2^{3/4} x^3 \,
   _2F_1\left(\frac{3}{4},\frac{3}{4};\frac{7}{4};\frac{1}{2} \left(1-\sqrt{1-x^4}\right)\right)\right. \\ \left.
   +6 \sqrt[4]{2} x \sqrt{1+\sqrt{1-x^4}}
   \, _2F_1\left(\frac{1}{4},\frac{1}{4};\frac{5}{4};\frac{1}{2}-\frac{\sqrt{1-x^4}}{2}\right)\right) \, dx 
\end{multline}
\end{example}
\begin{example}
From Eq. (4.521.6) in \cite{grad}
\begin{multline}
-\frac{1}{4} \pi  \left(\log (1+q)-2 \log \left(1+\sqrt{1+q}\right)\right)\\
=\int_0^1 \frac{1}{6
   x \left(1+q x^2\right) \left(1+\sqrt{1-x^4}\right)^{3/4}}\left(2^{3/4} x^3 \, _2F_1\left(\frac{3}{4},\frac{3}{4};\frac{7}{4};\frac{1}{2} \left(1-\sqrt{1-x^4}\right)\right)\right. \\ \left.
   +6
   \sqrt[4]{2} x \sqrt{1+\sqrt{1-x^4}} \, _2F_1\left(\frac{1}{4},\frac{1}{4};\frac{5}{4};\frac{1}{2}-\frac{\sqrt{1-x^4}}{2}\right)\right) \, dx 
\end{multline}
\end{example}
\begin{example}
From Eq. (4.521.7) in \cite{grad}
\begin{multline}
\frac{\pi  \left(1+q-\sqrt{1+q}\right)}{4 q (1+q)^{3/2}}\\
=\int_0^1 \frac{x}{6 \left(1+q x^2\right)^2
   \left(1+\sqrt{1-x^4}\right)^{3/4}} \left(2^{3/4} x^3 \,
   _2F_1\left(\frac{3}{4},\frac{3}{4};\frac{7}{4};\frac{1}{2} \left(1-\sqrt{1-x^4}\right)\right)\right. \\ \left.
   +6 \sqrt[4]{2} x \sqrt{1+\sqrt{1-x^4}}
   \, _2F_1\left(\frac{1}{4},\frac{1}{4};\frac{5}{4};\frac{1}{2}-\frac{\sqrt{1-x^4}}{2}\right)\right) \, dx
\end{multline}
\end{example}
\begin{example}
From Eq. (4.522.2) in \cite{grad}
\begin{multline}
\frac{-4 \left(-2+k^2\right) E\left(k^2\right)+\left(-1+k^2\right) \left(3 \sqrt{1-k^2} \pi +2 K\left(k^2\right)\right)}{18
   k^2}\\
=\int_0^1 \frac{x \sqrt{1-k^2 x^2}}{6
   \left(1+\sqrt{1-x^4}\right)^{3/4}} \left(2^{3/4} x^3 \,
   _2F_1\left(\frac{3}{4},\frac{3}{4};\frac{7}{4};\frac{1}{2} \left(1-\sqrt{1-x^4}\right)\right)\right. \\ \left.
   +6 \sqrt[4]{2} x \sqrt{1+\sqrt{1-x^4}}
   \, _2F_1\left(\frac{1}{4},\frac{1}{4};\frac{5}{4};\frac{1}{2}-\frac{\sqrt{1-x^4}}{2}\right)\right) \, dx
\end{multline}
\end{example}
\begin{example}
From Eq. (4.522.3) in \cite{grad}
\begin{multline}
\frac{\sqrt{p^2} \left(-4 \left(k^2+2 p^2\right) E\left(-\frac{k^2}{p^2}\right)+\left(k^2+p^2\right) \left(3
   \sqrt{1+\frac{k^2}{p^2}} \pi +2 K\left(-\frac{k^2}{p^2}\right)\right)\right)}{18 k^2}\\
   =\int_0^1
   \frac{x \sqrt{p^2+k^2 x^2} }{6
   \left(1+\sqrt{1-x^4}\right)^{3/4}}\left(2^{3/4} x^3 \, _2F_1\left(\frac{3}{4},\frac{3}{4};\frac{7}{4};\frac{1}{2}
   \left(1-\sqrt{1-x^4}\right)\right)\right. \\ \left.
   +6 \sqrt[4]{2} x \sqrt{1+\sqrt{1-x^4}} \,
   _2F_1\left(\frac{1}{4},\frac{1}{4};\frac{5}{4};\frac{1}{2}-\frac{\sqrt{1-x^4}}{2}\right)\right) \, dx 
   \end{multline}
\end{example}
\begin{example}
From Eq. (4.522.4) in \cite{grad}
\begin{multline}
-\frac{\sqrt{1-k^2} \pi -2 E\left(k^2\right)}{2 k^2}\\
=\int_0^1
   \frac{x}{6 \sqrt{1-k^2 x^2}
   \left(1+\sqrt{1-x^4}\right)^{3/4}} \left(2^{3/4} x^3 \, _2F_1\left(\frac{3}{4},\frac{3}{4};\frac{7}{4};\frac{1}{2} \left(1-\sqrt{1-x^4}\right)\right)\right. \\ \left.
   +6
   \sqrt[4]{2} x \sqrt{1+\sqrt{1-x^4}} \,
   _2F_1\left(\frac{1}{4},\frac{1}{4};\frac{5}{4};\frac{1}{2}-\frac{\sqrt{1-x^4}}{2}\right)\right) \, dx
\end{multline}
\end{example}
\begin{example}
From Eq. (4.522.6) in \cite{grad}
\begin{multline}
\frac{\sqrt{p^2} \left(\sqrt{1+\frac{k^2}{p^2}} \pi -2 E\left(-\frac{k^2}{p^2}\right)\right)}{2 k^2}\\
=\int_0^1 \frac{x}{6 \sqrt{p^2+k^2 x^2}
   \left(1+\sqrt{1-x^4}\right)^{3/4}} \left(2^{3/4} x^3 \, _2F_1\left(\frac{3}{4},\frac{3}{4};\frac{7}{4};\frac{1}{2}
   \left(1-\sqrt{1-x^4}\right)\right)\right. \\ \left.
   +6 \sqrt[4]{2} x \sqrt{1+\sqrt{1-x^4}} \,
   _2F_1\left(\frac{1}{4},\frac{1}{4};\frac{5}{4};\frac{1}{2}-\frac{\sqrt{1-x^4}}{2}\right)\right) \, dx
\end{multline}
\end{example}
\begin{example}
From Eq. (4.522.8) in \cite{grad}. The equation given is in error.
\begin{multline}
\int_0^1 \frac{x \sin ^{-1}(x)}{\sqrt{1-x^2} \left(x^2-\cos ^2(\lambda )\right)} \, dx\\
=\int_0^1 \frac{x}{6 \sqrt{1-x^2}
   \left(1+\sqrt{1-x^4}\right)^{3/4} \left(x^2-\cos ^2(\lambda )\right)} \left(2^{3/4} x^3 \,
   _2F_1\left(\frac{3}{4},\frac{3}{4};\frac{7}{4};\frac{1}{2} \left(1-\sqrt{1-x^4}\right)\right)\right. \\ \left.
   +6 \sqrt[4]{2} x \sqrt{1+\sqrt{1-x^4}}
   \, _2F_1\left(\frac{1}{4},\frac{1}{4};\frac{5}{4};\frac{1}{2}-\frac{\sqrt{1-x^4}}{2}\right)\right) \, dx
\end{multline}
\end{example}
\begin{example}
From Eq. (4.522.9) in \cite{grad}
\begin{multline}
\int_0^1 \frac{x \sin ^{-1}(x)}{\sqrt{\left(1-x^2\right) \left(1-k^2 x^2\right)}} \, dx\\
=\int_0^1 \frac{x}{6 \sqrt{\left(1-x^2\right)
   \left(1-k^2 x^2\right)} \left(1+\sqrt{1-x^4}\right)^{3/4}} \left(2^{3/4} x^3 \,
   _2F_1\left(\frac{3}{4},\frac{3}{4};\frac{7}{4};\frac{1}{2} \left(1-\sqrt{1-x^4}\right)\right)\right. \\ \left.
   +6 \sqrt[4]{2} x \sqrt{1+\sqrt{1-x^4}}
   \, _2F_1\left(\frac{1}{4},\frac{1}{4};\frac{5}{4};\frac{1}{2}-\frac{\sqrt{1-x^4}}{2}\right)\right) \, dx
\end{multline}
\end{example}
\begin{example}
From Eq. (4.523.1) in \cite{grad}
\begin{multline}
\frac{\pi }{2+4 n}-\frac{\sqrt{\pi } \Gamma (1+n)}{(1+2 n)^2 \Gamma \left(\frac{1}{2}+n\right)}\\
=\int_0^1 \frac{x^{2 n}}{6
   \left(1+\sqrt{1-x^4}\right)^{3/4}} \left(2^{3/4} x^3 \, _2F_1\left(\frac{3}{4},\frac{3}{4};\frac{7}{4};\frac{1}{2}
   \left(1-\sqrt{1-x^4}\right)\right)\right. \\ \left.
   +6 \sqrt[4]{2} x \sqrt{1+\sqrt{1-x^4}} \,
   _2F_1\left(\frac{1}{4},\frac{1}{4};\frac{5}{4};\frac{1}{2}-\frac{\sqrt{1-x^4}}{2}\right)\right) \, dx 
\end{multline}
\end{example}
\begin{example}
From Eq. (4.523.2) in \cite{grad}
\begin{multline}
\frac{n \pi  \Gamma (n)-\sqrt{\pi } \Gamma \left(\frac{1}{2}+n\right)}{4 n^2 \Gamma (n)}\\
=\int_0^1
   \frac{x^{-1+2 n} }{6
   \left(1+\sqrt{1-x^4}\right)^{3/4}}\left(2^{3/4} x^3 \, _2F_1\left(\frac{3}{4},\frac{3}{4};\frac{7}{4};\frac{1}{2}
   \left(1-\sqrt{1-x^4}\right)\right)\right. \\ \left.
   +6 \sqrt[4]{2} x \sqrt{1+\sqrt{1-x^4}} \,
   _2F_1\left(\frac{1}{4},\frac{1}{4};\frac{5}{4};\frac{1}{2}-\frac{\sqrt{1-x^4}}{2}\right)\right) \, dx
\end{multline}
\end{example}
\begin{example}
From Eq. (4.551.1) in \cite{grad}
\begin{multline}
-\frac{\pi  \left(e^{-b}-I_0(b)+\pmb{L}_0(b)\right)}{2 b}\\
=\int_0^1 \frac{e^{-b x} }{6
   \left(1+\sqrt{1-x^4}\right)^{3/4}}\left(2^{3/4} x^3 \,
   _2F_1\left(\frac{3}{4},\frac{3}{4};\frac{7}{4};\frac{1}{2} \left(1-\sqrt{1-x^4}\right)\right)\right. \\ \left.
   +6 \sqrt[4]{2} x \sqrt{1+\sqrt{1-x^4}}
   \, _2F_1\left(\frac{1}{4},\frac{1}{4};\frac{5}{4};\frac{1}{2}-\frac{\sqrt{1-x^4}}{2}\right)\right) \, dx
\end{multline}
\end{example}
\begin{example}
From Eq. (4.551.2) in \cite{grad}. The modified Struve function $\pmb{L}_n(z)$.
\begin{multline}
\frac{3 \left(-2+b^2\right) \pi  \pmb{L}_1(b)+b \left(4 b+3 \pi  \left(-\left((1+b) e^{-b}\right)+I_0(b)-b I_1(b)\right)-3 \pi
    \pmb{L}_2(b)\right)}{6 b^3}\\
    =\int_0^1 \frac{e^{-b x} x}{6
   \left(1+\sqrt{1-x^4}\right)^{3/4}} \left(2^{3/4} x^3 \,
   _2F_1\left(\frac{3}{4},\frac{3}{4};\frac{7}{4};\frac{1}{2} \left(1-\sqrt{1-x^4}\right)\right)\right. \\ \left.
   +6 \sqrt[4]{2} x \sqrt{1+\sqrt{1-x^4}}
   \, _2F_1\left(\frac{1}{4},\frac{1}{4};\frac{5}{4};\frac{1}{2}-\frac{\sqrt{1-x^4}}{2}\right)\right) \, dx 
\end{multline}
\end{example}
\begin{example}
From Eq. (4.571) in \cite{grad}
\begin{multline}
\int_0^{\frac{\pi }{2}} \sin ^{-1}(k \sin (x)) \, dx\\
=\int_0^{\frac{\pi }{2}} \frac{1}{6 \left(1+\sqrt{1-k^4 \sin ^4(x)}\right)^{3/4}}\left(2^{3/4} k^3 \,
   _2F_1\left(\frac{3}{4},\frac{3}{4};\frac{7}{4};\frac{1}{2} \left(1-\sqrt{1-k^4 \sin ^4(x)}\right)\right) \sin ^3(x)\right. \\ \left.
+6 \sqrt[4]{2} k
   \, _2F_1\left(\frac{1}{4},\frac{1}{4};\frac{5}{4};\frac{1}{2}-\frac{1}{2} \sqrt{1-k^4 \sin ^4(x)}\right) \sin (x)
   \sqrt{1+\sqrt{1-k^4 \sin ^4(x)}}\right) \, dx
\end{multline}
\end{example}
\begin{example}
From Eq. (4.591.1) in \cite{grad}
\begin{multline}
2-\frac{\pi }{2}-\log (2)\\
=\int_0^1 \frac{1}{6
   \left(1+\sqrt{1-x^4}\right)^{3/4}}\left(2^{3/4} x^3 \, _2F_1\left(\frac{3}{4},\frac{3}{4};\frac{7}{4};\frac{1}{2}
   \left(1-\sqrt{1-x^4}\right)\right)\right. \\ \left.
   +6 \sqrt[4]{2} x \sqrt{1+\sqrt{1-x^4}} \,
   _2F_1\left(\frac{1}{4},\frac{1}{4};\frac{5}{4};\frac{1}{2}-\frac{\sqrt{1-x^4}}{2}\right)\right) \log (x) \, dx
\end{multline}
\end{example}
\begin{example}
From Eq. (4.523.2) in \cite{grad}
\begin{multline}
\frac{2^{-1-n} \pi  ((-2+n)\text{!!})^2}{\left(\frac{1+n}{2}!\right)^2}
=\int_{-1}^1 \sin ^{-1}(x) P_n(x) \,
   dx\\
   =\int_{-1}^1 \frac{1}{6
   \left(1+\sqrt{1-x^4}\right)^{3/4}}\left(2^{3/4} x^3 \, _2F_1\left(\frac{3}{4},\frac{3}{4};\frac{7}{4};\frac{1}{2}
   \left(1-\sqrt{1-x^4}\right)\right)\right. \\ \left.
   +6 \sqrt[4]{2} x \sqrt{1+\sqrt{1-x^4}} \,
   _2F_1\left(\frac{1}{4},\frac{1}{4};\frac{5}{4};\frac{1}{2}-\frac{\sqrt{1-x^4}}{2}\right)\right) P_n(x) \, dx
\end{multline}
\end{example}
\begin{example}
From Eq. (8.415.1) in \cite{grad}
\begin{multline}
Y_0(t)+\frac{4 \int_1^{\infty } \frac{\log \left(x+\sqrt{-1+x^2}\right) \sin (t x)}{\sqrt{-1+x^2}} \, dx}{\pi
   ^2}
=\int_0^1 \frac{4 \sin ^{-1}(x) \sin (t x)}{\pi ^2 \sqrt{1-x^2}} \, dx\\
=\int_0^1 \frac{2}{3 \pi ^2
   \sqrt{1-x^2} \left(1+\sqrt{1-x^4}\right)^{3/4}} \left(2^{3/4} x^3 \,
   _2F_1\left(\frac{3}{4},\frac{3}{4};\frac{7}{4};\frac{1}{2} \left(1-\sqrt{1-x^4}\right)\right)\right. \\ \left.
   +6 \sqrt[4]{2} x \sqrt{1+\sqrt{1-x^4}}
   \, _2F_1\left(\frac{1}{4},\frac{1}{4};\frac{5}{4};\frac{1}{2}-\frac{\sqrt{1-x^4}}{2}\right)\right) \sin (t x) \, dx
\end{multline}
\end{example}
\begin{example}
From Eq. (2.7.3.1) in \cite{prud1}
\begin{multline}
e^{-i t} \left(-1+\sqrt{1-e^{2 i t}}+e^{i t} \sin ^{-1}\left(e^{i t}\right)\right)\\
=\int_0^1 \frac{1}{6 \left(1+\sqrt{1-e^{4 i t} x^4}\right)^{3/4}}\left(2^{3/4} e^{3 i t} x^3 \,
   _2F_1\left(\frac{3}{4},\frac{3}{4};\frac{7}{4};\frac{1}{2} \left(1-\sqrt{1-e^{4 i t} x^4}\right)\right)\right. \\ \left.
+6 \sqrt[4]{2} e^{i t} x
   \sqrt{1+\sqrt{1-e^{4 i t} x^4}} \, _2F_1\left(\frac{1}{4},\frac{1}{4};\frac{5}{4};\frac{1}{2}-\frac{1}{2} \sqrt{1-e^{4 i t}
   x^4}\right)\right) \, dx
\end{multline}
\end{example}
\begin{example}
From Eq. (2.7.3.2) in \cite{prud1}
\begin{multline}
\frac{\pi }{2 \alpha }-\frac{\sqrt{\pi } \Gamma \left(\frac{1+\alpha }{2}\right)}{\alpha ^2 \Gamma \left(\frac{\alpha
   }{2}\right)}\\
=\int_0^1 \frac{x^{-1+\alpha } }{6
   \left(1+\sqrt{1-x^4}\right)^{3/4}}\left(2^{3/4} x^3 \,
   _2F_1\left(\frac{3}{4},\frac{3}{4};\frac{7}{4};\frac{1}{2} \left(1-\sqrt{1-x^4}\right)\right)\right. \\ \left.
+6 \sqrt[4]{2} x \sqrt{1+\sqrt{1-x^4}}
   \, _2F_1\left(\frac{1}{4},\frac{1}{4};\frac{5}{4};\frac{1}{2}-\frac{\sqrt{1-x^4}}{2}\right)\right)\, dx
\end{multline}
\end{example}
\begin{example}
From Eq. (2.7.3.6) in \cite{prud1}
\begin{multline}
\frac{J_0(n \pi )-\cos (n \pi )}{2 n}\\
=\int_0^1 \frac{1}{6\left(1+\sqrt{1-x^4}\right)^{3/4}}\left(2^{3/4} x^3 \,
   _2F_1\left(\frac{3}{4},\frac{3}{4};\frac{7}{4};\frac{1}{2} \left(1-\sqrt{1-x^4}\right)\right)\right. \\ \left.
+6 \sqrt[4]{2} x \sqrt{1+\sqrt{1-x^4}}
   \, _2F_1\left(\frac{1}{4},\frac{1}{4};\frac{5}{4};\frac{1}{2}-\frac{\sqrt{1-x^4}}{2}\right)\right) \sin (n \pi  x) \, dx
\end{multline}
\end{example}
\begin{example}
From Eq. (2.7.3.9) in \cite{prud1}
\begin{multline}
\frac{\pi  (-I_0(a)+\cosh (a))}{a}\\
=\int_{-1}^1 \frac{e^{a x} }{6
   \left(1+\sqrt{1-x^4}\right)^{3/4}}\left(2^{3/4} x^3 \,
   _2F_1\left(\frac{3}{4},\frac{3}{4};\frac{7}{4};\frac{1}{2} \left(1-\sqrt{1-x^4}\right)\right)\right. \\ \left.
+6 \sqrt[4]{2} x \sqrt{1+\sqrt{1-x^4}}
   \, _2F_1\left(\frac{1}{4},\frac{1}{4};\frac{5}{4};\frac{1}{2}-\frac{\sqrt{1-x^4}}{2}\right)\right) \, dx
\end{multline}
\end{example}
\begin{example}
\begin{multline}
2-2 C+\frac{1}{2} \pi  (-2+\log (8))\\
=\int_0^1 \frac{1}{6
   \left(1+\sqrt{1-x^4}\right)^{3/4}} \left(2^{3/4} x^3 \,
   _2F_1\left(\frac{3}{4},\frac{3}{4};\frac{7}{4};\frac{1}{2} \left(1-\sqrt{1-x^4}\right)\right)\right. \\ \left.
+6 \sqrt[4]{2} x \sqrt{1+\sqrt{1-x^4}}
   \, _2F_1\left(\frac{1}{4},\frac{1}{4};\frac{5}{4};\frac{1}{2}-\frac{\sqrt{1-x^4}}{2}\right)\right) \log (1+x)\, dx
\end{multline}
\end{example}
\begin{example}
From Eq. (2.6.1.1) in \cite{brychkov}
\begin{multline}
-\frac{i e^{-\frac{1}{2} i \pi  s} \Gamma \left(1-\frac{s}{2}\right) \Gamma \left(\frac{1+s}{2}\right)}{\sqrt{\pi }
   s^2}\\
   =\int_0^{\infty } \frac{x^{-1+s}}{6
   \left(1+\sqrt{1-x^4}\right)^{3/4}} \left(2^{3/4} x^3 \,
   _2F_1\left(\frac{3}{4},\frac{3}{4};\frac{7}{4};\frac{1}{2} \left(1-\sqrt{1-x^4}\right)\right)\right. \\ \left.
   +6 \sqrt[4]{2} x \sqrt{1+\sqrt{1-x^4}}
   \, _2F_1\left(\frac{1}{4},\frac{1}{4};\frac{5}{4};\frac{1}{2}-\frac{\sqrt{1-x^4}}{2}\right)\right) \, dx
   \end{multline}
\end{example}
\begin{example}
From Eq. (2.6.1.7) in \cite{brychkov}
\begin{multline}
\frac{i i^{-s} a^{-s}}{2 \sqrt{\pi }} \left(\Gamma \left(-\frac{s}{2}\right) \Gamma \left(\frac{1+s}{2}\right) \,
   _3F_2\left(\frac{1}{2},\frac{1}{2},\frac{1}{2}+\frac{s}{2};\frac{3}{2},1+\frac{s}{2};1\right)\right. \\ \left.
+\frac{2 \Gamma
   \left(\frac{3}{2}-\frac{s}{2}\right) \Gamma \left(\frac{s}{2}\right) \,
   _3F_2\left(\frac{1}{2},\frac{1}{2}-\frac{s}{2},\frac{1}{2}-\frac{s}{2};1-\frac{s}{2},\frac{3}{2}-\frac{s}{2};1\right)}{(-1+s)^2}\right)\\
=\int_0^{\infty }
   \frac{x^{-1+s}}{6 \sqrt{1-a^2 x^2}
   \left(1+\sqrt{1-a^4 x^4}\right)^{3/4}} \left(2^{3/4} a^3 x^3 \, _2F_1\left(\frac{3}{4},\frac{3}{4};\frac{7}{4};\frac{1}{2} \left(1-\sqrt{1-a^4
   x^4}\right)\right)\right. \\ \left.
+6 \sqrt[4]{2} a x \sqrt{1+\sqrt{1-a^4 x^4}} \,
   _2F_1\left(\frac{1}{4},\frac{1}{4};\frac{5}{4};\frac{1}{2}-\frac{1}{2} \sqrt{1-a^4 x^4}\right)\right) \, dx
\end{multline}
\end{example}
\begin{example}
From Eq. (2.6.1.18) in \cite{brychkov}
\begin{multline}
a^{-1+s-\rho } b B(s,1-s+\rho ) \, _4F_3\left(\frac{1}{2},\frac{1}{2},\frac{1}{2}-\frac{s}{2}+\frac{\rho
   }{2},1-\frac{s}{2}+\frac{\rho }{2};\frac{3}{2},\frac{1}{2}+\frac{\rho }{2},1+\frac{\rho
   }{2};\frac{b^2}{a^2}\right)\\
=\int_0^{\infty } x^{-1+s} (a+x)^{-\rho } \sin ^{-1}\left(\frac{b}{a+x}\right) \,
   dx\\
=\int_0^{\infty } \frac{x^{-1+s} (a+x)^{-\rho }}{6 \left(1+\sqrt{1-\frac{b^4}{(a+x)^4}}\right)^{3/4}} \left(\frac{2^{3/4} b^3 \,
   _2F_1\left(\frac{3}{4},\frac{3}{4};\frac{7}{4};\frac{1}{2} \left(1-\sqrt{1-\frac{b^4}{(a+x)^4}}\right)\right)}{(a+x)^3}\right. \\ \left.
+\frac{6
   \sqrt[4]{2} b \sqrt{1+\sqrt{1-\frac{b^4}{(a+x)^4}}} \, _2F_1\left(\frac{1}{4},\frac{1}{4};\frac{5}{4};\frac{1}{2}-\frac{1}{2}
   \sqrt{1-\frac{b^4}{(a+x)^4}}\right)}{a+x}\right) \, dx
\end{multline}
\end{example}
\begin{example}
From Eq. (2.6.1.19) in \cite{brychkov}
\begin{multline}
a^{-2+s-\rho } b B(s,2-s+\rho ) \, _4F_3\left(1,1,\frac{1}{2} (2-s+\rho ),\frac{1}{2} (3-s+\rho
   );\frac{3}{2},\frac{2}{2}+\frac{\rho }{2},\frac{\rho +3}{2};\frac{b^2}{a^2}\right)\\
=\int_0^{\infty } \frac{x^{s-1}
   \left(\frac{2^{3/4} b^3 \, _2F_1\left(\frac{3}{4},\frac{3}{4};\frac{7}{4};\frac{1}{2}
   \left(1-\sqrt{1-\frac{b^4}{(a+x)^4}}\right)\right)}{(a+x)^3}+\frac{6 \sqrt[4]{2} b
   \sqrt{1+\sqrt{1-\frac{b^4}{(a+x)^4}}} \, _2F_1\left(\frac{1}{4},\frac{1}{4};\frac{5}{4};\frac{1}{2}-\frac{1}{2}
   \sqrt{1-\frac{b^4}{(a+x)^4}}\right)}{a+x}\right)}{\left(\sqrt{(x+a)^2-b^2} (x+a)^{\rho }\right) \left(6
   \left(1+\sqrt{1-\frac{b^4}{(a+x)^4}}\right)^{3/4}\right)} \, dx
\end{multline}
\end{example}
\begin{example}
From Eq. (2.6.1.20) in \cite{brychkov}
\begin{multline}
a^{s-\rho } b B(s+1,-s+\rho ) \,
   _4F_3\left(\frac{1}{2},\frac{1}{2},\frac{s+1}{2},\frac{s+2}{2};\frac{3}{2},\frac{1}{2}+\frac{\rho }{2},\frac{\rho
   +2}{2};b^2\right)\\
=\int_0^{\infty } \frac{x^{s-1} \left(\frac{2^{3/4} b^3 x^3 \,
   _2F_1\left(\frac{3}{4},\frac{3}{4};\frac{7}{4};\frac{1}{2} \left(1-\sqrt{1-\frac{b^4
   x^4}{(a+x)^4}}\right)\right)}{(a+x)^3}+\frac{6 \sqrt[4]{2} b x \sqrt{1+\sqrt{1-\frac{b^4 x^4}{(a+x)^4}}} \,
   _2F_1\left(\frac{1}{4},\frac{1}{4};\frac{5}{4};\frac{1}{2}-\frac{1}{2} \sqrt{1-\frac{b^4
   x^4}{(a+x)^4}}\right)}{a+x}\right)}{(x+a)^{\rho } \left(6 \left(1+\sqrt{1-\frac{b^4
   x^4}{(a+x)^4}}\right)^{3/4}\right)} \, dx
\end{multline}
\end{example}
\begin{example}
From Eq. (2.6.1.21) in \cite{brychkov}
\begin{multline}
a^{s-\rho } b B(s+1,-s+\rho ) \, _4F_3\left(1,1,\frac{s+1}{2},\frac{s+2}{2};\frac{3}{2},\frac{1}{2}+\frac{\rho
   }{2},\frac{\rho +2}{2};b^2\right)\\
=\int_0^{\infty } \frac{\left(x^{s-1} (x+a)^{-\rho }\right) }{\sqrt{1-\left(\frac{b x}{x+a}\right)^2} \left(6 \left(1+\sqrt{1-\frac{b^4
   x^4}{(a+x)^4}}\right)^{3/4}\right)}\\
\left(\frac{2^{3/4} b^3
   x^3 \, _2F_1\left(\frac{3}{4},\frac{3}{4};\frac{7}{4};\frac{1}{2} \left(1-\sqrt{1-\frac{b^4
   x^4}{(a+x)^4}}\right)\right)}{(a+x)^3}\right. \\ \left.
+\frac{6 \sqrt[4]{2} b x \sqrt{1+\sqrt{1-\frac{b^4 x^4}{(a+x)^4}}} \,
   _2F_1\left(\frac{1}{4},\frac{1}{4};\frac{5}{4};\frac{1}{2}-\frac{1}{2} \sqrt{1-\frac{b^4
   x^4}{(a+x)^4}}\right)}{a+x}\right) \, dx
\end{multline}
\end{example}
\begin{example}
From Eq. (2.6.1.22) in \cite{brychkov}
\begin{multline}
\frac{1}{2} a^{s-2 \rho -1} b B\left(\frac{s+1}{2},\frac{1}{2} (1-s+2 \rho )\right)\\
 \,
   _4F_3\left(\frac{1}{2},\frac{1}{2},\frac{s+1}{2},\frac{1}{2} (1-s+2 \rho );\frac{3}{2},\frac{\rho +1}{2},\frac{\rho
   +2}{2};\left(\frac{b}{2 a}\right)^2\right)\\
=\int_0^{\infty } \frac{x^{s-1}}{\left(x^2+a^2\right)^{\rho } \left(6 \left(1+\sqrt{1-\frac{b^4
   x^4}{\left(a^2+x^2\right)^4}}\right)^{3/4}\right)}\\
 \left(\frac{2^{3/4} b^3 x^3 \,
   _2F_1\left(\frac{3}{4},\frac{3}{4};\frac{7}{4};\frac{1}{2} \left(1-\sqrt{1-\frac{b^4
   x^4}{\left(a^2+x^2\right)^4}}\right)\right)}{\left(a^2+x^2\right)^3}\right. \\ \left.
+\frac{6 \sqrt[4]{2} b x
   \sqrt{1+\sqrt{1-\frac{b^4 x^4}{\left(a^2+x^2\right)^4}}} \,
   _2F_1\left(\frac{1}{4},\frac{1}{4};\frac{5}{4};\frac{1}{2}-\frac{1}{2} \sqrt{1-\frac{b^4
   x^4}{\left(a^2+x^2\right)^4}}\right)}{a^2+x^2}\right) \, dx
\end{multline}
\end{example}
\begin{example}
From Eq. (2.6.1.23) in \cite{brychkov}
\begin{multline}
\frac{1}{2} \left(a^{s-2 \rho -1} b\right) B\left(\frac{s+1}{2},\frac{1}{2} (1-s+2 \rho )\right)\\ \,
   _4F_3\left(1,1,\frac{s+1}{2},\frac{1}{2} (1-s+2 \rho );\frac{3}{2},\frac{1}{2}+\frac{\rho }{2},\frac{\rho
   +2}{2};\left(\frac{b}{2 a}\right)^2\right)\\
=\int_0^{\infty } \frac{\left(x^{s-1} \left(x^2+a^2\right)^{-\rho }\right)
  }{\sqrt{1-\left(\frac{b x}{x^2+a^2}\right)^2} \left(6
   \left(1+\sqrt{1-\frac{b^4 x^4}{\left(a^2+x^2\right)^4}}\right)^{3/4}\right)}\\
 \left(\frac{2^{3/4} b^3 x^3 \, _2F_1\left(\frac{3}{4},\frac{3}{4};\frac{7}{4};\frac{1}{2} \left(1-\sqrt{1-\frac{b^4
   x^4}{\left(a^2+x^2\right)^4}}\right)\right)}{\left(a^2+x^2\right)^3}\right. \\ \left.
+\frac{6 \sqrt[4]{2} b x
   \sqrt{1+\sqrt{1-\frac{b^4 x^4}{\left(a^2+x^2\right)^4}}} \,
   _2F_1\left(\frac{1}{4},\frac{1}{4};\frac{5}{4};\frac{1}{2}-\frac{1}{2} \sqrt{1-\frac{b^4
   x^4}{\left(a^2+x^2\right)^4}}\right)}{a^2+x^2}\right) \, dx
\end{multline}
\end{example}
\begin{example}
From Eq. (2.6.1.24) in \cite{brychkov}
\begin{multline}
a^{s-\rho -1} b B(s,1-s+\rho ) \, _3F_2\left(1,1,1-s+\rho ;\frac{3}{2},\rho
   +1;\frac{b^2}{a}\right)\\
=\int_0^{\infty } \frac{\left((x+a)^{-\rho } x^{s-1}\right) }{\sqrt{a-b^2+x} \left(6
   \left(1+\sqrt{1-\frac{b^4}{(a+x)^2}}\right)^{3/4}\right)}\\
\left(\frac{2^{3/4} b^3 \,
   _2F_1\left(\frac{3}{4},\frac{3}{4};\frac{7}{4};\frac{1}{2}
   \left(1-\sqrt{1-\frac{b^4}{(a+x)^2}}\right)\right)}{(a+x)^{3/2}}\right. \\ \left.
+\frac{6 \sqrt[4]{2} b
   \sqrt{1+\sqrt{1-\frac{b^4}{(a+x)^2}}} \, _2F_1\left(\frac{1}{4},\frac{1}{4};\frac{5}{4};\frac{1}{2}-\frac{1}{2}
   \sqrt{1-\frac{b^4}{(a+x)^2}}\right)}{\sqrt{a+x}}\right) \, dx
\end{multline}
\end{example}
\begin{example}
From Eq. (42.1.4) in \cite{hansen}
\begin{multline}
\frac{3 \tan ^{-1}\left(\sinh \left(\frac{\pi  x}{2 \sqrt{1-x^2}}\right)\right)}{\sqrt[4]{2} x}
=\sum
   _{k=0}^{\infty } \frac{(-1)^k }{(1+2 k)^3 \left(1+\frac{x^2}{(1+2 k)^2}\right)^{3/2}
   \left(1+\sqrt{1-\frac{x^4}{\left((1+2 k)^2+x^2\right)^2}}\right)^{3/4}}\\
\left(\sqrt{2} x^2 \, _2F_1\left(\frac{3}{4},\frac{3}{4};\frac{7}{4};\frac{1}{2}
   \left(1-\sqrt{1-\frac{x^4}{\left((1+2 k)^2+x^2\right)^2}}\right)\right)\right. \\ \left.
+6 \left((1+2 k)^2+x^2\right)
   \sqrt{1+\sqrt{\frac{(1+2 k)^4+2 (x+2 k x)^2}{\left((1+2 k)^2+x^2\right)^2}}} \right. \\ \left.
\,
   _2F_1\left(\frac{1}{4},\frac{1}{4};\frac{5}{4};\frac{1}{2} \left(1-\sqrt{1-\frac{x^4}{\left((1+2
   k)^2+x^2\right)^2}}\right)\right)\right)
\end{multline}
\end{example}
\begin{example}
From Eq. (42.1.10) in \cite{hansen}
\begin{multline}
-\tan ^{-1}\left(\sinh \left(\frac{\pi  a}{y}\right) \csc \left(\frac{\pi  x}{2 a}\right)\right)+\tan
   ^{-1}\left(\frac{2 a^2}{x y}\right)\\
=\sum _{k=1}^{\infty } \frac{4 (-1)^k\left(1+\sqrt{1-\frac{256 a^8 x^4 y^4}{\left(16 a^8+8 a^4 \left(4 a^2 k^2+x^2\right)y^2+\left(-4 a^2 k^2+x^2\right)^2 y^4\right)^2}}\right)^{-3/4}}{3 \left(4
   a^4+\left(4 a^2 k^2-x^2\right) y^2\right)^3 \left(1+\frac{x^2}{\left(\frac{a^2}{y}+k^2 y-\frac{x^2 y}{4
   a^2}\right)^2}\right)^{3/2} }\\ \left(8\times 2^{3/4} a^6 x^3 y^3 \,_2F_1\left(\frac{3}{4},\frac{3}{4};\frac{7}{4};\frac{1}{2} \left(1-\sqrt{1-\frac{256 a^8 x^4 y^4}{\left(16 a^8+8a^4 \left(4 a^2 k^2+x^2\right) y^2+\left(-4 a^2 k^2+x^2\right)^2 y^4\right)^2}}\right)\right)\right. \\ \left.
+3 \sqrt[4]{2} a^2
   x y \left(16 a^4 x^2 y^2+\left(4 a^4+\left(4 a^2 k^2-x^2\right) y^2\right)^2\right) \right. \\ \left.
\sqrt{1+\sqrt{1-\frac{256
   a^8 x^4 y^4}{\left(16 a^8+8 a^4 \left(4 a^2 k^2+x^2\right) y^2+\left(-4 a^2 k^2+x^2\right)^2 y^4\right)^2}}}\right. \\ \left. \,
   _2F_1\left(\frac{1}{4},\frac{1}{4};\frac{5}{4};\frac{1}{2} \left(1-\sqrt{1-\frac{256 a^8 x^4 y^4}{\left(16 a^8+8
   a^4 \left(4 a^2 k^2+x^2\right) y^2+\left(-4 a^2 k^2+x^2\right)^2 y^4\right)^2}}\right)\right)\right)
\end{multline}
\end{example}
\begin{example}
From Eq. (42.3.1) in \cite{hansen}
\begin{multline}
-\tan ^{-1}\left(x^n\right)
=\sum _{k=0}^{n-1} \frac{(-1)^k \left(1+\sqrt{1-\frac{x^4 \sin ^4\left(\frac{\pi +2 k \pi }{2 n}\right)}{\left(1+x^2-2 x \cos \left(\frac{\pi +2 k \pi }{2 n}\right)\right)^2}}\right)^{-3/4}}{6 \left(1-x \cos \left(\frac{\pi
   +2 k \pi }{2 n}\right)\right)^3 \left(1+\frac{x^2 \sin ^2\left(\frac{\pi +2 k \pi }{2 n}\right)}{\left(-1+x \cos \left(\frac{\pi +2 k \pi }{2 n}\right)\right)^2}\right)^{3/2}
   }\\
\times \left(2^{3/4} x^3 \, _2F_1\left(\frac{3}{4},\frac{3}{4};\frac{7}{4};\frac{1}{8} \left(4-4 \sqrt{1-\frac{x^4 \sin
   ^4\left(\frac{\pi +2 k \pi }{2 n}\right)}{\left(1+x^2-2 x \cos \left(\frac{\pi +2 k \pi }{2 n}\right)\right)^2}}\right)\right) \sin ^3\left(\frac{\pi +2 k \pi }{2 n}\right)\right. \\ \left.
+6
   \sqrt[4]{2} x \left(1+x^2-2 x \cos \left(\frac{\pi +2 k \pi }{2 n}\right)\right) \right. \\ \left.
\, _2F_1\left(\frac{1}{4},\frac{1}{4};\frac{5}{4};\frac{1}{8} \left(4-4 \sqrt{1-\frac{x^4 \sin
   ^4\left(\frac{\pi +2 k \pi }{2 n}\right)}{\left(1+x^2-2 x \cos \left(\frac{\pi +2 k \pi }{2 n}\right)\right)^2}}\right)\right)\right. \\ \left.
 \sin \left(\frac{\pi +2 k \pi }{2 n}\right)
   \sqrt{1+\sqrt{1-\frac{x^4 \sin ^4\left(\frac{\pi +2 k \pi }{2 n}\right)}{\left(1+x^2-2 x \cos \left(\frac{\pi +2 k \pi }{2 n}\right)\right)^2}}}\right)
\end{multline}
\end{example}
\begin{example}
From Eq. (42.1.11) in \cite{hansen}
\begin{multline}
-\frac{\pi }{4}=\sum _{k=1}^{\infty } \frac{2 (-1)^k \sqrt[4]{2} }{3 \sqrt{1+\frac{4}{k^4}} k^2
   \left(4+k^4\right) \left(1+\sqrt{\frac{k^4 \left(8+k^4\right)}{\left(4+k^4\right)^2}}\right)^{3/4}}\\
\left(2 \sqrt{2} \, _2F_1\left(\frac{3}{4},\frac{3}{4};\frac{7}{4};\frac{1}{2} \left(1-\sqrt{\frac{k^4
   \left(8+k^4\right)}{\left(4+k^4\right)^2}}\right)\right)\right. \\ \left.
+3 \left(4+k^4\right) \sqrt{1+\sqrt{\frac{k^4 \left(8+k^4\right)}{\left(4+k^4\right)^2}}} \,
   _2F_1\left(\frac{1}{4},\frac{1}{4};\frac{5}{4};\frac{1}{2}-\frac{1}{2} \sqrt{\frac{k^4 \left(8+k^4\right)}{\left(4+k^4\right)^2}}\right)\right)
\end{multline}
\end{example}
\begin{example}
From Eq. (11.201) in \cite{wheelon}
\begin{multline}
\frac{\pi }{4}=\sum _{n=1}^{\infty } \frac{1}{3\times 2^{3/4} \left(1+n+n^2\right)^3 \left(1+\frac{1}{\left(1+n+n^2\right)^2}\right)^{3/2}
   \left(1+\sqrt{1-\frac{1}{\left(1+n^2\right)^2 (2+n (2+n))^2}}\right)^{3/4}}\\
\left(\sqrt{2} \, _2F_1\left(\frac{3}{4},\frac{3}{4};\frac{7}{4};\frac{1}{2} \left(1-\sqrt{1-\frac{1}{\left(1+n^2\right)^2 (2+n
   (2+n))^2}}\right)\right)\right. \\ \left.
+6 \left(1+n^2\right) (2+n (2+n)) \sqrt{1+\sqrt{1-\frac{1}{\left(1+n^2\right)^2 (2+n (2+n))^2}}}\right. \\ \left.
 \, _2F_1\left(\frac{1}{4},\frac{1}{4};\frac{5}{4};\frac{1}{2}
   \left(1-\sqrt{1-\frac{1}{\left(1+n^2\right)^2 (2+n (2+n))^2}}\right)\right)\right)
\end{multline}
\end{example}
\begin{example}
From Eq. (42.2.3) in \cite{hansen}
\begin{multline}
\tan ^{-1}\left(n x^n\right)\\
=\sum _{k=0}^{n-1} \frac{2^{3/4} x^{3 k} (k (-1+x)+x)^3 }{6 \left(1+k (1+k)
   x^{1+2 k}\right)^3 \left(1+\frac{x^{2 k} (k (-1+x)+x)^2}{\left(1+k (1+k) x^{1+2 k}\right)^2}\right)^{3/2}
   \left(1+\sqrt{1-\frac{x^{4 k} (k (-1+x)+x)^4}{\left(1+k^2 x^{2 k}\right)^2 \left(1+(1+k)^2 x^{2+2
   k}\right)^2}}\right)^{3/4}}\\
\left(\,
   _2F_1\left(\frac{3}{4},\frac{3}{4};\frac{7}{4};\frac{1}{2} \left(1-\sqrt{1-\frac{x^{4 k} (k
   (-1+x)+x)^4}{\left(1+k^2 x^{2 k}\right)^2 \left(1+(1+k)^2 x^{2+2 k}\right)^2}}\right)\right)\right. \\ \left.
+6 \sqrt[4]{2} x^k (k
   (-1+x)+x) \left(x^{2 k} (k (-1+x)+x)^2+\left(1+k (1+k) x^{1+2 k}\right)^2\right)\right. \\ \left.
 \sqrt{1+\sqrt{1-\frac{x^{4 k} (k
   (-1+x)+x)^4}{\left(1+k^2 x^{2 k}\right)^2 \left(1+(1+k)^2 x^{2+2 k}\right)^2}}}\right. \\ \left.
 \,
   _2F_1\left(\frac{1}{4},\frac{1}{4};\frac{5}{4};\frac{1}{2} \left(1-\sqrt{1-\frac{x^{4 k} (k
   (-1+x)+x)^4}{\left(1+k^2 x^{2 k}\right)^2 \left(1+(1+k)^2 x^{2+2 k}\right)^2}}\right)\right)\right)
\end{multline}
\end{example}
\begin{example}
From Eq. (42.2.1) in \cite{hansen}
\begin{multline}
2^n \tan ^{-1}\left(2^{-n} x\right)-\tan ^{-1}(x)\\
=\sum _{k=0}^{n-1} \frac{2^{-1+k} }{3 \left(2^{2+3 k}+3\times 2^k x^2\right)^3
   \left(1+\frac{x^6}{\left(2^{2+3 k}+3\times 2^k x^2\right)^2}\right)^{3/2}
   \left(1+\sqrt{1-\frac{x^{12}}{\left(4^k+x^2\right)^2 \left(4^{1+k}+x^2\right)^4}}\right)^{3/4}}\\
\left(2^{3/4} x^9 \,
   _2F_1\left(\frac{3}{4},\frac{3}{4};\frac{7}{4};\frac{1}{2} \left(1-\sqrt{1-\frac{x^{12}}{\left(4^k+x^2\right)^2
   \left(4^{1+k}+x^2\right)^4}}\right)\right)\right. \\ \left.
+6 \sqrt[4]{2} x^3 \left(4^k+x^2\right) \left(4^{1+k}+x^2\right)^2
   \sqrt{1+\sqrt{1-\frac{x^{12}}{\left(4^k+x^2\right)^2 \left(4^{1+k}+x^2\right)^4}}}\right. \\ \left.
 \,
   _2F_1\left(\frac{1}{4},\frac{1}{4};\frac{5}{4};\frac{1}{2} \left(1-\sqrt{1-\frac{x^{12}}{\left(4^k+x^2\right)^2
   \left(4^{1+k}+x^2\right)^4}}\right)\right)\right)
\end{multline}
\end{example}
\begin{example}
From Eq. (42.1.18) in \cite{hansen}
\begin{multline}
4 n (-1)^{n-1} \tan ^{-1}(\exp (-\pi  n))\\
=\sum _{k=0}^{\infty } \frac{8 (-1)^k }{3 (1+2 k)^5 \left(1+\frac{64 n^4}{(1+2
   k)^4}\right)^{3/2} \left(1+\sqrt{1-\frac{4096 n^8}{(1+2 k)^8 \left(1+\frac{64 n^4}{(1+2
   k)^4}\right)^2}}\right)^{3/4}}\\
\left(32\times 2^{3/4} n^6 \,
   _2F_1\left(\frac{3}{4},\frac{3}{4};\frac{7}{4};\frac{1}{2} \left(1-\sqrt{1-\frac{4096 n^8}{(1+2 k)^8
   \left(1+\frac{64 n^4}{(1+2 k)^4}\right)^2}}\right)\right)\right. \\ \left.
+3 \sqrt[4]{2} n^2 \left((1+2 k)^4+64 n^4\right)
   \sqrt{1+\sqrt{1-\frac{4096 n^8}{(1+2 k)^8 \left(1+\frac{64 n^4}{(1+2 k)^4}\right)^2}}}\right. \\ \left. \,
   _2F_1\left(\frac{1}{4},\frac{1}{4};\frac{5}{4};\frac{1}{2} \left(1-\sqrt{1-\frac{4096 n^8}{(1+2 k)^8
   \left(1+\frac{64 n^4}{(1+2 k)^4}\right)^2}}\right)\right)\right)
\end{multline}
\end{example}
\begin{example}
From Eq. (42.1.18) in \cite{hansen}
\begin{multline}
4 n (-1)^{n-1} \tan ^{-1}(\exp (-\pi  n))\\
=\sum _{k=0}^{\infty } \frac{8 (-1)^k }{3 (1+2 k)^5 \left(1+\frac{64 n^4}{(1+2
   k)^4}\right)^{3/2} \left(1+\sqrt{1-\frac{4096 n^8}{(1+2 k)^8 \left(1+\frac{64 n^4}{(1+2
   k)^4}\right)^2}}\right)^{3/4}}\\
\left(32\times 2^{3/4} n^6 \,
   _2F_1\left(\frac{3}{4},\frac{3}{4};\frac{7}{4};\frac{1}{2} \left(1-\sqrt{1-\frac{4096 n^8}{(1+2 k)^8
   \left(1+\frac{64 n^4}{(1+2 k)^4}\right)^2}}\right)\right)\right. \\ \left.
+3 \sqrt[4]{2} n^2 \left((1+2 k)^4+64 n^4\right)
   \sqrt{1+\sqrt{1-\frac{4096 n^8}{(1+2 k)^8 \left(1+\frac{64 n^4}{(1+2 k)^4}\right)^2}}}\right. \\ \left. \,
   _2F_1\left(\frac{1}{4},\frac{1}{4};\frac{5}{4};\frac{1}{2} \left(1-\sqrt{1-\frac{4096 n^8}{(1+2 k)^8
   \left(1+\frac{64 n^4}{(1+2 k)^4}\right)^2}}\right)\right)\right)
\end{multline}
\end{example}
\begin{example}
From Eq. (42.3.2) in \cite{hansen}
\begin{multline}
\tan ^{-1}\left(x^n\right)-\frac{1}{4} \pi  \left(1+(-1)^n\right)\\
=\sum _{k=0}^{n-1} \frac{(-1)^k \left(x-\cos
   \left(\frac{\pi +2 k \pi }{2 n}\right)\right) \csc \left(\frac{\pi +2 k \pi }{2 n}\right)\left(1+\left(\cot
   \left(\frac{\pi +2 k \pi }{2 n}\right)-x \csc \left(\frac{\pi +2 k \pi }{2 n}\right)\right)^2\right)^{-3/2}}{3\times 2^{3/4} \left(1+\sqrt{1-\frac{\left(x-\cos \left(\frac{(1+2 k) \pi
   }{2 n}\right)\right)^4 \csc ^4\left(\frac{(1+2 k) \pi }{2 n}\right)}{\left(1+\left(x-\cos \left(\frac{(1+2 k) \pi
   }{2 n}\right)\right)^2 \csc ^2\left(\frac{(1+2 k) \pi }{2 n}\right)\right)^2}}\right)^{3/4} }\\
 \left(\sqrt{2} \left(\cot
   \left(\frac{\pi +2 k \pi }{2 n}\right)-x \csc \left(\frac{\pi +2 k \pi }{2 n}\right)\right)^2\right. \\ \left.
 \,_2F_1\left(\frac{3}{4},\frac{3}{4};\frac{7}{4};\frac{1}{8} \left(4-2 \sqrt{2} \sqrt{\frac{\left(3+4 x^2-8 x \cos
   \left(\frac{\pi +2 k \pi }{2 n}\right)+\cos \left(\frac{\pi +2 k \pi }{n}\right)\right) \sin ^2\left(\frac{\pi +2 k
   \pi }{2 n}\right)}{\left(1+x^2-2 x \cos \left(\frac{\pi +2 k \pi }{2 n}\right)\right)^2}}\right)\right)\right. \\ \left.+6
   \sqrt{1+\sqrt{1-\frac{\left(x-\cos \left(\frac{(1+2 k) \pi }{2 n}\right)\right)^4 \csc ^4\left(\frac{(1+2 k) \pi
   }{2 n}\right)}{\left(1+\left(x-\cos \left(\frac{(1+2 k) \pi }{2 n}\right)\right)^2 \csc ^2\left(\frac{(1+2 k) \pi
   }{2 n}\right)\right)^2}}}\right. \\ \left.
 \left(1+\left(\cot \left(\frac{\pi +2 k \pi }{2 n}\right)-x \csc \left(\frac{\pi +2 k \pi
   }{2 n}\right)\right)^2\right)\right. \\ \left.
 \, _2F_1\left(\frac{1}{4},\frac{1}{4};\frac{5}{4};\frac{1}{8} \left(4-2 \sqrt{2}
   \sqrt{\frac{\left(3+4 x^2-8 x \cos \left(\frac{\pi +2 k \pi }{2 n}\right)+\cos \left(\frac{\pi +2 k \pi
   }{n}\right)\right) \sin ^2\left(\frac{\pi +2 k \pi }{2 n}\right)}{\left(1+x^2-2 x \cos \left(\frac{\pi +2 k \pi }{2
   n}\right)\right)^2}}\right)\right)\right)
\end{multline}
\end{example}
\begin{example}
From Eq. (42.3.4) in \cite{hansen}
\begin{multline}
-\tan ^{-1}\left(\sinh \left(2 n \sinh ^{-1}(x)\right) \text{sech}\left((2 n+1) \sinh
   ^{-1}(x)\right)\right)\\
=\sum _{k=0}^{2 n-1} \frac{x }{3\times 2^{3/4} \left(1+\sqrt{1-\frac{4 x^4}{\left(1+2 x^2+\cos
   \left(\frac{2 (\pi +2 k \pi )}{1+4 n}\right)\right)^2}}\right)^{3/4} \left(1+x^2 \sec ^2\left(\frac{\pi +2 k \pi
   }{1+4 n}\right)\right)^{3/2}}\\
\left(\sqrt{2} x^2 \,
   _2F_1\left(\frac{3}{4},\frac{3}{4};\frac{7}{4};\frac{1}{2} \left(1-\sqrt{1-\frac{4 x^4}{\left(1+2 x^2+\cos
   \left(\frac{2 (\pi +2 k \pi )}{1+4 n}\right)\right)^2}}\right)\right)\right. \\ \left.
+3 \left(1+2 x^2+\cos \left(\frac{2 (1+2 k)
   \pi }{1+4 n}\right)\right) \sqrt{1+\sqrt{1-\frac{4 x^4}{\left(1+2 x^2+\cos \left(\frac{2 (\pi +2 k \pi )}{1+4
   n}\right)\right)^2}}}\right. \\ \left.
 \, _2F_1\left(\frac{1}{4},\frac{1}{4};\frac{5}{4};\frac{1}{2} \left(1-\sqrt{1-\frac{4
   x^4}{\left(1+2 x^2+\cos \left(\frac{2 (\pi +2 k \pi )}{1+4 n}\right)\right)^2}}\right)\right)\right) \sec
   ^3\left(\frac{\pi +2 k \pi }{1+4 n}\right)
\end{multline}
\end{example}
\begin{example}
From Eq. (42.2.3) in \cite{hansen}
\begin{multline}
\tan ^{-1}\left(\sinh \left(2 n \sinh ^{-1}(x)\right) \text{sech}\left((2 n+1) \sinh
   ^{-1}(x)\right)\right)\\
=\sum _{k=1}^{2 n} \frac{\sqrt{2} x^3}{3\times 2^{3/4} \left(1+\sqrt{1-\frac{4 x^4}{\left(1+2 x^2+\cos
   \left(\frac{4 k \pi }{1+4 n}\right)\right)^2}}\right)^{3/4} \left(1+x^2 \sec ^2\left(\frac{2 k \pi }{1+4
   n}\right)\right)^{3/2}}\\
\left(  \,
   _2F_1\left(\frac{3}{4},\frac{3}{4};\frac{7}{4};\frac{1}{2} \left(1-\sqrt{1-\frac{4 x^4}{\left(1+2 x^2+\cos
   \left(\frac{4 k \pi }{1+4 n}\right)\right)^2}}\right)\right) \sec ^3\left(\frac{2 k \pi }{1+4 n}\right)\right. \\ \left.
+6 x
   \sqrt{1+\sqrt{1-\frac{4 x^4}{\left(1+2 x^2+\cos \left(\frac{4 k \pi }{1+4 n}\right)\right)^2}}}\right. \\ \left.
 \,
   _2F_1\left(\frac{1}{4},\frac{1}{4};\frac{5}{4};\frac{1}{2} \left(1-\sqrt{1-\frac{4 x^4}{\left(1+2 x^2+\cos
   \left(\frac{4 k \pi }{1+4 n}\right)\right)^2}}\right)\right)\right. \\ \left.
 \sec \left(\frac{2 k \pi }{1+4 n}\right) \left(1+x^2
   \sec ^2\left(\frac{2 k \pi }{1+4 n}\right)\right)\right)
\end{multline}
\end{example}
\begin{example}
From Eq. (6.8.12.4) in \cite{brychkov}
\begin{multline}
\sum _{k=1}^{\infty } (-1)^k J_{\frac{k}{2}}(k z){}^2
=-\frac{1}{2}+\frac{1}{2 \sqrt{1-4 z^2}}-\frac{\sin ^{-1}(2
   z)}{\pi  \sqrt{1-4 z^2}}\\
=-\frac{1}{2}+\frac{1}{2 \sqrt{1-4 z^2}}-\frac{\frac{4\ 2^{3/4} z^3 \,
   _2F_1\left(\frac{3}{4},\frac{3}{4};\frac{7}{4};\frac{1}{2} \left(1-\sqrt{1-16 z^4}\right)\right)}{3
   \left(1+\sqrt{1-16 z^4}\right)^{3/4}}+\frac{2 \sqrt[4]{2} z \,
   _2F_1\left(\frac{1}{4},\frac{1}{4};\frac{5}{4};\frac{1}{2}-\frac{1}{2} \sqrt{1-16 z^4}\right)}{\sqrt[4]{1+\sqrt{1-16
   z^4}}}}{\pi  \sqrt{1-4 z^2}}
\end{multline}
\end{example}
%
%\begin{example}
%From Eq. () in \cite{}
%
%\end{example}
%
%\section{Discussion}
%%
%In this paper, we have presented a few methods for deriving a new Arctangent integral transform along with some interesting definite integrals similar to those published by Oberhettinger, using contour integration. The results presented were numerically verified for both real and imaginary and complex values of the parameters in the integrals using Mathematica by Wolfram.
%%

%
\end{document}